\documentclass[centertags,12pt]{amsart}
\usepackage{latexsym}
\usepackage{amssymb}

\numberwithin{equation}{section}

\makeatletter
\renewcommand{\subsection}{\@startsection
{subsection}{2}{0mm}{\baselineskip}{-0.25cm}{\normalfont\normalsize\rm}}
\makeatother

\newtheorem{theorem}{Theorem}[section]
\newtheorem{proposition}[theorem]{Proposition}
\newtheorem{lemma}[theorem]{Lemma}
\newtheorem{corollary}[theorem]{Corollary}
\newtheorem{scholium}[theorem]{Scholium}
{\theoremstyle{definition}
\newtheorem{example}[theorem]{Example}
\newtheorem*{conjecture*}{Conjecture}
\newtheorem{problem}[theorem]{Problem}
\newtheorem*{definition*}{Definition}}
\theoremstyle{remark}
\newtheorem{remark}[theorem]{Remark}
\newtheorem*{claim*}{Claim}

\def\C{\mathbf C}
\def\F{\mathbf F}
\def\K{\mathbf K}
\def\N{\mathbf N}
\def\Q{\mathbf Q}
\def\P{\mathbf P}
\def\R{\mathbf R}
\def\Z{\mathbf Z}
\def\cA{\mathcal A}
\def\cB{\mathcal B}
\def\cC{\mathcal C}
\def\cD{\mathcal D}
\def\cE{\mathcal E}
\def\cF{\mathcal F}
\def\cH{\mathcal H}
\def\cJ{\mathcal J}
\def\cK{\mathcal K}
\def\cL{\mathcal L}
\def\cM{\mathcal M}
\def\cO{\mathcal O}

\def\cR{\mathcal R}
\def\cS{\mathcal S}
\def\cV{\mathcal V}
\def\cW{\mathcal W}
\def\cX{\mathcal X}
\def\cY{\mathcal Y}
\def\cZ{\mathcal Z}
\def\l{\ell}
\def\fl{\mathbf F_\ell}
\def\fls{\mathbf F_{\ell^2}}
\def\fq{\mathbf F_q}
\def\fqi{\mathbf F_{q^i}}
\def\aut{{\rm Aut}}
\def\dim{{\rm dim}}
\def\div{{\rm div}}
\def\Div{{\rm Div}}
\def\deg{{\rm deg}}
\def\det{{\rm det}}
\def\sq{\sqrt{q}}
\def\fro{{\mathbf \Phi_q}}

\def\supp{{\rm Supp}}
\begin{document}

\author[F.~Torres]{Fernando Torres}\thanks{Author's address: 
IMECC-UNICAMP, Cx. P. 6065, Campinas, 13083-970-SP-Brazil}
\thanks{ftorres@ime.unicamp.br}

\title[St\"ohr-Voloch's approach to the Hasse-Weil bound and
applications]{The approach of St\"ohr-Voloch\\ to the Hasse-Weil bound\\
with applications\\ to optimal curves and plane arcs}

\address{}

\maketitle

{\bf Contens.}

\begin{enumerate}
\item Linear series on curves
  \subitem1.1. Terminology and notation
  \subitem1.2. Morphisms from linear series; Castelnuovo's genus bound
  \subitem1.3. Linear series from morphisms
  \subitem1.4. Relation between linear series and morphisms
  \subitem1.5. Hermitian invariants; Weierstrass semigroups I

\item Weierstrass point theory 
  \subitem2.1. Hasse derivatives
  \subitem2.2. Order sequence; Ramification divisor
  \subitem2.3. $\cD$-Weierstrass points
  \subitem2.4. $\cD$-osculating spaces
  \subitem2.5. Weierstrass points; Weierstrass semigroups II

\item Frobenius orders

\item Optimal curves
  \subitem4.1. A $\fq$-divisor from the Zeta Function
  \subitem4.2. The Hermitian curve
  \subitem4.3. The Suzuki curve

\item Plane arcs
  \subitem5.1. B. Segre's fundamental theorem: Odd case
  \subitem5.2. The work of Hirschfeld, Korchm\'aros and Voloch
   \end{enumerate}

\section*{Introduction}

The objective of this paper is to report applications of the approach of
St\"ohr-Voloch to the Hasse-Weil bound \cite{sv}, to the investigation of
the uniqueness of certain optimal curves, as well as to the search of
upper bounds for the second largest size that a complete plane arc (in a
projective plane of odd order) can have.

Let $\cX$ be a (projective, geometrically irreducible, non-singular
algebraic) curve of genus $g$ defined over a finite field $\fq$ of $q$
elements. Weil \cite{weil} showed that
   \begin{equation*}
|\#\cX(\fq)-(q+1)|\leq 2\sq g\, ,\tag{$*$}
   \end{equation*} 
being this bound sharp as Example \ref{ex4.0} 
here shows. Goppa \cite{goppa} constructed linear codes from
curves defined over $\fq$. These codes were used by Tsfasman, Vladut and
Zink \cite{tvz} to show that the Gilbert-Varshamov bound can
be improved whenever $q$ is a square and $q\geq 49$. This was an
unexpected result for coding theorist.

The length and the minimum distance of Goppa codes are related with the
number of $\fq$-rational points in the underlying curve. Then Goppa's
construction provided motivation and in fact reawakened the interest in
the study of rational points of curves which, despite of this motivation,
is an interesting mathematical problem by its own.

Serre \cite{serre} noticed that $(*)$ can be improved by replacing
$2\sq$ by $\lfloor 2\sq\rfloor$. A refined version of Ihara \cite{ihara}
shows that
   $$
g>\frac{q^2-q}{2\sq^2+2\sq-2q}\quad\Rightarrow\quad 
\#\cX(\fq)<q+1+\lfloor2\sq\rfloor g\, ,
   $$
and in this case Serre \cite{serre}, \cite{serre2} upper bounded
$\#\cX(\fq)$ via explicit formulae.

A geometric point of view to bound $\#\cX(\fq)$ was introduced by St\"ohr
and Voloch \cite{sv}: Suppose that $\cX$ admits a base-point-free
linear series $g^r_d$ defined over $\fq$; then
   $$
\#\cX(\fq)\leq\frac{\sum_{i=0}^{r-1}\nu_i(2g-2)+(q+r)d}{r}\, ,
   $$ 
where $\nu_0,\ldots,\nu_{r-1}$ are certain $\fq$-invariants associated to
$g^r_d$ (see Theorem \ref{thm3.1} here). By an appropriate choice of
$g^r_d$ this result implies $(*)$ \cite[Cor. 2.14]{sv}, and in several
cases one obtains improvements on $(*)$. We write an exposition of
St\"ohr-Voloch's approach in Sect. \ref{s3}. For the sake of completeness
we include an expository account on Weierstrass point theory of linear series
on curves: Sects. \ref{s1}, \ref{s2}.

Next we discuss two applications of \cite{sv} studied here. The first one
is concerning the uniqueness of certain optimal curves. The most well
known example of a $\fq$-maximal curve is the {\em Hermitian curve}
(Example \ref{ex4.0} here) whose genus is $\sq(\sq-1)/2$; i.e., the
biggest one that a $\fq$-maximal curve can have according to the
aforementioned Ihara's result. R\"uck and Stichtenoth \cite{r-sti} showed
that this property characterize Hermitian curves up to $\fq$-isomorphic.
In Sect. \ref{s4.1} we equip the curve $\cX$ with a linear series
$\cD_\cX$ obtained from its Zeta Function provided that
$\cX(\fq)\neq\emptyset$. It turns out that $\cD_\cX=|(\sq+1)P_0|$,
$P_0\in\cX(\fq)$, whenever $\cX$ is $\fq$-maximal. Then applying \cite{sv}
to $\cD_\cX$ we prove a stronger version of R\"uck-Stichtenoth's result;
see Theorem \ref{thm4.21} here. Further properties of $\fq$-maximal were
proved via an interplay of St\"ohr-Voloch's paper \cite{sv}, and results
on linear series such as Castelnuovo's genus bound and Halphen's theorem
applied to $\cD_\cX$; see \cite{fgt}, \cite{ft2},\cite{kt1},\cite{kt2}. A
characterization result is also proved for the Suzuki curve (Theorem
\ref{thm4.31}), which in fact is optimal with genus $q_0(q-1)$ and
$(q^2+1)$ $\fq$-rational points.

The second application of \cite{sv} studied here is the bounding of the
size $k$ of a complete plane arc $\cK$ in $\P^2(\fq)$ which indeed is a
basic problem in Finite Geometry. What it makes this possible is the fact
that associated to $\cK$ there is a (possible singular) plane curve $\cC$.
A fundamental result of B. Segre \cite{segre} (see Theorem \ref{thm5.11}
here for the odd case) allows then to upper bound $k$ via \cite{sv}
applied to certain linear series defined on the non-singular model of an
irreducible component of $\cC$. Details of the following discussion can be
seen in Sect. \ref{s5}. The largest $k$ is already well known and so the
problem is concerning the second largest size $m_2'(2,q)$. Let $q$ be a
square. If $q$ is even, then $m'(2,q)=q-\sq+1$ and a similar result is
expected for $q$ odd, $q\geq 49$. Let $q$ be odd. Applying $(*)$ B. Segre
showed that $m'(2,q)\leq q-\sq/4+7/4$. One obtains the same bound by using
\cite{sv};  see Proposition \ref{prop5.21} here. If in addition, for $q$
large, one takes into consideration a bound for the number of
$\fq$-rational of plane
curves due to Hirschfeld and Korchm\'aros \cite{kt2} (see Theorem
\ref{thm5.21} here) one finds the currently best upper bound for
$m'(2,q)$, namely
   $$
m'(2,q)\leq q-\frac{\sq}{2}+\frac{5}{2}\, .
   $$ 
So far, for $\sq\not\in\N$, the best upper bound for
$m'(2,q)$ is due to Voloch \cite{voloch1}, \cite{voloch2}; see Lemmas
\ref{lemma5.23}, \ref{lemma5.24} here.

This paper is an outgrowth and a considerable expanded of lectures given
at the University of Essen in April 1997 and the University of Perugia in
February 1998.

{\bf Convention.} The word {\em curve} will mean a projective,
irreducible, non-singular algebraic curve. 

\section{Linear series on curves}\label{s1} 

The purpose of this section is to summarize relevant material regarding
linear series on curves. Standard references are
Arbarello-Cornalba-Griffiths-Harris \cite{acgh}, Griffiths
\cite{griffiths}, Griffiths-Harris \cite{griffiths-harris}, Hartshorne 
\cite{hartshorne}, Namba \cite{namba}, Seidenberg \cite{seidenberg},
Stichtenoth \cite{sti}. 

Let $\cX$ be a curve over an algebraically closed field $\F$; set
$\P^r:=\P^r(\F)$.

    \subsection{Terminology and notation}\label{s1.1} We start by fixing
some terminology and notation.

1.1.1.\ We denote by $\Div(\cX)$ the group of {\it divisors} on $\cX$; 
i.e., the $\Z$-free abelian group generated by the points of $\cX$. Let
$D=\sum n_PP \in \Div(\cX)$.  The {\em multiplicity} of $D$ at $P$ is
$v_P(D):=n_P$. The divisor $D$ is called {\it effective} (notation:
$D\succeq 0$) if $v_P(D)\ge 0$ for each $P$. For $D, E\in \Div(\cX)$, we
write $D\succeq E$ if $D-E\succeq 0$. The {\it degree} of $D$ is the
number $\deg(D):= \sum v_P(D)$, and the {\it support} of $D$ is the set
$\supp(D):=\{P\in X:  v_P(D)\neq 0\}$. 

1.1.2.\ Let $\F(\cX)$ denote the field of rational functions on $\cX$. 
Associated to $f\in \F(\cX)^*:=\F(\cX)\setminus\{0\}$ we have the divisor
  $$ 
\div(f):= \sum v_P(f)P\, ,
  $$
where $v_P$ stands for the {\em valuation} at $P\in \cX$. Recall that 
$v_P$ satisfies: $v_P(0):=+\infty$, $v_P(f+g)\ge \min(v_P(f),v_P(g))$, and
$v_P(fg)=v_P(f)+v_P(g)$ for $f,g\in \F(\cX)$. 

For $f\in \F^*:=\F\setminus\{0\}$, $\div(f)=0$ and for $f\in
\F(\cX)\setminus \F$, $\div(f)=\div_0(f)-\div_\infty(f)$, where
$\div_0(f):=\sum_{v_P(f)>0}v_P(f)P$ and
$\div_\infty(f):=\sum_{v_P(f)<0}(-v_P(f))P$ are respectively the {\it
zero} and the {\it polar} divisor of $f$. Moreover, $\deg(\div(f))=0$ and
$\div(fg)=\div(f)+\div(g)$. 

Associated to $D\in \Div(\cX)$ we have the $\F$-linear space  
  $$
L(D):=\{f\in \F(\cX)^*: D+\div(f)\succeq 0\}\cup\{0\}\, ,
  $$ 
where $\ell(D):=\dim_\F L(D)\le \deg(D)+1$. For $D, E\in \Div(\cX)$
such that $L(D)\subseteq L(E)$, we have
  $$
\ell(E)-\ell(D)\leq \deg(E)-\deg(D)\, .
  $$
The Riemann-Roch theorem computes $\ell(D)$: If $C$ is a canonical divisor
on $\cX$ and $g$ is the genus of $\cX$, then
  $$
\ell(D)=\deg(D)+1-g+\ell(C-D)\, .
  $$
In particular, $C$ is characterized by the properties: $\deg(C)=2g-2$ and
$\ell(C)\geq g$.

A {\em local parameter} at $P\in\cX$ is a rational function $t\in\F(\cX)$
such that $v_P(t)=1$. Associated to $f\in\F(\cX)^*$ we have its {\em local
expansion at $P$}, $\sum_{i=v_P(f)}^\infty a_it^i$, where
$a_{v_P(f)}\neq 0$. Let $f\in\F(\cX)$ be a {\em separating variable of
$\F(\cX)|\F$}; i.e., let the field extension $\F(\cX)|\F(f)$ be separable.
Then we have the divisor of the {\em differential of $f$}, namely
$\div(df)$ where $v_P(\div(df))$ equals the minimum integer
$i$ such that $ia_i\neq 0$. It holds that $\deg(\div(f))=2g-2$.

1.1.3.\ Two divisors $D, E\in \Div(\cX)$ are called {\it linearly
equivalent} (notation: $D\sim E$) if there exists $f\in \F(\cX)^*$ such
that
$D=E+\div(f)$. In this case, $\deg(D)=\deg(E)$ and $L(D)$ is
$\F$-isomorphic to $L(E)$ via the map $g\mapsto fg$. For $E\in \Div(\cX)$,
let
  $$
|E|:= \{D\in \Div(\cX): D\succeq 0,\ D\sim E\}\, ;
  $$
i.e.,
  $$
|E|=\{E+\div(f): f\in L(E)\setminus\{0\}\}\, .
  $$
Since, for $f,g \in \F(X)^*$, $\div(f)=\div(g)$ if and only if there
exists $a \in \F^*$ such that $f=ag$, the set $|E|$ is equipped with a
structure of projective space by means of the map $E+\div(f)\in |E|\mapsto
[f]\in \P(L(E))$; notation: $|E|\cong \P(L(E))$.

A {\em linear series} $\cD$ on $\cX$ is a subset of some $|E|$, of type
  $$
\{E+\div(f): f\in \cD'\setminus\{0\}\}\, ,
  $$
with $\cD'$ being a $\F$-linear subspace of $L(E)$. The numbers
$d=\deg(\cD):=\deg(E)$ and $r=\dim(\cD):=\dim_\F(\cD')-1$ are called
respectively the {\em degree} and the (projective) {\em dimension} of
$\cD$. We say
that $\cD$ is a $g^r_d$ on $\cX$. $\cD$ is called {\em complete} if
$\cD=|E|$. Observe that, under the
identification $|E|\cong \P(L(E))$, $\cD$ corresponds to $\P(\cD')$;
notation: $\cD\cong \P(\cD')\subseteq |E|$. A linear series $\cD_1\cong
\P(\cD_1')\subseteq |E_1|$ will be called a {\em subspace} of
$\cD\cong
\P(\cD')\subseteq |E|$ if $L(E_1)\subseteq L(E)$ and $\cD_1'\subseteq
\cD'$.

1.1.4.\ Let $P\in\cX$ and $f\in \F(\cX)$ {\em regular} at $P$; i.e.,
$v_P(f)\ge 0$. Then there exists a unique $a_f\in \F$ such that
$v_P(f-a_f)>0$. We set $f(P):=a_f$. For $f, g\in \F(\cX)$ regular at $P$,
$(f+g)(P)=f(P)+g(P)$ and $(fg)(P)=f(P)g(P)$. A point of the $r$-projective
space $\P^r$ will be denoted by $(a_0:\dots:a_r)$.

Let $\phi: \cX\to \P^r$ be a morphism; i.e., let $f_0,\ldots,f_r\in
\F(\cX)$, not all zero, such that
   $$
\phi(P)=((t^{e_P}f_0)(P):\ldots:(t^{e_P}f_r(P))\, ,
   $$
where $t$ is a local parameter at $P$, and 
   $$
e_P:=-{\rm min}\{v_P(f_0),\ldots,v_P(f_r)\}\, .
   $$ 
Observe that each $t^{e_P}f_i$ is regular at $P$. The
rational functions $f_0, \ldots,f_r$ are called (homogeneous) {\em 
coordinates} of $\phi$. We set
   $$ 
\phi=(f_0:\ldots:f_r)\, .
   $$
The coordinates $f_0,\ldots,f_r$ are uniquely determinated by $\phi$ up to
a factor in $\F(\cX)^*$; so $\phi$ corresponds to a point of
$\P^r(\F(\cX))$. If $\phi$ is non-constant, the image $\phi(\cX)$ is a
(possible singular)  algebraic curve in $\P^r$ whose function field is
$\F(\phi(\cX))=\F(f_0,\ldots,f_r)$. The curve $\cX$ can be thought as a
parametrized curve in $\P^r$, or $\phi(\cX)$ as being a concrete
manifestation of $\cX$ in $\P^r$. For $Q\in \phi(\cX)$, the points of the
fiber $\phi^{-1}(Q)$ will be called the {\em branches} of $\phi(\cX)$
centered at $Q$. The {\em degree} of $\phi$ is $\deg(\phi):=
[\F(\cX):\F(\phi(\cX))]$. 

   \begin{example}\label{ex1.11} Each rational function $f\in\F(\cX)$ can
be seen as a morphism $f:\cX\to \P^1=\F\cup\{\infty\}$, such that
$P\mapsto f(P)$ if
$P\not\in\div_\infty(f)$; $P\mapsto \infty$ otherwise. If $f\not\in\F$, we
have $d:=\deg(f)=[\F(\cX):\F(f)]=\deg(\div_\infty(f))$. Moreover, if
$\F(\cX)|\F(f)$ is separable, the genus $g$ of $\cX$ can be computed via
the so-called Riemann-Hurwitz formula:
  $$
2g-2=d(-2)+\deg(R_f)\, ,
  $$ 
where $R_f=\div(df)+2\div_\infty(f)$ is the ramification divisor of $f$.
If ${\rm char}(\F)$ does not divide the ramification index
$e_P$ of $P$ over $f(P)$, then $v_P(R_f)=e_P-1$
otherwise $v_P(R_f)>e_P-1$. We have the product formula 
  $$
\sum_{P\in f^{-1}(f(P))} e_P=d\, .
  $$
   \end{example}

For all but finitely many $Q\in\phi(\cX)$, $\#\phi^{-1}(Q)$ equals the
separable degree of $\F(\cX)|\F(\phi(\cX))$.  $\phi$ is called {\em
birational} (resp. {\em embedding}) if $\deg(\phi)=1$ (resp. $\cX$ is
$\F$-isomorphic to $\phi(\cX)$); in both cases, $\cX$ is a (the)  
non-singular model of $\phi(\cX)$.

Let $H$ be a hyperplane in $\P^r$ such that $\phi(\cX)\not\subseteq H$.
Then $\#\phi(\cX)\cap H$ is finite. To each $P\in \cX$ one associates a
number $I_P(H)=I(\phi(\cX),H;P)$, called the {\em intersection
multiplicity} of $\phi(\cX)$ and $H$ at $P$, in such a way that $I_P=0
\Leftrightarrow P\not\in \phi(\cX)\cap H$ and that $\sum I_P(H)$ is
independent of $H$; i.e., if $H'$ is another hyperplane in $\P^r$ such
that $\phi(\cX)\not\subseteq H'$, then $\sum I_P(H)=\sum I_P(H')$. This
number is called the {\em degree} $\deg(\phi(\cX))$ of $\phi(\cX)$. If
$\phi(\cX)\subseteq \P^2$, the degree of $\phi(\cX)$ equals the degree of
the polynomial that defines $\phi(\cX)$.

A morphism $\phi:\cX\to \P^r$ is called {\em non-degenerate} if
$\phi(\cX)\not\subseteq H$ for each hyperplane $H$ in $\P^r$. A curve
$\cX\subseteq \P^r$ is called {\em non-degenerate} if the inclusion
morphism $\cX\hookrightarrow\P^r$ is so.

     \begin{lemma}\label{lemma1.11} A morphism
$\phi=(f_0:\ldots:f_r):\cX\to \P^r$ is non-degenerate if and only if
$f_0,\ldots,f_r$ are $\F$-linearly independent. 
     \end{lemma}
     \begin{proof} There exists a hyperplane $H$ in $\P^r$ such that
$\phi(\cX)\subseteq H$ if and only if there exist $a_0,\ldots, a_r\in\F$,
not all zero, such that $\sum_ia_if_i(P)=0$ for all but finitely many
$P\in \cX$. The last condition is equivalent to $\sum_ia_if_i=0$, as a
non-zero rational function has only finitely many zeros (cf. Sect. 1.1.2);
now the result follows.
      \end{proof}

For $V\subseteq \F(\cX)$, $\langle V\rangle$ stands for the $\F$-vector 
space in $\F(\cX)$ generated by $V$. 

    \subsection{Morphisms from linear series; Castelnuovo's
genus bound}\label{s1.2} Let $\cD$ be a
$r$-dimensional linear series on $\cX$, say $\cD\cong \P(\cD')\subseteq
|E|$. The following subsets will provide information on the geometry of
$\cX$.

     \begin{definition*} For $P\in \cX$ and $i\in \N_0$,
$$
\cD_i(P):=\{D\in \cD: D\succeq iP\}\, .
$$
     \end{definition*}

Clearly $\cD_i(P)\supseteq \cD_{i+1}(P)$ and $\cD_i(P)=\emptyset$ if 
$i>d$. 

     \begin{lemma}\label{lemma1.21}
\begin{enumerate}
   \item[\rm(1)] $\cD_i(P)$ is a linear series$;$
   \item[\rm(2)] $\cD_i(P)$ is a subspace of $\cD;$
   \item[\rm(3)] $\dim(\cD_i(P))\le \dim(\cD_{i+1}(P))+1.$
\end{enumerate}
     \end{lemma}
     \begin{proof} Set $\cD_j:=\cD_j(P)$ and let $f\in
\cD'\setminus\{0\}$. Then $E+\div(f)\in \cD_i$ if and only if
$v_P(E)+v_P(f)\ge i$;  i.e., $\cD_i\cong \P(\cD_i')$, where
   $$
\cD_i':= \cD'\cap L(E-iP)\, .
   $$
This shows parts (1) and (2). Now $\cD_i^{'}/\cD_{i+1}^{'}$ 
is $\F$-isomorphic to a $\F$-subspace of
$\cL:=L(E-iP)/L(E-(i+1)P)$. Since $\dim_\F\cL\le 1$ (see Sect. 1.1.2),
part (3) follows. 
     \end{proof} 
     \begin{definition*} The {\em multiplicity} of $\cD$ at $P\in\cX$ is
defined by
   $$
b(P):=\ {\rm min}\{v_P(D): D\in \cD\}\, .
   $$
     \end{definition*} 
We have $b(P)>0$ if and only if $P\in \supp(D)$ for all $D\in \cD$;
so $b(P)\neq 0$ for finitely many $P\in X$. Consequently, we can define
the effective divisor $B=B^{\cD}$ on $\cX$ by setting
  $$
v_P(B):= b(P)\, .
  $$
     \begin{definition*} The divisor $B$ is called the {\em base
locus of $\cD$}. A point $P\in \supp(B)$ is called a {\em base point
of $\cD$}. If $B=0$, $\cD$ is called {\em base-point-free}.
     \end{definition*}

Thus $\cD$ is base-point-free if and only if for each $P\in \cX$ there
exists $f\in\cD'\setminus\{0\}$ such that $v_P(E+\div(f))=0$. Now, since
$D\succeq B$ for each $D\in \cD$, $\cD'\subseteq L(E-B)$ and
  $$
\cD^B:=\{D-B: D\in \cD\} \subseteq |E-B|
  $$ 
is a subspace of $\cD$ such that $\cD^B\cong\P(\cD')\subseteq |E-B|$. We
have $B^{\cD^B}=0$; i.e., $\cD^B$ is a $g^r_{d-\deg(B)}$ base-point-free
on $\cX$.

     \begin{lemma}\label{lemma1.22} 
Let $\cD\cong \P(\cD')\subseteq |E|$ be a linear series, where 
$\cD'=\langle f_0,\ldots,f_s\rangle$. Then $E$ is
determinated by $\cD$; i.e,
    $$
v_P(E)= b(P)-{\rm min}\{v_P(f_0),\ldots,v_P(f_s)\}\, .
    $$
     \end{lemma}
     \begin{proof} Since $\cD'\subseteq L(E-B)$, $v_P(E)-b(P)+v_P(f_i)\geq
0$ for each $i$ and each $P$ so that $v_P(E)\ge b(P) -{\rm
min}\{v_P(f_0),\ldots,v_P(f_s)\}$. On other hand, as $\cD^B$ is
base-point-free, for each $P$ there exists $(a_0:\ldots:a_s)\in
\P^s(\F)$ such that $v_P(E-B+\div(\sum_ia_if_i))=0$; now the result
follows.
     \end{proof}

Next we associate a morphism to $\cD$. For $P\in\cX$ we have  
$\cD=\cD_{b(P)}(P)\supsetneqq\cD_{b(P)+1}(P)$, so that 
$\dim(\cD_{b(P)+1})=\dim(\cD)-1$ by Lemma \ref{lemma1.21}. Thus we have
the following map
   $$
\phi_{\cD}:\ \cX\to \cD^* \cong \P(\cD')^* \, ,\ \  \qquad
P\mapsto \cD_{b(P)+1}\, .
   $$
Homogeneous coordinates of $\phi_\cD$ are given as follows. 
Let $\{f_0,\ldots,f_r\}$ be a $\F$-base of $\cD'$, $t$ a local parameter
at $P$, and $f\in \cD'\setminus \{0\}$. Then 
$v_P(t^{v_P(E)-b(P)}f)\ge 0$ and
    $$
E+\div(f)\in \cD_{b(P)+1}\ \Leftrightarrow\ v_P(t^{v_P(E)-b(P)}f)\ge 1
\ \Leftrightarrow\ (t^{v_P(E)-b(P)}f)(P)=0\, .
    $$
Since $f=\sum_i a_if_i$ with $(a_0:\ldots:a_r)\in \P^r$, we have 
      \begin{align*}   
\cD_{b(P)+1}& \cong \{(a_0:\ldots:a_r)\in \P^r:
\sum_{i=0}^r(t^{v_P(E)-b(P)}f_i)(P)a_i=0\}\in {\P^r}^*\\
            & \cong
((t^{v_P(E)-b(P)}f_0)(P):\ldots:(t^{v_P(E)-b(P)}f_r)(P))\in \P^r\, .
      \end{align*}
Hence from Lemma \ref{lemma1.22} the morphism $\phi_{f_0,\ldots,f_r}:=
(f_0:\ldots:f_r)$ gives a coordinate description of $\phi_\cD$, and it
will be referred as {\em a morphism associated to $\cD$}. If
$\phi_{g_0,\ldots,g_r}$ is another morphism associated to $\cD$, then
$\phi_{g_0,\ldots,g_r}=T\circ\phi_{f_0,\ldots,f_r}$, with $T\in
\aut(\P^r(\F))$; i.e., a morphism associated to $\cD$ is uniquely
determinated by $\cD$, up to projective equivalence. Observe that
$\phi_\cD$ and $\phi_{\cD^B}$ have the same coordinate description. We
summarize the above discussion as follows. 

       \begin{lemma}\label{lemma1.23} Let $\cD\cong \P(\cD')$ be a
$r$-dimensional linear series on $\cX$. Each $\F$-base $f_0,\ldots,f_r$ of
$\cD'$ defines a non-degenerate morphism
$\phi_{f_0,\ldots,f_r}=(f_0:\ldots:f_r):\cX\to \P^r$. If $g_0,\ldots,g_r$
is another $\F$-base of $\cD'$, then there exists $T\in \aut(\P^r)$ such
that $\phi_{g_0,\ldots,g_r}=T\circ\phi_{f_0,\ldots,f_r}$. 
       \end{lemma}

At this point we recall Castelnuovo's genus bound. Let $g$ be the genus of
$\cX$.

     \begin{definition*} A linear series $\cD$ is called {\em simple} if a
(any) morphism associated to $\cD$ is birational.
     \end{definition*}
Let $\cD$ be a simple $g^r_d$, $r\ge 2$, on $\cX$. Let
$d':=d-\deg(B^\cD)$, and let $\epsilon$ be the unique integer with $0\leq
\epsilon\leq r-2$ and $d'-1\equiv \epsilon\pmod{(r-1)}$. Define
Castelnuovo's number $c_0(d',r)$ by
   $$
c_0(d',r)=\frac{d'-1-\epsilon}{2(r-1)}(d'-r+\epsilon)\, .
   $$

   \begin{lemma}\label{lemma1.24} {\rm (Castelnuovo's genus bound for
curves in projective spaces, \cite{castelnuovo}, \cite[p. 116]{acgh},
\cite[IV, Thm. 6.4]{hartshorne}, \cite[Cor. 2.8]{rathmann})}
   $$
g\le c_0(d',r)\, .
   $$
   \end{lemma}

   \begin{remark}\label{rem1.21}
   $$
c_0(d',r)\le \begin{cases}
(d'-1-(r-1)/2)^2/2(r-1) & \text{for $r$ odd,}\\
(d'-1-(r-1)/2)^2-1/4)2/(r-1) & \text{for $r$ even.}\end{cases}
   $$
   \end{remark}
   \begin{remark}\label{rem1.22} Any curve $\cX$ of genus $g$ admits a
simple $g^2_d$ (i.e., a birational plane model) such that
   $$
g=d(d-1)/2-\sum_{P}\delta_P\, ,
   $$
where the $\delta_P$'s are the $\delta$-invariants of the plane curve
$\phi(\cX)$ with $\phi$ being a morphism associated to $g^2_d$. We have
that $\delta_P>0$ if and only if $\phi(\cX)$ is singular at $P$. A nice method
to compute $\delta_P$ was recently noticed by Beelen and Pellikaan
\cite{beelen-pellikaan}.
   \end{remark}
    \subsection{Linear series from morphisms}\label{s1.3} Let
$\phi=(f_0:\ldots:f_r):\cX\to \P^r$ be a morphism on $\cX$. In Sect.
1.1.4 we defined
   $$
e_P=-{\rm min}\{v_P(f_0),\ldots,v_P(f_r)\}\, ,\qquad P\in\cX\, .
   $$ 
Then $e_P\neq 0$ for finitely many $P\in\cX$, and so we have a divisor
$E=E_{f_0,\ldots,f_r}$ defined by 
   $$
v_P(E):=e_P\, .
   $$ 
Observe that $f_i\in L(E)$ for each $i$. Let
   $$
\cD':= \langle f_0,\ldots,f_r\rangle\subseteq L(E)\, .
   $$ 
Then we have the following linear series on $\cX$
   $$
\cD_{f_0,\ldots,f_r}:= \{E+\div(f): f\in \cD'\setminus\{0\}\}\subseteq
|E|\, ,
   $$
which is base-point-free. Indeed, $v_P(E+\div(f_{i_0}))=0$ where $i_0$ is
defined by $e_P=-v_P(f_{i_0})$. In addition, 
if $\phi_1=(g_0:\ldots:g_r)=T\circ\phi$ with $T\in \aut(\P^r)$, then
   $$
\min\{v_P(g_0),\ldots,v_P(g_r)\}=\min\{v_P(f_0),\ldots,v_P(f_r)\}\, ,
   $$
and hence $\cD_{g_0,\ldots,g_r}=\cD_{f_0,\ldots,f_r}$. Moreover, if 
$h\in \F(\cX)^*$, then
   \begin{align*}
E_{f_0h,\ldots,f_rh}  & =E_{f_0,\ldots,f_r}-\div(h)\\
\intertext{and so}
\cD_{f_0h,\ldots,f_rh} & =\cD_{f_0,\ldots,f_r}\, .
   \end{align*}
Consequently, the linear series $\cD_\phi:=\cD_{f_0,\ldots,f_r}$ is
uniquely determinated by $\phi$ and it is invariant under projective
equivalence of morphisms. Summarizing we have the following.

    \begin{lemma}\label{lemma1.31}
Associated to a morphism $\phi=(f_0:\ldots:f_r):X\to \P^r$, there
exists a base-point-free linear series 
$\cD_{\phi}\subseteq |E|,$ where $E$ is defined by
   $$
v_P(E):=-\min\{v_P(f_0),\ldots,v_P(f_r)\}\, .
   $$
If $\phi$ is non-degenerate, then $\dim(\cD_\phi)=r$. If
$\phi_1=T\circ\phi$, $T\in \aut(\P^r)$, then $\cD_{\phi_1}=\cD_\phi$.
    \end{lemma}
In the remaining part of this subsection, we let $\phi=(f_0:\ldots:f_r)$ 
be a non-degenerate morphism on $\cX$. Then $\cD_\phi$ is given by
   $$
\cD_\phi=\{E+\div(\sum_{i=0}^{r}a_if_i): (a_0:\ldots:a_r)\in
\P^r\}\, ,
   $$
because $\sum_ia_if_i=0 \Leftrightarrow a_i=0$ for each $i$ by Lemma
\ref{lemma1.11}. Therefore, since 
the point $(a_0:\ldots:a_r)$ can be identify with the hyperplane $H$ of  
equation $\sum_ia_iX_i=0$, 
      \begin{equation}\label{eq1.1}
\cD_\phi=\{\phi^*(H): H\ {\rm hyperplane\ in}\ \P^r\}\, ,
      \end{equation}
where $\phi^*(H) = E+\div(\sum_ia_if_i)$ is the pull-back of $H$ by
$\phi$.

   \begin{lemma}\label{lemma1.311} We have
$\phi^*(H)=(T\circ\phi)^*(T(H))$, where $T\in \aut(\P^r)$ and $H$ is a
hyperplane in $\P^r$.
   \end{lemma}
   \begin{proof} The result follows from the facts that
$E_\phi=E_{T\circ\phi}$ and that $T(H):\sum_i b_iY_i=0$, where
$(b_0,\ldots,b_r)=(a_0,\ldots,a_r)A^{-1}$, $A$ being the matrix 
defining $T$ and $H: \sum_ia_iX_i=0$.
   \end{proof}

   \begin{lemma}\label{lemma1.32} With the aforementioned notation,
    \begin{enumerate}
  \item[\rm(1)] $P\in \supp(\phi^* (H))\Leftrightarrow
\phi(P)\in H$; i.e, $\supp(\phi^*(H))=\phi^{-1}(\phi(\cX)\cap H);$

\item[\rm(2)] For $P_1\in \phi^{-1}(\phi(P))$,
$P_1\in\supp(\phi^*(H))\Leftrightarrow \phi^{-1}(\phi(P))\subseteq
\supp(\phi^*(H));$

\item[\rm(3)] $d:=\deg(\cD)=\deg(\phi)\deg(\phi(\cX)).$
    \end{enumerate}
     \end{lemma}
     \begin{proof} Let $t$ be a local parameter at $P\in \cX$. 

(1) The proof follows from the equivalences 
   $$
P\in \supp(\phi^*(H)) \Leftrightarrow v_P(\div(\sum_ia_it^{e_P}f_i))\ge 1
\Leftrightarrow (\sum_ia_it^{e_P}f_i)(P)=0\, .
   $$
(2) The implication ($\Leftarrow$) is trivial. ($\Rightarrow$): Let
$P_2\in \phi^{-1}(\phi(P))$. Then $\phi(P_1)=\phi(P_2)$ which belong to
$H$ by part (1). Thus, once again by (1) we conclude that $P_2\in
\supp(\phi^*(H))$. 

(3) Let $H_1$ be a hyperplane in $\P^r$ such that $\phi(\cX)\cap H\cap
H_1=\emptyset$. Denote by $h/h_1$ the rational function on $\P^r$, 
obtained by dividing the equation of $H$ by the one of $H_1$. Then 
we obtain a rational function on $\cX$, namely $\varphi:=(h/h_1)\circ
\phi$
(i.e., the pull-back of $h/h_1$ by $\phi$). The function $h/h_1$ is
regular on $\P^r\setminus H_1$ and hence $\varphi$ is regular on 
$\phi^{-1}(\P^r\setminus H_1)$. Moreover, by the election of $H_1$, we
have that $v_P(\varphi)\ge 1 \Leftrightarrow \phi(P)\in H$ and therefore 
from part (1) we conclude that $v_P(\varphi)\ge 1\Leftrightarrow P\in
\supp(\phi^*(H))$. From the definition of $\varphi$ we even conclude that
$\phi^*(H)=\div_0(\varphi)$.

Now suppose that $\phi(P)=Q\in \phi(\cX)\cap H$ is non-singular; let $u$ 
be a local parameter at $Q$ and set $i_P:=v_P(u)$ (the ramification
index at $P$). By considering $h/h_1$ as a function on $\phi(\cX)$ we have    
$v_P(\phi^{-1}(H))=v_P(\varphi)=i_Pv_Q(h/h_1)$, and by the 
product formula we also have
   $$
\sum_{P\in \phi^{-1}(Q)}v_P(\phi^{-1}(H))=\deg(\phi)v_Q(h/h_1)\, .
   $$
Now take $H$ such that every point in $\phi(X)\cap H$ is non-singular
(this is possible because $\phi(\cX)$ has a finite number of singular
points and so we can apply Bertini's theorem). Then from the above
equation, 
   $$
d=\deg(\phi)\sum_{Q\in \phi(\cX)\cap H}v_Q(h/h_1)\, .
   $$ 
It turns out that $v_Q(h/h_1)=I(\phi(\cX),H;Q)$ (cf. 
\cite[Ex.6.2]{hartshorne}), and the result follows.
     \end{proof}
From this lemma and its proof we obtain:

    \begin{corollary}\label{cor1.31} Let $\phi:\cX\to \P^r$ be a
non-degenerate morphism.  
    \begin{enumerate} 
    \item[\rm(1)] If $\phi$ is birational; i.e., $\deg(\phi)=1$, then
$\deg(\cD_\phi)=\deg(\phi(\cX))$.

    \item[\rm(2)] If $\cX\subseteq \P^r$ and $\phi$ is the inclusion
morphism, then 
   $$
\cD_\phi=\{\cX\cdot H: H\ {\rm hyperplane\ in}\ \P^r\}\, ,
   $$
where $\cX\cdot H=\sum_P I(\cX,H;P)$ is the intersection divisor of $\cX$
and $H$. 
    \end{enumerate}
    \end{corollary}

\subsection{Relation between linear series and morphisms}\label{s1.4} 
Define the following sets:
   \begin{itemize}
\item $\cL=\cL_r:=\{\cD^B: \text{$\cD$ linear series with 
$\dim(\cD)=r$}\};$

\item $\cM=\cM_r:=\{\langle \phi\rangle: \text{$\phi:\cX\to \P^r$ 
non-degenerate morphism}\}$, where\\
$\langle\phi\rangle:=\{T\circ\phi: T\in \aut(\P^r)\}$ denotes the
projective equivalent class of $\phi$.
   \end{itemize}
From Sects. \ref{s1.2} and \ref{s1.3} we have two maps, namely
   \begin{align*}
M=M_r & : \cL\to \cM;\qquad
\cD^B\mapsto \langle \text{coordinate representation of $\phi_{\cD^B}$}
\rangle\, ,\\
\intertext{and}
L=L_r & : \cM\to \cL;\qquad \langle\phi\rangle\mapsto \cD_\phi\, .
   \end{align*}  
We have $M\circ L = {\rm id}_{\cM}$ by definition, and $L\circ M = {\rm
id}_{\cL}$ by Lemma \ref{lemma1.22}. Therefore,
    
    \begin{lemma}\label{lemma1.41} 
The set of base-point-free linear series of dimension $r$ is equivalent to
the set of projective equivalent class of non-degenerate morphism from
$\cX$ to $\P^r$.
    \end{lemma} 

    \begin{remark}\label{rem1.41} The fact that $(L\circ M)(\cD^B)=\cD^B$ 
means that
  $$
\cD^B= \{\phi^*(H): H\ {\rm hyperplane\ in}\ \P^r\}\subseteq |E-B|\, ,
  $$
where $\phi:\cX\to \P^r$ is the non-degenerate morphism 
determinated, up to an automorphism of $\P^r$, by a base of $\cD'$. 
    \end{remark}

   \subsection{Hermitian invariants; Weierstrass semigroups I}\label{s1.5}

Let $\cD$ be a $g^r_d$ on $\cX$, say $\cD\cong \P(\cD')\subseteq |E|$, and
$P\in \cX$. We continue the study of the linear series $\cD_i(P)$ started
in Sect. \ref{s1.2}. Recall that $\cD_i(P)'=\cD'\cap L(E-iP)$ and that
$\cD_i(P)\supseteq \cD_{i+1}(P)$.

    \begin{definition*} A non-negative integer $j$ is called a {\em
$(\cD,P)$-order} (or an {\em Hermitian $P$-invariant}), if
$\cD_{j}(P)\supsetneqq \cD_{j+1}(P)$.
    \end{definition*}

From Lemma \ref{lemma1.21}, there exist $r+1$ $(\cD,P)$-orders, say
   $$
j_0(P)=j^{\cD}_0(P)<\ldots<j_r(P)=j^{\cD}_r(P)\, .
   $$
For $i=0,\ldots,r$, 
   $$
j_i(P)={\rm min}\{v_P(E)+v_P(f): f\in \cD_{j_i(P)}(P)'\}\, , 
   $$
and thus $\cD_{j_i}(P)$ is a $g^{r-i}_{d}$ on $\cX$. 

    \begin{lemma}{\rm (Esteves-Homma \cite[Lemma
1]{esteves-homma})}\label{lemma1.51} For $P, Q\in \cX$, $P\neq Q$, 
   $$
j_i(P)+j_{r-i}(Q)\le d\, .
   $$
    \end{lemma}    
     \begin{proof}
Since $\dim(\cD_{j_i(P)}(P)\cap \cD_{j_{r-i}(Q)}(Q))\ge 0$, there exists
$D\in \cD_{j_i(P)}(P)\cap \cD_{j_{r-i}(Q)}(Q)$ and the result follows.
     \end{proof}

This result will be complemented by Corollary \ref{cor2.23}.

     \begin{remark}\label{rem1.51} (i) Since $j_0(P)$ equals $b(P)$,
$\cD$ is base-point-free if and only if $j_0(P)=0$ for each $P\in \cX$.
Moreover, $j$ is a $(\cD,P)$-order if and only if $j-b(P)$ is a
$(\cD^B,P)$-order.

(ii) $j_r(P)\le d$ as $\cD_i(P)=\emptyset$ for $i>d$.  

(iii) Let $j\in \N_0$. From Lemma \ref{lemma1.21}, the following
statements are equivalent:
     \begin{enumerate}
\item[\rm(1)] $j$ is a $(\cD,P)$-order;
\item[\rm(2)] $\exists\ D\in \cD$ such that $v_P(D)=j$;
\item[\rm(3)] $\exists\ f\in \cD'$ such that $v_P(E)+v_P(f)=j$;
\item[\rm(4)] $\exists\ f\in \cD'$ such that $f\in L(E-jP)\setminus
L(E-(j+1)P)$;
\item[\rm(5)] $\dim_\F(\cD_j'(P))=\dim_\F(\cD_{j+1}'(P))+1$;
\item[\rm(6)] $\dim(\cD_j(P))=\dim(\cD_{j+1}(P))+1$.
     \end{enumerate}

(iv) Let $\cD=|E|$; i.e., $\cD'=L(E)$, $C$ 
a canonical divisor on $\cX$, and $j\in \N_0$. From 
$\cD_j'(P)=L(E-jP)$, the Riemann-Roch theorem, and part(iii)(5) above, 
the following statements are equivalent:
\begin{enumerate}
\item[\rm(1')] $j$ is a $(|E|,P)$-order;
\item[\rm(2')] $\exists\ f\in L(E)$ such that $v_P(E)+v_P(f)=j$;
\item[\rm(3')] $\exists\ f\in L(E-jP)\setminus L(E-(j+1)P)$;
\item[\rm(4')] $L(C-E+(j+1)P)=L(C-E+jP)$;
\item[\rm(5')] $\not\exists\ f\in L(C-E+(j+1)P)$ such that
$v_P(C-E)+v_P(f)=-(j+1)$.
\end{enumerate}
   \end{remark}

   \begin{example}\label{ex1.51} 
Let $g$ be the genus of $\cX$, and $\cD:=|E|$ with $d=\deg(E)\ge 2g$.
For $P\in \cX$, we compute some $(\cD,P)$-orders. We have $j_i(P)=i$ for
$0\le i\le d-2g$. Indeed for such an $i$, $\deg(C-E+(i+1)P)<0$ and then
Remark \ref{rem1.51}(iv(4')) is trivially satisfied. In particular, $\cD$
is base-point-free. 
   \end{example} 

    \begin{example}\label{ex1.52} We claim that for a given sequence of
non-negative integers $\l_0<\ldots<\l_r$, there exists a curve $\cY$, a
point $P_0\in\cY$, and a linear series $\cF$ on $\cY$ such that the
sequence equals the $(\cF,P_0)$-orders. Indeed, let $\cY:=\P^1(\F)$ and  
$x$ a transcendental element over $\F$. Set $P_\infty:=(0:1)$, and
$P_a:=(1:a)$ for $a\in \F$. We assume 
$\div(x)=P_0-P_\infty$, $v_{P_a}(x-a)=1$ for $a\in \F$. Define
   $$
E:=\l_rP_\infty,\qquad\text{and}\qquad \cF':=\langle
x^{\l_0},\ldots,x^{\l_r}\rangle\subseteq \F(x)\, .
   $$ 
Then $\cF:=\{E+\div(f):f\in \cF'\}$ is a $g^r_{\l_r}$ on $\cY$. We have
$E+\div(x^{\l_i})=\l_iP_0+(\l_r-\l_i)P_\infty$ and hence the
$(\cF,P_0)$-orders are $\l_0,\ldots,\l_r$. In addition, we have that
$j_0^\cF(P)=0$ for $P\neq P_0$; i.e., the base locus of $\cF$ is
$B^\cF=\l_0P_0$. Moreover, for the morphism associated to $\cF$
$\phi=(x^{\l_0}:\ldots:x^{\l_r})$ we have $E_\phi=\l_rP_\infty-\l_0P_0$. 
If $\ell_r=r$, then $\cF$ is complete and base-point-free, and
the curve
$\phi(\cY)$ is the so-called rational normal curve in $\P^r$. Conversely,
if $\cF$ is complete, say $\cF=|E_1|$, then $E_1=E$ by Lemma
\ref{lemma1.22}, and so $\ell=r$.
    \end{example}

We will introduce next the so-called Weierstrass semigroup. To begin with
we state a definition which is motivated by Remark \ref{rem1.51}(iv)(5').

     \begin{definition*} Let $D\in \Div(\cX)$ and $\ell\in \N_0$. We say
that $\ell$ is a {\em $(D,P)$-gap} if does not exist $f\in L(D+\ell P)$
such that $v_P(D)+v_P(f)=-\ell$.
      \end{definition*}

We have that
   $$
\text{$\ell$ is a $(D,P)$-gap}\quad\text{if and only if}\quad 
\text{$\ell-1$ is a $(|C-D|,P)$-order}\, ,
   $$ 
where $C$ is a canonical divisor on $\cX$. Denote by
$\cK=\cK_{\cX}:=|C|$ the canonical linear series on $\cX$.

   \begin{definition*} 
The $(0,P)$-gaps are called the {\em Weierstrass
gaps} at $P$. The {\em Weierstrass semigroup} at $P$ is the set
  $$
H(P):=\N_0\setminus G(P)\, ,
  $$
where
  $$
G(P):= \{\ell\in \Z^+: \ell\ \text{Weierstrass gap at}\ P\}\, .
  $$
The elements of $H(P)$ are called {\em Weierstrass non-gaps} at $P$.
    \end{definition*}

    \begin{lemma}\label{lemma1.52}
Let $g$ be the genus of $\cX$. Then 
   \begin{enumerate}
\item[\rm(1)] $\# G(P)=g$ (Weierstrass gap theorem)$;$
\item[\rm(2)] For $h\in \N_0$, the following statements are equivalent$:$
   \subitem\rm(i) $h\in H(P);$

   \subitem\rm(ii) $\exists\ f_h\in L(hP)$ such that $v_P(f_h)=-h;$
 
   \subitem\rm(iii) $\exists\ f_h\in k(X)$ such that
$\div_\infty(f_h)=hP;$

   \subitem\rm(iv) $\ell(hP)=\ell((h-1)P)+1.$
    \end{enumerate}
    \end{lemma}
    \begin{proof} Since $\dim(\cK)=g-1$ and
   $$
G(P)=\{j_0^\cK(P)+1,\ldots,j_{g-1}^\cK(P)+1\}\, ,
   $$ 
part (1) follows. Remark \ref{rem1.51}(iv) implies part (2).
     \end{proof}

We see now that $H(P)$ is indeed a semigroup.
     
     \begin{corollary}\label{cor1.51} 
The set $H(P)$ is a sub-semigroup of $(\N_0,+)$ such that
   $$
H(P)\supseteq \{2g,2g+1,2g+2,\ldots\}\, ,
   $$ 
where $g$ is the genus of $\cX$.
      \end{corollary}
      \begin{proof} It follows from 
Lemma \ref{lemma1.52}(2.(iii)) and $j_{g-1}^\cK(P)\le \deg(\cK)=2g-2$.
      \end{proof}

Let $(n_i(P):i=0,1,\ldots)$ denote the strictly increasing sequence that
enumerates the Weierstrass semigroup $H(P)$. From Lemma
\ref{lemma1.52}(2)(iv), 
$\ell(n_i(P)P)=i+1$ and from Corollary \ref{cor1.51}, 
$n_i(P)=g+i$ for $i\ge g$.

   \begin{remark}\label{rem1.511} For $g=0$, $\cK=\emptyset$ and hence
$H(P)=\N_0$ for any $P\in \cX$. If $g=1$, then $\dim(\cK)=0$ and hence
$H(P)=\{0,2,3,\ldots\}$ for any $P\in \cX$.
   \end{remark}    

     \begin{corollary}\label{cor1.52} If $\cX$ is a curve of genus $g\ge
1$, then $\cK$ is base-point-free.
      \end{corollary}
   \begin{proof} We have to show that $j_0(P):=j_0^{\cK}(P)=0$ for
each $P\in \cX$. Suppose that $j_0(P_0)\ge 1$ for some $P_0\in \cX$. Then
$1\in H(P_0)$ and hence $H(P_0)=\N_0$. This implies $g=0$.
   \end{proof}

   \begin{example}\label{ex1.53} We consider complete linear series on
$\cX$ arising from Weierstrass non-gaps which will be useful for
applications to optimal curves. Let $P\in \cX$, set $n_i:=n_i(P)$ and
consider $\cD:=|n_rP|$. Then
   \begin{enumerate}
 \item[\rm(1)] $\cD$ is a $g^r_{n_r}$ base-point-free on $\cX$;
 \item[\rm(2)] The $(\cD,P)$-orders are $n_r-n_i$, $i=0,\ldots,r$.
   \end{enumerate}

In fact, we already noticed that $\dim(\cD)=r$; $P$ cannot be a base point
of $\cD$ by Lemma \ref{lemma1.52}(2)(iv); if $Q\neq P$, then $D:=
n_rP+\div(1)\in \cD$ and $v_Q(D)=0$. This prove (1).  To see (2), let
$f_i\in \F(\cX)$ such that $\div(f_i)=\div_0(f_i)-n_iP$; cf. Lemma
\ref{lemma1.52}(2)(iii). Then
   $$
n_rP+\div(f_i)=(n_r-n_i)P+\div_0(f_i)\, ,
  $$
and the result follows.
    \end{example}
    \begin{lemma}\label{lemma1.520} Let $f\in\F(\cX)$ such that
$\div_\infty(f)=n_1(P)P$. Then $f$ is a separating variable of
$\F(\cX)|\F$.
    \end{lemma} 
    \begin{proof} If $\F(\cX)|\F(f)$ were not separable, then $f=g^p$,
$g\in\F(\cX)$ by \cite[Prop. III.9.2]{sti}. Then $n_1(P)/p$ would be a
non-gap at $P$, a contradiction.
    \end{proof}    

By definition, a Weierstrass semigroup $H(P)$ belongs to the class of {\em
numerical semigroup}; i.e., it is a sub-semigroup $H$ of $(\N_0,+)$ whose
complement in $\N_0$, $G(H):=\N_0\setminus H$, is finite. For such a
semigroup $H$, $g(H):=\# (\N_0\setminus H)$ is called the {\em genus} of
$H$. We let $(n_i(H):i\in\N)$ (resp. $(\ell_i(H):i= 1,\ldots,g(H))$)
denote the strictly increasing sequence that enumerates $H$ (resp.
$G(H)$). Clearly $n_i(H)=g(H)+i$ for $i\ge g(H)$, and $n_i(H)=2i$ for
$i=1,\ldots, g(H)$ whenever $n_1(H)=2$. $H$ is called {\em hyperellitpic}
if $2\in H$ (note that $2\in H$ if and only if $n_1(H)=2$, whenever
$g(H)\ge 1$). This definition is motivated by the so-called {\em
hyperelliptic curves}, namely those curves admitting a $g^1_2$, or
equivalently those admitting rational functions of degree two. Indeed,
$\cX$ is hyperelliptic if and only if there exists $P\in\cX$ such that
$2\in H(P)$ (see Example \ref{ex2.51}).

   \begin{lemma} {\rm (Buchweitz \cite[I.3]{buchweitz0}, Oliveira
\cite[Thm. 1.1]{oliveira})}\label{lemma1.521} If $n_1(H)\ge 3,$ then
$n_i(H)\ge 2i+1$ for $i=1,\ldots,g(H)-2$. In particular, $n_{g-1}(H)\ge
2g(H)-2.$
   \end{lemma}

The {\em weight} of $H$ is $w(H):=\sum_{i=1}^{g(H)}(\ell_i(H)-i)$. It is
easy to see that

   \begin{equation}\label{eq1.51} 
w(H)=(3g(H)^2+g(H))/2-\sum_{i=1}^{g(H)}n_i(H)\, ,
   \end{equation}

and that $w(H)=g(H)(g(H)-1)/2$ if $H$ is hyperelliptic. Now Lemma
\ref{lemma1.521} and (\ref{eq1.51}) imply:

   \begin{corollary}\label{cor1.521} 
   \begin{enumerate}
\item[\rm(1)] $0\leq w(H)\leq g(H)(g(H)-1)/2;$

\item[\rm(2)] $w(H)= g(H)(g(H)-1)/2$ if and only if $H$ is
hyperelliptic$;$

\item[\rm(3)] $w(H)\le (g(H)^2-3g(H)+4)/2$ if $n_1(H)\ge 3.$
   \end{enumerate}
   \end{corollary}

   \begin{remark}\label{rem1.52} {\rm (Kato \cite{kato})} If $n_1(H)\ge
3$, we indeed have $w(H)\le g(H)(g(H)-1)/3$, for $g(H)=3,4,6,7,9,10$ and 
$w(H)\le (g(H)^2-5g(H)+10)/2$, otherwise.
   \end{remark}

   \begin{definition*} A numerical semigroup $H$ is called {\em
Weierstrass} if there exist a curve $\cX$ and a point $P\in\cX$
such that $H$ equals the Weierstrass semigroup $H(P)$ at $P$.
   \end{definition*}

   \begin{remark}\label{rem1.53} If $H$ is Weierstrass, say $H=H(P)$ on a
curve $\cX$ of genus $g=g(H)$, then Lemma \ref{lemma1.521} follows from
Castelnuovo's genus bound (Lemma \ref{lemma1.24}): We want to show that
$n_i:=n_i(P)\geq 2i+1$ provided that $n_1:=n_1(P)\geq 3$ and $1\leq i\le
g-2$. Let $i$ be the least integer for which $n_i\leq 2i$. Then $i\ge 2$, 
$n_{i-1}=2i-1$, and $n_i=2i$. Thus $\cD:=|n_iP|$ is a simple $g^i_{n_i}$
on $\cX$; therefore Castelnuovo's genus bound implies $g\le i+1$, a
contradiction.
   \end{remark}

A numerical semigroup $H$ is Weierstrass if any of the following
conditions hold:

   \begin{itemize} 
\item either $g(H)\le 7$, or $g(H)=8$ and $2n_1(H)>\ell_g(H)$; see
Komeda \cite{komeda3}; 

\item $n_1(H)\le 5$; see Komeda \cite{komeda1}, \cite{komeda4}, Maclachlan
\cite[Thm. 4]{maclachlan};

\item either $w(H)\leq g(H)/2$ or $g(H)/2<w(H)\le g(H)-1$ and
$2n_1(H)>\ell_g(H)$; see Eisenbud-Harris \cite{eisenbud-harris}, Komeda
\cite{komeda2};
   \end{itemize}

We remark that the underlying curve in these examples is defined over the
complex numbers. 

In 1893, Hurwitz \cite{hurwitz} asked about the characterization of
Weierstrass semigroups; see \cite[p. 32]{buchweitz} and \cite[p.
499]{eisenbud-harris} for further historical information. Long after that,
in 1980 Buchweitz (see Corollary \ref{cor1.53}) showed the existence
of a non-Weierstrass semigroup as a consequence of the following.

    \begin{lemma} {\rm (Buchweitz's necessary condition, \cite[p.
33]{buchweitz})}\label{lemma1.53} Let $H$ be a numerical semigroup. For an
integer $n\ge 2$, let $nG(H)$ be the set of all sums of $n$
elements of $G(H)$. If $H$ is Weierstrass, then
   \begin{equation}\label{eq1.52}
\# nG(H)\le (2n-1)(g(H)-1)\, .
   \end{equation}
    \end{lemma}
    \begin{proof} We have that $g:=g(H)$ is the genus of the underlying
curve, say $\cX$. For a canonical divisor $C$ on $\cX$, we observe that
$\ell(nC)=(2n-1)(g-1)$ by the Riemann-Roch theorem. Let
$\ell:=\ell_1+\ldots+\ell_n\in nG(H)$. From Remark \ref{rem1.51}(iv)(2'),
there exists $f_i\in L(C)$ such that $v_P(C)+v_P(f_i)=\ell_i-1$ for
$i=1,\ldots,n$. Then $f_\ell:=f_1\ldots f_n\in L(nC)$ and being the map
$\ell\mapsto f_\ell$ injective, the result follows.
    \end{proof}

    \begin{corollary} {\rm (\cite[p. 31]{buchweitz})}\label{cor1.53} 
$\{1,\ldots, 12,19,21,24,25\}$ is the set of gaps of a numerical
semigroup $H$ of genus 16 which is not Weierstrass.
    \end{corollary}
    \begin{proof} We apply the case $n=2$ in Lemma \ref{lemma1.53}. An
easy computations shows that $2G(H)=[2,50]\setminus\{39,41,47\}$. Then 
$\#2G(H)=46>3g-3=45$ and so $H$ cannot be Weierstrass.
    \end{proof} 

In addition, Buchweitz (loc. cit.) showed that for every integer $n\ge 2$
there exist numerical semigroups which do not satisfy (\ref{eq1.52}).
Further examples of such semigroups were given in \cite[Sect.
4.1]{torres2} and Komeda \cite{komeda5}. On the other hand, what can we
say about semigroups $H$ that satisfy (\ref{eq1.52}) for each $n\ge 2$ ?
In fact, there exist at least two classes of such semigroups, namely {\em
symmetric semigroups} (resp. {\em quasi-symmetric semigroups}); i.e.,
those $H$ with $\ell(H)=2g(H)-1$ (resp. $\ell(H)=2g(H)-2$). Indeed,
equality in (\ref{eq1.52}) for each $n$ characterize symmetric semigroups
(see Oliveira \cite[Thm. 1.5]{oliveira}), and Oliveira and St\"ohr
\cite[Thm. 1.1]{oliveira-stohr} noticed that $\#nG(H)=(2n-1)(g-1)-(n-2)$
whenever $H$ is quasi-symmetric. In 1993, St\"ohr \cite[Scholium
3.5]{torres1} constructed symmetric semigroups which are not Weierstrass.
Indeed, symmetric non-Weierstrass semigroups of any genus larger than 99
can be constructed (loc. cit.) by using the Buchweitz's semigroup
(Corollary \ref{cor1.53}) as a building block. A similar result was
obtained for quasi-symmetric semigroups \cite[Thm. 5.1]{oliveira-stohr}
and these examples were generalized in \cite[Sect. 4.2]{torres2}. We
stress that any symmetric (resp. quasi-symmetric) semigroup is a
Weierstrass semigroup on a Gorenstein (resp. reducible Gorenstein) curve;
see \cite{stohr} (resp. \cite{oliveira-stohr}).

Finally, we mention that Hurwitz's question for numerical semigroups that
satisfy (\ref{eq1.52}) for each $n\ge 2$ is currently an open problem.

    \section{Weierstrass point theory}\label{s2}

In this section we study Weierstrass Point Theory of linear series on
curves from St\"ohr-Voloch's paper \cite[\S1]{sv}. Other references are
Farkas-Kra \cite[III.5]{farkas-kra}, Homma \cite[Sects. 1,2]{homma0},
Laksov \cite{laksov}, F.K. Schmidt \cite{sch0}, \cite{sch}. 

Let $\cX$ be a curve over an algebraically closed field $\F$
of characteristic $p\ge 0$. Let $\cD$ be a $g^r_d$ on $\cX$, say $\cD\cong
\P^r(\cD')\subseteq |E|$. 

In Sect. \ref{s1.5}, to any point $P\in \cX$ we have assigned a sequence
of $(r+1)$ integers, namely the $(\cD,P)$-orders. Here we study the
behaviour of such sequences for general points of $\cX$; i.e, for points
in an open Zariski subset of $\cX$. In order to do that we use
``wronskians" on $\cX$; i.e., certain functions in $\F(\cX)$ defined via
derivatives. To avoid restrictions on the characteristic $p$, we use Hasse
derivatives.
     \subsection{Hasse derivatives}\label{s2.1} 

Let $x$ be a trascendental element over $\F$. For $i,j \in \N_0$, set 
   $$
D^i_xx^j:=\binom{j}{i}x^{j-i}\, , 
   $$ 
and extend it $\F$-linearly on $\F[x]$. The $\F$-linear map $D^i_x$ is
called the {\it $i$-th Hasse derivative} on $\F[x]$. $i!\
D^i_xx^j$ is the usual $i$-th derivative $\frac{d^i}{dx^i}$, and 
$D^i_x\neq 0$, as $D^i_xx^i=1$, but $\frac{d}{d^ix}=0$ for $i\ge
p>0$. 

   \begin{remark}\label{rem2.11} 
For $f(x)\in \F[x]$, $D^i_xf(x)$ is the 
coefficient of $u^i$ in the expansion of $f(x+u)$ as a polynomial in $u$. 
   \end{remark}

The $\F$-linear maps $D^i_x$, $i\in\N_0$, satisfy the following four
properties:
     \begin{enumerate} 
\item[\bf(H1)] $D^0_x={\rm id}$;
\item[\bf(H2)] ${D^i_x}_{|\F} =0$ for $i\ge 1$; 
\item[\bf(H3)] $D^i_x(fg)=\sum_{j=0}^{i}D^j_xfD^{i-j}_xg$ (Product Rule);
\item[\bf(H4)] $D^i_x\circ D^j_x = \binom{i+j}{i}D^{i+j}_x$. 
     \end{enumerate}

Properties (H1), (H2) and (H4) easily follow from the definition of
$D^i_x$, while (H3) follows by comparing the coefficients of $(fg)(x+u)$
and $f(x+u)g(x+u)$.

Next one extends $D^i_x$ to $\F(x)$ and then to each finite separable
extension of $\F(x)$. This is done in just one way; moreover, the extended
map remains $\F$-linear and still satisfies the four
aforementioned properties. The
extension on $\F(x)$ is constructed as follows. By (H1) and (H3) it is
enough to define $D^i_x(1/f)$ for $i\ge 1$ and $f\neq 0$. From $f(1/f)=1$,
(H2) and (H3) one finds the following recursive formula:
    $$
\sum_{j=0}^{i}D^j_x(1/f)D^{i-j}_xf=0\, .
    $$
For $i=1$ one obtains the expected relation 
$D^1_x(1/f)=-(D^1_xf)/f^2$, and in general \cite[p. 119]{niederreiter}
     \begin{equation*}
D^i_x(1/f)=\sum_{j=1}^i\frac{(-1)^j}{f^{j+1}}\sum_{
i_1,\ldots,i_j\ge 1;\ i_1+\ldots+i_j=i} D^{i_1}_xf
\ldots D^{i_j}_xf\, . 
     \end{equation*}
     \begin{remark}\label{rem2.12} The maps $D^i_x$ on $\F(x)$,
$i\in\N_0$, are characterized by the following four properties: 
    \begin{enumerate}
\item[\rm(i)] they are $\F$-linear;
\item[\rm(ii)] they satisfy (H1) and (H3) above;
\item[\rm(iii)] $D^1_xx=1$;
\item[\rm(iv)] $D^i_xx=0$ for $i\ge 2$. 
   \end{enumerate}
To see this, let $\eta_i$, $i\in\N_0$, be maps on $\F(x)$
satisfying (i), (ii), (iii) and (iv). From the formula for $D^i_x(1/f)$
above, is enough to show that $\eta_i(x^j)=D^i_xx^j\ (*)$ for
$i,j\in\N_0$. Now, since the $\eta_i$'s satisfy (H3), it follows  
\cite[Lemma 3.11]{hefez}
    \begin{equation}\label{eq2.11}
\eta_i(x^j)=
jx^{j-1}\eta_i(x)+\sum_{\ell=2}^{j}\sum_{m=1}^{i-1}x^{j-\ell}(\eta_m(x))
(\eta_{i-m}(x^{\ell-1}))\, ,
    \end{equation}
and we obtain $(*)$ by induction on $i$ and $j$.	
     \end{remark}
     \begin{remark}\label{rem2.13} The maps $D^i_x$, $i\in\N_0$, on
$\F(x)$ have also a unique extension to the Laurent series $\F((x))$ which
satisfy (H1), (H2), (H3), and (H4) above. One sets
$D^i_x(\sum_ja_jx^j):=\sum_j\binom{j}{i}a_jx^{i-j}$, see \cite[p.  
12]{hefez}.
     \end{remark}
Next we extend $D^i_x$ to a finite separable extension $\K|\F(x)$. Let
$y\in \K$ be such that $\K=\F(x,y)$, and $F(x)[Y]$ the minimal
polynomial of $y$ over 
$\F(x)$. Then we define $D^i_xy^m$ by using $F(x,y)=0$ and
(\ref{eq2.11}). For example, for $i=1$ we obtain
     \begin{equation}\label{eq2.12}
F_Y(x,y)D^1_xy+\sum_j(D^1_xa_j(x))y^j=0\, ,
     \end{equation}
so that $D^1_xy$ is well defined as $F_Y(x,y)\neq0$. Notice that 
these extensions satisfy (H1), (H2), (H3) and (H4) above and depend on
the element $y$. However, it is a matter of
fact that the $\F$-linear maps $D^i_x$ on $\F(x)$ admit a unique
extension to $\F$-linear maps on $\K$ satisfying the aforementioned (H1),
(H2), (H3), and (H4); see \cite{hasse-schmidt}. 

Therefore, $\F(\cX)$ is equipped with $\F$-linear maps $D^i_x$ such 
that (H1), (H2), (H3) and (H4) above hold true, with $x$ being a
separating variable of $\F(\cX)|\F$. If $y$ is another separating variable
of $\F(\cX)|\F$, relations among the $D^i_x$'s and the $D^j_y$'s are given
by the so called {\it chain rule}; see (\ref{eq2.13}) and
(\ref{eq2.14}). 
     \begin{remark}\label{rem2.14} 
For $i\in\N_0$, let $D^i$ be $\F$-linear maps on a $\F$-algebra
$\K$ satisfying (H1), (H2), (H3) and (H4) above. From (H4), 
   $$
i!\ D^i=(D^1)^i:=D^1\circ \ldots\circ D^1\quad \text{$i$ times}\, ,
   $$ 
so that each $D^i$ is determinated by $D^1$ provided that $p=0$. 
Suppose now $p>0$.
     \begin{claim*} Let $0\le a, b<p$, $\alpha,\
\beta \in \N$. Then
   \begin{enumerate} 
\item[\rm(1)] $D^{ap^\alpha+bp^\beta}=D^{ap^\alpha}\circ
D^{bp^\beta}$.
\item[\rm(2)] $D^{ap^\alpha}=(D^{p^\alpha})^a/a!$.
   \end{enumerate}
     \end{claim*}
     \begin{proof} The statements are consequence of (H4) and the
following property of binomial numbers: if 
$i=\sum_\alpha a^\alpha p^\alpha$, $j=\sum_\alpha b^\alpha p^\alpha$ are  
the $p$-adic expansion of $i, j\in \N$, then 
$\binom{i}{j}=\prod_\alpha\ \binom{a_\alpha}{b_\alpha}$.
      \end{proof} 
Therefore in positive characteristic the $D^i$'s are determinated by $D^1,
D^p, D^{p^2},\ldots $. 
     \end{remark}
A $\F$-linear map $D$ on $\F(\cX)$ satisfying $D(fg)=fD(g)+gD(f)$,   
is called a {\em $\F$-derivation} on $\F(\cX)$. For example, $D^1_x$ is a
derivation on $\F(\cX)$, where $x$ is a separating
variable of $\F(\cX)|\F$. From (\ref{eq2.11})  
follows that two $\F$-derivations $\delta_1$ and $\delta_2$ on $\F(\cX)$
are equal if $\delta_1(x)=\delta_2(x)$.

Now let $y$ be another separating variable of $\F(\cX)|\F$. Since the
$\F$-derivations $\delta_1:=D^1_y$ and $\delta_2:=D^1_y(x)D^1_x$ satisfy 
$\delta_1(x)=\delta_2(x)$, we obtain the usual chain rule, namely
   \begin{equation}\label{eq2.13}
D^1_y\ =\ D^1_y(x)D^1_x\, .
   \end{equation}
To generalize this relation to higher derivatives, let $T$ be a
trascendental element over $\F(\cX)$. The maps $D^i_x$ and $D^j_y$ can be 
read off from the homomorphisms of $\F$-algebras $\eta_x, \eta_y$:  
$\F(X)\to \F(X)[[T]]$ defined respectively by 
   
   $$
\eta_x(f):= \sum_{i\ge 0}D^i_x(f)T^i\, ,\quad\text{and}\quad 
\eta_y(f):=\sum_{i\ge 0}D^i_y(f)T^i\, .   
   $$

Let $h:\F(\cX)[[T]]\to \F(\cX)[[T]]$ be the $\F$-homomorphism defined by 
$h_{|\F(\cX)}={\rm id}_{|\F(\cX)}$ 
and $h(T):= \sum_{i\ge 1}D^i_y(x)T^i$. Since $D^1_y(x)\neq 0$ by
(\ref{eq2.13}), $h$ is an 
automorphism of $\F(\cX)[[T]]$. Consider the $\F$-homomorphism
$\eta:\F(X)\to
\F(X)[[T]]$ given by $\eta:= h^{-1}\circ \eta_y$. 
For $f\in \F(\cX)$, set $\eta(f):=\sum_{i\ge 0}\eta_i(f)T^i$. Then the
maps  
$\eta_i$ are $\F$-linear on $\F(\cX)$ and satisfy properties 
(H1) and (H3) above. Write $h(T)=TU$, $U=D^1_y(x)+D^2_y(x)T+\ldots$. 
     \begin{claim*} Let $i\in \N_0$ and $f\in \F(\cX)$. Then
$\eta_0(f)=D^0_y(f)$ and for $i\ge 1$ the following holds
   $$
D^i_y(f)=\sum_{j=1}^{i}a_j\eta_j(f)\, ,
   $$
where $a_j$ is the coefficient of $T^{i-j}$ in $U^j$. In particular,
$a_1=D^i_y(x)$, $a_i=(D^1_yx)^i$.
     \end{claim*}
     \begin{proof}
Write $\eta_y=h\circ\eta$. The coefficient of $T^i$ in $(h\circ\eta)(f)$
can be read off from $\sum_{j=0}^{i}a_j(f)(TU)^j$, and the claim follows.
     \end{proof}
Then we have $\eta_1(x)=1$ and $\eta_i(x)=0$ for $i\ge 2$. Therefore from 
Remark \ref{rem2.12}, $\eta_i=D^i_x$ on $\F(x)$ and hence also on
$\F(\cX)$. This implies the generalized chain rule:
   $$
\eta_y=h\circ\eta_x\, ,
   $$
or equivalently 
       \begin{equation}\label{eq2.14}
D^i_y=\sum_{j=1}^{i}f_jD^j_x\, ,\qquad i=1,2\ldots\, ,
      \end{equation}
where $f_j\in \F(\{D^m_y(x):m=1,2,\ldots\})$. Observe that $f_1=D^i_y(x)$
and $f_i=(D^1_yx)^i$.

    \begin{remark}\label{rem2.15} We mention two further properties of
Hasse derivatives regarding prime powers of rational functions. Let
$f\in\F(\cX)$, $x$ a separating variable of $\F(\cX)|\F$, and $q$ a power
of $p={\rm char}(\F)>0$. We have
  \begin{enumerate} \item[\rm(i)] $D_x^if^q=(D_x^{i/q} f)^q$ if $q$
divides $i$, and $D_x^if^q=0$ otherwise;
  \item[\rm(ii)] (\cite[Satz 10]{hasse-schmidt}) $\exists\ g\in\F(\cX)$
such that $f=g^q$ if and only if $D^i_x(f)=0$ for $i=1,\ldots,q-1$.
  \end{enumerate}
   \end{remark}

    \begin{definition*}
A {\em wronskian} on $\cX$ is a rational function of type
   $$
W^{\l_0,\ldots,\l_r}_{f_0,\ldots,f_r;x}:=
\det((D^{\l_i}_xf_j))\, ,
   $$
where $\l_0<\ldots<\l_r$ is a sequence of non-negative
integers, $x$ is a separating variable of $\F(\cX)|\F$, and    
$f_0,\ldots,f_r\in \F(\cX)$. We set
   $$
\cA(f_0,\ldots,f_r;x):=\{(m_0,\ldots,m_r)\in \N_0^{r+1}:
m_0<\ldots<m_r;\ W^{m_0,\ldots,m_r}_{f_0,\ldots,f_r;x}\neq 0\}\, .
   $$
   \end{definition*}
\subsection{Order sequence; Ramification divisor}\label{s2.2} 

Let $P\in \cX$ and $t$ be a local parameter at $P$. Let 
   $$
j_0=j_0(P)<\ldots <j_r=j_r(P)
   $$ 
denote the $(\cD,P)$-orders. From Remark \ref{rem1.51}(iii)(3)  
there exists $f_\l\in \F(\cX)$ such that
   $$
v_P(t^{v_P(E)}f_\l)=j_\l\, ,\qquad \l=0,\ldots,r\, .
   $$
      \begin{claim*} $\{f_0,\ldots,f_r\}$ is a $\F$-base of $\cD'$.
      \end{claim*}
      \begin{proof} 
If there exists a non-trivial relation $\sum_ia_if_i=0$ 
with $a_i\in \F$, then we would have $v_P(f_i)=v_P(f_\l)$ for $i\neq
\l$ and so $j_i=j_\l$, a contradiction.
      \end{proof}

      \begin{definition*} The aforementioned $\F$-base
$\{f_0,\ldots,f_r\}$ is called a {\em $(\cD,P)$-base} (or {\em
$(\cD,P)$-Hermitian base}).
     \end{definition*}

     \begin{remark}\label{rem2.21} 
Let $\{f_0,\ldots,f_r\}$ be a $(\cD,P)$-base. For 
$i=0,\ldots,r$, $\cD_i'(P)=\cD'\cap L(E-j_iP)$ so that 
    $$
\cD_{j_i}'(P)=\langle f_i,\ldots,f_r\rangle\, ,
    $$
or equivalently
    $$
\cD_{j_i}(P) =\{E+\div(\sum_{\l=i}^{r}a_\l f_\l): (a_i:\ldots:a_r)\in
\P^{r-i}(\F)\}\, .
    $$
Thus
    $$
j_i(P) ={\rm min}\{v_P(\sum_{\l=i}^{r}a_\l f_\l t^{v_P(E)}):
(a_i:\ldots:a_r)\in \P^{r-i}(\F)\}\, .
    $$
     \end{remark}
Let $\{f_0,\ldots,f_r\}$ be a $(\cD,P)$-base. Set $g_\l:=t^{v_P(E)}f_\l$.

     \begin{lemma}\label{lemma2.21} If $m_0<\ldots<m_r$ is a sequence of
non-negative integers such that $\det(\binom{j_\l}{m_i})\not\equiv
0\pmod{p}$, then $(m_0,\ldots,m_r)\in \cA(g_0,\ldots,g_r;t)$. In
particular, $(j_0,\ldots,j_r)\in \cA(g_0,\ldots,g_r;t)$.
     \end{lemma}
     \begin{proof} 
Let $g_\l=\sum_{s=j_\l}^{\infty}c^{\l}_s t^s$, $c^\l_{j_\l}\neq 0$, be the
local expansion of $g_\l$ at $P$. Set $C:=\prod_{\l=0}^{r}c^{\l}_{j_\l}$.
Then 
     \begin{align*}  
W^{m_0,\ldots,m_r}_{g_0,\ldots,g_r;t} & =\det(\sum_{s=j_\l}^{\infty}
\binom{s}{m_i}c^{\l}_st^{s-m_i})\\
                           {}           & =Ct^{-\sum_im_i}
\det(\sum_{s=j_\l}^{\infty}\binom{s}{m_i}\frac{c^{\l}_s}{c^\l_{j_\l}}t^s)\\
      {}  & =C\det(\binom{j_\l}{m_i})t^{\sum_i(j_i-m_i)}+\ldots\neq 0\, ,
     \end{align*}
and the result follows.
     \end{proof}

For $\ell\in\N_0$, set $D^{\l}_x\phi:= (D^{\l}_xg_0,\ldots,D^{\l}_xg_r)$.
Since each coordinate of this vector is regular at $P$, we also set
$D^{\l}_x\phi(P):= (D^{\l}_xg_0(P),\ldots,D^{\l}_xg_r(P))$.

Then, for $0\leq m_0<\ldots<m_r$, $(m_0,\ldots,m_r)\in
\cA(g_0,\ldots,g_r;t)$ if and only if
$D^{m_0}_t\phi,\ldots,D^{m_r}_t\phi$ are $\F(\cX)$-linearly independent.

    \begin{scholium}\label{scholium2.21}
\begin{enumerate}
\item[\rm(1)] Set $j_{-1}:=0$. For $i=0,\ldots,r$,
    $$
j_i=j_i^\cD(P)=\min\{s> j_{i-1}: \text{$(D^{j_0}_t\phi)(P),
\ldots,(D^{j_{i-1}}_t\phi)(P),(D^s_t\phi)(P)$ are $\F$-l.i.}\}\, ;
    $$
\item[\rm(2)] Let $m_0<\ldots<m_{r'}$ be non-negative integers, with
$r'\leq r$, such that the 
vectors $(D^{m_0}_t\phi)(P),\ldots,(D^{m_{r'}}_t\phi)(P)$
are $\F$-linearly independent. Then $j_i\le m_i$ for $i=0,\ldots,r'.$
\end{enumerate}
    \end{scholium}
    \begin{proof} (1) From Lemma \ref{lemma2.21} and its proof, the
vectors $(D^{j_0}_t\phi)(P),\ldots,
(D^{j_i}_t\phi)(P)$ are $\F$-linearly independent and
   \begin{equation*}
D^{j_i}_tg_\l(P)=\begin{cases}
0 & \text{if $\l>i$}\, ,\\
c^\l_{j_\l} & \text{if $\l=i$}\, ,\\
c^{\l}_{j_i} & \text{if $\l<i$}\, .\end{cases}
   \end{equation*}
Let $j_{i-1}<s<j_i$. For $\l=0,\ldots,i-1$, we have vectors of type 
   $$
(D^{j_{\l}}_t\phi)(P)=(*,\ldots,*,c^{\l}_{j_{\l}},0,\ldots,0)\, ,
   $$
with $(r-\l)$ zeros and where $*$ denotes an element of $\F$. Since the 
last $(r-i+1)$ entries of the vector $(D^s_t\phi)(P)$ are zeroes, 
(1) follows.

(2) From (1), $\dim_\F \langle\{(D^s\phi)(P): s=0,\ldots,j_i-1\}
\rangle =i$ so that $j_i-1<m_i$.
    \end{proof}

In $\Z^{r+1}$ we have a partial order given by the so-called {\em
lexicographic order} $<$. For $\alpha, \beta\in \Z^{r+1}$, $\alpha<\beta$
if in the vector  $\beta-\alpha$ the left most non-zero entry is
positive. This order is a well-ordering on $\N^{r+1}$, see e.g. \cite[p.  
55]{cox}. Let
   $$
\cE:=(\epsilon_0,\ldots,\epsilon_r)
   $$
be the minimum (in the lexicographic order) of $\cA(g_0,\ldots,g_r;t)$. 

     \begin{lemma}\label{lemma2.22} 
   \begin{enumerate}
\item[\rm(1)] $\epsilon_0=0;$
\item[\rm(2)] $\epsilon_1=1$ whenever $p$ does not divide
$\deg(\cD)-\deg(B^\cD);$
\item[\rm(3)] For $i=1,\ldots, r,$
   $$
\epsilon_i=\min\{s>\epsilon_{i-1}:
\text{$D^{\epsilon_0}_t\phi,\ldots,D^{\epsilon_{i-1}}_t\phi, D^s_t\phi$
are $\F(\cX)$-l.i.}\}\, .
   $$
   \end{enumerate}
     \end{lemma}
     \begin{proof} (1) Suppose that $\epsilon_0>0$. Then
$D^0_t\phi=\sum_{j=1}^{r}h_jD^{\epsilon_j}_t\phi$ with some $h_{j_0}\in
\F(\cX)^*$, because 
$(0,\epsilon_1,\ldots,\epsilon_r)<\cE$. 
Then  we replace the row $D^{\epsilon_{j_0}}_t\phi$ by $D^0_t\phi$ 
in $W^{\epsilon_0,\ldots,\epsilon_r}_{g_0,\ldots,g_r;t}$ 
so that 
$(0,\epsilon_0,\ldots,\epsilon_{j_0-1},\epsilon_{j_0+1},\ldots,\epsilon_r)\in
\cA(g_0,\ldots,g_r;t)$, a contradiction to the minimality of $\cE$.

(2) As in part (1) we have that $\epsilon_1=0$ if and only if
$D_t^1g_\l=0$ (or equivalently $D_t^ig_\l=0$ for $1\le i<p$) for any
$\l=0,\ldots,r$. Then each $g_\l$ is a $p$-power by Remark
\ref{rem2.15}(ii), and so $p$ divides $v_P(E)-b(P)$ by Lemma
\ref{lemma1.22}; i.e., $p$ divides
$\deg(\cD)-\deg(B^\cD)$.

(3) Clearly $D^{\epsilon_0}_t\phi,\ldots,D^{\epsilon_i}_t\phi$ are
$\F(\cX)$-linearly independent. Let $\epsilon_{i-1}<s<\epsilon_i$. Since 
$(\epsilon_0,\ldots,\epsilon_{i-1},s,\epsilon_{i+1},\ldots,\epsilon_r)<\cE$, 
there exists a relation of type
   $$
D^s_t\phi=\sum_{j=0}^{i-1}h_jD^{\epsilon_j}_t\phi+
\sum_{j=i+1}^{r}h_jD^{\epsilon_j}_t\phi\, ,
   $$
with $h_j\in \F(\cX)$. We claim that $h_j=0$ for $j\ge i+1$. Indeed,
suppose that $h_{j_0}\neq 0$ for some $j_0\ge i+1$. Then by replacing
$D^{\epsilon_{j_0}}_t\phi$ by $D^s_t\phi$ in
$W^{\epsilon_0,\ldots\epsilon_r}_{g_0,\ldots,g_r;t}$ we would have that 
$(\epsilon_0,\ldots,\epsilon_{i-1},s,\epsilon_i,\ldots,
\epsilon_{j_0-1},\epsilon_{j_0+1},\ldots,\epsilon_r)\in
\cA(g_0,\ldots,g_r;t)$, a contradiction to the minimality of $\cE$. This
finish the proof.
     \end{proof}

     \begin{corollary}\label{cor2.21} 
\begin{enumerate} 
   \item[\rm(1)] Let $(m_0,\ldots,m_r)\in
\cA(g_0,\ldots,g_r;t)$. Then for each $i$, $\epsilon_i\le m_i$.
In particular, $\epsilon_i\le j_i=j_i(P);$
   \item[\rm(2)] If $0\le m_0<\ldots<m_r$ are integers such that
$\det(\binom{j_i}{m_\l})\not\equiv 0\pmod{p}$, then
$\epsilon_i\le m_i$ for each $i.$
   \end{enumerate}
     \end{corollary}
     \begin{proof} From Lemma \ref{lemma2.22},  
     \begin{equation}\label{eq2.21}
\langle\{D^{\l}_t\phi:\l=0,\ldots,\epsilon_i-1\}\rangle=\langle
\{D^{\epsilon_j}_t\phi: j=0,\ldots,i-1\}\rangle\, .
     \end{equation}
If $\epsilon_i>m_i$, we would have 
   $$
\dim_{\F(\cX)}(\{D^{\l}_t\phi:\l=0,\ldots,\epsilon_i-1\})\ge
\dim_{\F(\cX)}(\{D^{m_{\l}}_t\phi: \l=0,\ldots,i\})\geq i+1\, ,
    $$
a contradiction. This proves (1). Now (2) follows from Lemma
\ref{lemma2.21} and (1).
     \end{proof}

    \begin{proposition}\label{prop2.21} 
\begin{enumerate}
   \item[\rm(1)] If $h_i=\sum a_{ij}g_j$ with
$(a_{ij})\in M_{r+1}(\F)$, then
    $$
W^{\epsilon_0,\ldots,\epsilon_r}_{h_0,\ldots,h_r;t}
=\det((a_{ij}))
W^{\epsilon_0,\ldots,\epsilon_r}_{g_0,\ldots,g_r;t}\, ;
    $$

\item[\rm(2)] If $f\in \F(\cX)$, then
    $$
W^{\epsilon_0,\ldots,\epsilon_r}_{fg_0,\ldots,fg_r;t}=f^{r+1}
W^{\epsilon_0,\ldots,\epsilon_r}_{g_0,\ldots,g_r;t}\, ;
    $$

\item[\rm(3)] Let $x$ be any separating variable of $\F(\cX)|\F$. Then
   $$
W^{\epsilon_0,\ldots,\epsilon_r}_{g_0,\ldots,g_r;x}=(D^1_xt)^{\sum_i
\epsilon_i}
W^{\epsilon_0,\ldots,\epsilon_r}_{g_0,\ldots,g_r;t}\, .
   $$
    \end{enumerate}
    \end{proposition}
    \begin{proof} (1) It follows from $D^{\epsilon_{\l}}_th_i=\sum
a_{ij}D^{\epsilon_\l}_tg_j$. Note that this result does not depend on the
minimality of $\cE$. 

(2) By the product rule (cf. Sect. \ref{s2.1}), we have
   $$
D^{\epsilon_i}_t(fg_j)= 
\sum_{\l=0}^{\epsilon_i}D^\l_tfD^{\epsilon_i-\l}_tg_j\, .
   $$
Then 
   $$
(D^{\epsilon_i}_t fg_0,\ldots,D^{\epsilon_i}_t fg_r)=
fD^{\epsilon_i}_t\phi+\sum_{\l=1}^{\epsilon_i} D^\l_tf
D^{\epsilon_i-\l}_t\phi\, .
   $$
By (\ref{eq2.21}) we can factor out $f$ in each row of 
$W^{\epsilon_0,\ldots,\epsilon_r}_{fg_0,\ldots,fg_r;t}$, and (2) follows.

(3) The proof is similar to (2) but here we use the chain
rule (\ref{eq2.14}) instead of the product rule. We have
  $$
D^{\epsilon_i}_x g_j=\sum_{\l=1}^{\epsilon_i}f_\l D^{\l}_tg_j\, ,
  $$ 
where $f_\l \in \F(\cX)$ and $f_{\epsilon_i}=(D^1_xt)^{\epsilon_i}$. Hence
  $$
D^{\epsilon_i}_x\phi=(D^1_xt)^{\epsilon_i}D^{\epsilon_i}_t\phi
+\sum_{\l=1}^{\epsilon_i-1}f_\l D^\l_t\phi\, ,
  $$
and again by (\ref{eq2.21}) we can factor out $(D^1_xt)^{\epsilon_i}$ in
each row of $W^{\epsilon_0,\ldots,\epsilon_r}_{g_0,\ldots,g_r;x}$.
      \end{proof}

Now we see that $\cE$ depends only on $\cD$: Let $f_0',\ldots,f_r'$
be any $\F$-base of $\cD'$ and $x$ any separating variable of
$\F(\cX)|\F$; since $g_{\l}=t^{v_P(E)}f_{\l}$, from Proposition
\ref{prop2.21}(1)(2)  $\cE$ is the minimum for $\cA(f_0',\ldots,f_r';t)$.
Moreover by part (3) of that proposition, $\cE$ is also the minimum for
$\cA(g_0,\ldots,g_r;x)$. Finally, from part (2), $\cE$ is also the minimum
for $\cA(f_0',\ldots,f_r';x)$.

     \begin{definition*} $\cE=\cE_\cD$ is called the {\em order
sequence} of $\cD$. The {\em order sequence} of a morphism $\phi$ is the
order sequence of $\cD_\phi$.
     \end{definition*}

     \begin{remark}\label{rem2.22} Let $m_0<\ldots<m_r$ be a sequence of
non-negative integers such that $\det(\binom{j_\l}{m_i})\not\equiv
0\pmod{p}$. Then $\epsilon_i\leq m_i$ for each $i$ by Corollary
\ref{cor2.21}(2). We shall discuss the best election of the $m_i$'s. In
Example \ref{ex1.52} we have seen that the $(\cD,P)$-orders
$j_0<\ldots<j_r$
are the $(\cD_\phi,P_0)$-orders for
$\phi=(x^{j_0}:\ldots:x^{j_r}):\P^1(\F)\to\P^{j_r}$ and $P_0=(1:0)$.
Observe
that
   \begin{equation}\label{eq2.22}
W^{n_0,\ldots,n_r}_{x^{j_0},\ldots,x^{j_r};x}=
\det(\binom{j_\l}{n_i})x^{\sum_i(j_i-n_i)}\, .
   \end{equation}
Let $\eta_0,\ldots,\eta_r$ be the $\cD_\phi$-orders. Then
   \begin{enumerate}
\item[\rm(1)] $\det(\binom{j_i}{\eta_\l})\not\equiv 0\pmod{p}$ by
(\ref{eq2.22})
with $n_i=\eta_i$, and the definition of $\cD_\phi$-orders;

\item[\rm(2)] $\eta_\l\leq m_\l$ for each $\l$ by (\ref{eq2.22}) with
$n_i=m_i$, and Corollary 2.8(2).
   \end{enumerate}
This shows that the best way to upper bound the $\epsilon_i$'s is by means
of 
the sequence $\eta_0,\ldots,\eta_r$. In addition, from
(\ref{eq2.22}) and Lemma \ref{lemma2.22} applied to $\cD_\phi$, we
obtain the following.
   \end{remark}    

     \begin{corollary}\label{cor2.22} 
Let $i\in\{0,\ldots,r\}$ and let $m_0<\ldots<m_i$ be non-negative 
integers, such 
that the vectors $(\binom{j_0}{m_\l},\ldots,\binom{j_r}{m_\l})$,    
$\l=0,\ldots,i$ are $\F_{p}$-linearly independent. Then $\epsilon_\l\le
m_\l$ for $\l=0,\ldots,i$.
     \end{corollary}
 
     \begin{corollary} {\rm (Esteves, \cite{esteves})}\label{cor2.23}
  $$
\epsilon_i+j_\l(P)\leq j_{i+\l}(P)\, ,\qquad i+\l\leq r\, .
  $$
     \end{corollary}
     \begin{proof} (Following Homma \cite{homma1}) By means of suitable 
central projections \cite[Lemma 2]{esteves} one can assume that
$i+\l=r$. Let
$\cD_\phi$ be the linear series on $\P^1(\F)$ in Remark \ref{rem2.22}, and
$\eta_0,\ldots,\eta_r$ the $\cD_\phi$-orders. By Example \ref{ex1.52},
$j_r-j_r, j_r-j_{r-1},\ldots,j_r-j_0$ are the $(\cD_\phi,(0:1))$-orders.
Then, for each $i$, $j_r-j_{r-i}\geq \eta_i\geq \epsilon_i$ by Corollary
\ref{cor2.21}(1) and Remark \ref{rem2.22}, and the result follows.
     \end{proof}

     \begin{remark}\label{rem2.23} Corollary \ref{cor2.23} was first
noticed by Homma \cite{homma} for $\cD$-orders; see also \cite{garcia2}
and \cite{homma1}.
     \end{remark}

Now we define the so-called ramification divisor of
$\cD$. Let $f_0',\ldots, f_r'$ be any base of $\cD'$ and $x$ any
separating variable of $\F(\cX)|\F$. As before let $P\in\cX$, $t$ a local
parameter at $P$, $\{f_0,\ldots,f_r\}$ a
$(\cD,P)$-base; set $g_\l=t^{v_P(E)}f_\l$. We have a matrix $(a_{ij})\in
GL(r+1,\F)$ such
that
$f_i'=\sum_ja_{ij}f_j$ for each $i$. Proposition \ref{prop2.21} implies
   \begin{align*}
W^{\epsilon_0,\ldots,\epsilon_r}_{f_0',\ldots,f_r';x} 
= &\ \det(a_{ij})W^{\epsilon_0,\ldots,\epsilon_r}_{f_0,\ldots,f_r;x}=\  
\det(a_{ij})t^{-(r+1)v_P(E)}W^{\epsilon_0,\ldots,\epsilon_r 
}_{g_0,\ldots,g_r;x}\\
 = &\ \det(a_{ij})t^{-(r+1)v_P(E)}
(D^1_xt)^{\sum_i\epsilon_i}W^{\cE}_{g_0,\ldots,g_r;t}\, ;
   \end{align*}
i.e.,
    \begin{equation}\label{eq2.23}
W^{\epsilon_0,\ldots,\epsilon_r 
}_{f_0',\ldots,f_r';x}(D_t^1x)^{\sum_i\epsilon_i}t^{(r+1)v_P(E)}=
\det(a_{ij}) W^{\epsilon_0,\ldots,\epsilon_r}_{g_0,\ldots,g_r;t}\, .
     \end{equation}
Thus the divisor 
   $$ 
R=R^{\cD}:= \div(W^{\epsilon_0,\ldots,\epsilon_r}_{f_0',\ldots,f_r';x})+
(\sum_{i=0}^r\epsilon_i)\div(dx)+(r+1)E\ ,
   $$ 
just depends on $\cD$ and locally is given by (\ref{eq2.23}).

    \begin{definition*} $R$ is called the {\em
ramification divisor} of $\cD$. The {\em ramification divisor} of a
morphism $\phi$ is the ramification divisor of $\cD_\phi$.
    \end{definition*} 

    \begin{example}\label{ex2.21} Let $x$ be a separating variable of
$\F(\cX)|\F$ and consider the morphism $\phi=(1:x):\cX\to\P^1(\F)$. Then
$E_\phi=\div_\infty(x)$; moreover, as $\#
x^{-1}(x(P))=\deg(\div_\infty(x))$ for infinitely many $P\in\cX$, the
$\cD_\phi$-orders are 0,1. Then
   $$
R^{\cD_\phi}=\div(dx)+2\div_\infty(x)\, ;
   $$
i.e., it coincides with the ramification divisor $R_x$ of $x$, see Example
\ref{ex1.11}.
    \end{example}
    \begin{lemma} {\rm (Garcia-Voloch \cite[Thm.
1]{garcia-voloch})}\label{lemma2.23} Let $\phi=(f_0:\ldots:f_r)$ be a
morphism associated to $\cD$, and $q'$ a power of ${\rm char}(\F)>0$. Then
$\epsilon_r\geq q'$ if and only if there exist $z_0,\ldots,z_r\in\F(\cX)$,
not all zero, such that 
   $$
z_0^{q'}f_0+\ldots+z_r^{q'}f_r=0\, .
   $$
      \end{lemma}
    \begin{corollary}\label{cor2.24} Let $P\in\cX$. Under the hypothese of
the previous lemma, there exist $i,\l\in \{0,\ldots,r\}$, $i\neq \l$, such
that $j_i(P)\equiv j_\l(P)\pmod{q'}$.
    \end{corollary}
    \begin{proof} We can assume that $f_0,\ldots.f_r$ is a $(\cD,P)$-base.
Now there exist $0\leq i<\l\leq r$ such that
$v_P(z_i^{q'}f_i)=v_P(z_\l^{q'}f_\l)$ and the result follows.
    \end{proof}
\subsection{$\cD$-Weierstrass points}\label{s2.3}
 
Let us keep the notation of the previous subsection. Now we  
study $R$ locally at $P$ via (\ref{eq2.23}); i.e., we study
   $$
v_P(R) = v_P(W^{\epsilon_0,\ldots,\epsilon_r}_{g_0,\ldots,g_r;t})\, .
   $$
We observe that $v_P(R)\ge 0$ since $g_\l$ is regular at $P$ for each
$\l$.

     \begin{theorem}\label{thm2.31} 
\begin{enumerate}
  \item[\rm(1)] $v_P(R)\ge \sum_{i=0}^{r}(j_i(P)-\epsilon_i);$
  \item[\rm(2)] $v_P(R)=\sum_{i=0}^{r}(j_i(P)-\epsilon_i)\
\Leftrightarrow\ \det(\binom{j_\l(P)}{\epsilon_i})\not\equiv 0\pmod{p}.$
\end{enumerate}
      \end{theorem}
      \begin{proof} 
Set $j_i:=j_i(P)$. From the proof of Lemma \ref{lemma2.21} with
$m_i=\epsilon_i$ we have a local expansion of type
   $$
W^{\epsilon_0,\ldots,\epsilon_r}_{g_0,\ldots,g_r;t}=
   C\det(\binom{j_\l}{\epsilon_i})t^{\sum_i(j_i-\epsilon_i)}+\ldots\, ,
   $$ 
with $C\in \F^*$ and the result follows.
      \end{proof} 
We have already observed that $R$ is an effective divisor which also
follows from $j_i(P)\ge \epsilon_i$ (cf. Corollary 
\ref{cor2.21}(1)). Moreover, the following is clear from the theorem.

   \begin{corollary}\label{cor2.31} $v_P(R)=0$ if and only if
$j_i(P)=\epsilon_i$ for each $i$. In particular, for all but finitely many
$P\in\cX$, the $(\cD,P)$-orders equal $\epsilon_0,\ldots,\epsilon_r.$
   \end{corollary}   

   \begin{definition*} The {\em $\cD$-Weierstrass points} of $\cX$ are
those of $\supp(R)$. The {\em $\cD$-weight} of $P$ is $v_P(R)$.
    \end{definition*}

Thus the number of $\cD$-Weierstrass points of $\cX$, counted with their
weighs, equals
   $$
\deg(R)=(\sum_{i=0}^r\epsilon_i)(2g-2)+(r+1)d\, .
   $$
     \begin{lemma} {\rm ($p$-adic criterion)}\label{lemma2.31} Let
$\epsilon$ be a $\cD$-order and let $\mu$ be an integer such that
$\binom{\epsilon}{\mu}\not\equiv 0\pmod{p}$. Then $\mu$ is also a
$\cD$-order. In particular, $0,1,\ldots,\epsilon-1$ are $\cD$-orders
provided that $p>\epsilon.$
     \end{lemma}
     \begin{proof} Let $\l\in\{0,\ldots,r-1\}$ be such that
$\epsilon_\l<\mu\le\epsilon_{\l+1}\le \epsilon$. We apply Corollary
\ref{cor2.22} to a point $P\not\in\supp(R)$; i.e., such that
$j_i(P)=\epsilon_i$ for
each $i$. Let $m_0=\epsilon_0, \ldots,m_\l=\epsilon_\l, m_{\l+1}:=\mu$.
Then the vectors
$(\binom{\epsilon_0}{m_s},\ldots,\binom{\epsilon_r}{m_s})$,
$s=0\ldots,\l+1$, are $\F_p$-linearly independent and the result follows.
     \end{proof}

    \begin{definition*} The curve $\cX$ is called {\em classical with
respect to $\cD$}, or the linear series $\cD$ is called {\em classical},
if the $\cD$-orders are $0,\ldots,r$. A morphism $\phi$ is called
{\em classical} if $\cD_\phi$ is classical.
     \end{definition*}

     \begin{lemma}\label{lemma2.32} Suppose that 
$\prod_{i>\l}\frac{j_i(P)-j_\l(P)}{i-\l}\not\equiv 0\pmod{p}$. Then
   \begin{enumerate}
\item[\rm(1)] $\cD$ is classical$;$
\item[\rm(2)] $v_P(R)=\sum_{i=0}^{r}(j_i(P)-i).$ 
    \end{enumerate}
    \end{lemma}
     \begin{proof} (1) Set $j_i=j_i(P)$. We have 
   $$
\det(\binom{j_i}{\l})=\prod_{i>\l}\frac{j_i-j_\l}{i-\l}\not\equiv 0 
\pmod{p}\, ,
   $$
by hypothesis. Then $\epsilon_i\le i$ by Corollary
\ref{cor2.21}(2); i.e, $\epsilon_i=i$ for each $i$.

(2) Follows from Theorem \ref{thm2.31}(2).
     \end{proof}

In particular, as $j_r(P)\le d=\deg(\cD)$, we obtain:

     \begin{corollary}\label{cor2.32} If $p=0$ or $p>d=\deg(\cD),$ then 
    \begin{enumerate}
\item[\rm(1)] $\cD$ is classical$;$
\item[\rm(2)] For each $P\in \cX$, $v_P(R)=\sum_i(j_i(P)-i).$
    \end{enumerate}
     \end{corollary}
\subsection{$\cD$-osculating spaces}\label{s2.4} 

Assume that $\cD$ is base-point-free, $\cD=g^r_d\cong\P^r(\cD')\subseteq
|E|$.
From
Remark \ref{rem1.41},
   $$
\cD=\{\phi^*(H): \text{$H$ hyperplane in $\P^r$}\}\, ,
   $$ 
where $\phi=(f_0:\ldots:f_r)$, and where $\{f_0,\ldots,f_r\}$ is a
$\F$-base of $\cD'$. Let $P\in\cX$ with $(\cD,P)$-orders $j_0<\ldots<j_r$.
From Lemma \ref{lemma1.22},
  $$
v_P(E)=-{\rm min}\{v_P(f_0),\ldots,v_P(f_r)\}\, .
  $$
For $i=0,\ldots,r-1$, let 
$L_i^{f_0,\ldots,f_r}(P)$ be the intersection of the 
hyperplanes $H$ in $\P^r$ such that $v_P(\phi^*(H))\ge j_{i+1}$. If 
$g_0,\ldots,g_r$ is another base of $\cD'$, there exists $T\in
\aut(\P^r(\F))$ such that $\phi_1:=(g_0:\ldots:g_r)=T\circ\phi$; thus 

    \begin{equation}\label{eq2.41}
L_i^{g_0,\ldots,g_r}(P)=T(L_i^{f_0,\ldots,f_r}(P))\, .
    \end{equation}

We conclude then that $L_i^{f_0,\ldots,f_r}(P)$ is uniquely determinated
by $\cD$ up to projective equivalence. 

    \begin{definition*}
$L_i(P)=L_i^{f_0,\ldots,f_r}(P)$ is called the {\em $i$-th osculating
space} at $P$ (with respect to the base $\{f_0,\ldots,f_r\})$.
    \end{definition*}

Clearly we have:
   $$
L_0(P)\subseteq\ldots\subseteq L_{r-1}(P)\, .
   $$

    \begin{lemma}\label{lemma2.41} $L_i^{f_0,\ldots,f_r}(P)$ is an
$i$-dimensional space generated by the vectors $(D^{j_s}_t\phi')(P)$,
$s=0,\ldots, i,$ where $\phi'=(t^{v_P(E)}f_0:\ldots:t^{v_P(E)}f_r).$
    \end{lemma}
    \begin{proof} From Lemma \ref{lemma1.311} and (\ref{eq2.41}) we can
assume that $f_0,\ldots,f_r$ is a $(\cD,P)$-base. Let 
$H_i$ be the hyperplane corresponding to $X_i=0$, where $X_0,\ldots,X_r$ 
are homogeneous coordinates of $\P^r$. Let $H: \sum_ia_iX_i=0$ be a
hyperplane. Then $v_P(\phi^*(H))\geq j_{i+1}$ if and only if $a_0=\ldots
a_i=0$, since $v_P(t^{v_P(E)}f_\l)=j_\l$ for each $\l$. Thus
   $$
L_i^{f_0,\ldots,f_r}(P)=H_{i+1}\cap\ldots\cap H_r\, ;
   $$ 
i.e., it has dimension $i$. In addition, it is generated by the vectors 
$(D^{j_s}_t\phi')(P)$ by the proof of Scholium \ref{scholium2.21} 
    \end{proof}

From the proof above we obtain:

   \begin{scholium}\label{scholium2.41} $H\supseteq L_i(P)$ if and only if
$v_P(\phi^*(H))\geq j_{i+1}.$
    \end{scholium}

    \begin{remark}\label{rem2.41} If $\cD$ has base points, the
$i$-osculating spaces for $\cD$ are, by definition, those of $\cD^B$.
    \end{remark}

    \begin{definition*} The 1-osculating (resp. $(r-1)$-osculating)  
space at $P$ is called the {\em tangent line} (resp. {\em osculating
hyperplane }) at $P$.
    \end{definition*}

A consequence of Lemma \ref{lemma2.41} is the following.

    \begin{corollary}\label{cor2.41} The osculating hyperplane at $P$
(with respect to the base $\{f_0,\ldots,f_r\}$) is given by the equation
   $$
{\rm det}\begin{pmatrix} X_0& \ldots & X_r \\
                (D^{j_0}_tg_0)(P) & \ldots & (D^{j_0}_tg_r)(P) \\ 
                  \vdots & \vdots & \vdots \\
            (D^{j_{r-1}}_tg_0)(P) & \ldots & (D^{j_{r-1}}_tg_r)(P)
\end{pmatrix}=0\, ,
   $$
where $g_\l:=t^{v_P(E)}f_\l,$ $ \l=0,\ldots,r.$
    \end{corollary}

    \subsection{Weierstrass points; Weierstrass semigroups II}\label{s2.5}
In this sub-section we consider Weierstrass Point Theory for the canonical
linear series $\cK=\cK^\cX$ on the curve $\cX$ of genus $g$. By Remark
\ref{rem1.511} we can assume $g\geq 2$. The special feature in
the canonical case is the existence of a (numerical) semigroup, namely the
Weierstrass semigroup $H(P)$ at $P\in\cX$ (cf. Sect. 1.5) which is closely
related to the $(\cK,P)$-orders. We stress the following.

  \begin{definition*} \begin{enumerate} 
  \item[\rm(1)] The {\em Weierstrass points} of the curve $\cX$ is the set
$\cW=\cW_\cX$ of its $\cK$-Weierstrass points; i.e., $\cW=\supp(R^\cK)$.
The $\cK$-weight of $P$ is called the {\em Weierstrass weight} $\omega_P$
of $P$; i.e., $\omega_P=v_P(R^\cK).$

  \item[\rm(2)] We set $w_P:=\sum_{i=0}^{g-1}(j_i^\cK(P)-i)$; i.e., $w_P$
is the weight of the Weierstrass semigroup $H(P)$ at $P$.

   \item[\rm(3)] The curve $\cX$ is called {\em classical} if it is
classical with respect to the canonical linear series $\cK$.
   \end{enumerate}
   \end{definition*}

In particular, since $\cK$ has dimension $g-1$ and degree $2g-2$, the
number of Weierstrass points $P\in \cW$ counted with their weights
$\omega_P$ equals

   \begin{equation}\label{eq2.51}
\deg(R^\cK)= (\sum_{i=0}^{g-1}\epsilon_i)(2g-2)+g(2g-2)\, ,
   \end{equation}

where $\epsilon_0<\ldots<\epsilon_{g-1}$ are the $\cK$-orders. From
Theorem \ref{thm2.31}(1) we have
   
   $$
\omega_P\geq \sum_{i=0}^{g-1}(j_i^\cK(P)-\epsilon_i)\, .
   $$

In general, $\omega_P>\sum_i(j_i^\cK(P)-\epsilon_i)$ and $\omega_P\neq
w_P$ (see Example \ref{ex2.51}); however, if either $p=0$ or $p>2g-2$,
then the curve is classical and $\omega_P=\sum_i(j_i^\cK(P)-i)=w_P$ by
Corollary \ref{cor2.32}; in this case the curve has $g(g^2-1)$ Weierstrass
points (counted with their weights) by (\ref{eq2.51}).

   \begin{example} {\rm (Hyperelliptic curves)}\label{ex2.51} Let $\cX$
be hyperelliptic with $g^1_2=|\div_\infty(f)|$, $f\in\F(\cX)$ of degree
two. Note that $f$ is a separating variable since $g>0$. We have
$\cK=|(g-1)\div_\infty(f)|$, where $\cK'$ is generated by $1,f,\ldots,
f^{g-1}$. Then $W^{0,1,\ldots,g-1}_{1,f,\ldots,f^{g-1};f}=1$; i.e., $\cX$
is classical. 

The ramification divisor of $\cK$ is thus
  $$
R^\cK= \frac{g(g-1)}{2}\div(df)+g(g-1)\div_\infty(f)\, ,
  $$
so that $R^\cK=\frac{g(g-1)}{2}R_f$ by Example \ref{ex2.21}. Note that $f$
has $\deg(R_f)=2g+2$ ramifications points (counted with multiplicity), and
that $P\in \supp(R_f)$ if and only if $e_P=2$; see Example \ref{ex1.11}.
Therefore the following conditions are equivalent:

  \begin{itemize}
\item $P\in \cW$;
\item $P\in \supp(R_f)$;
\item $e_P=2$;
\item $2\in H(P)$;
\item the $(\cK,P)$-orders are $0,2,\ldots,2g-2$.
  \end{itemize}

If $P\not\in \cW$, then the $(\cK,P)$-orders are $0,1,\ldots,g-1$; i.e.,
$H(P)=\{0,g+1,\ldots\}$. In particular, a hyperellitpic curve has only two
types of Weierstrass semigroups.

If $p=0$ or $p>2$, and $P\in\supp(R_f)$, then $v_P(R_f)=1$ and hence $\cX$
has $2g+2$ Weierstrass points $P$ such that $\omega_P=g(g-1)/2$. In
particular, here we have $\omega_P=\sum_i(j_i^\cK-i)=w_P\ (*)$.

If $p=2$, then $(*)$ is in general not true as the following
example
shows. Let $\cX$ be the non-singular model of the plane curve of equation
  $$
y^2+y=x^{q+1}\, ,
  $$ 
over $\F$ of characteristic two, and where $q=2^a$, $a\geq 2$. Then
$x\in\F(\cX)$ has degree two an so $\cX$ is hyperellitpic. There are two
different points in $\cX$ over each $a\in\F$, since $Y^2+Y=a$ has two
different
solutions. Let $P$ over $x=\infty$. Then $2v_P(y)=-(q+1)e_P$ so that
$e_P=2$; hence there is just one point $P_\infty$ over $x=\infty$;
i.e., $\#\supp(R_x)=1$. In particular, $P_\infty$ is the only Weierstrass
point of $\cX$ and thus its weight is $\omega_P=\deg(R^\cK)=g(g^2-1)>
\sum_i(j_i^\cK(P)-i)=w_P=g(g-1)/2$ because $g>1$ as we see below.

To compute the genus of $\cX$ we use the fact that $P_\infty$ is the only
ramified point for $x$: We have $2g-2=\deg(dx)=v_{P_\infty}(dx)=q-2$ and
so $g=q/2>1$.
   \end{example}

   \begin{lemma}\label{lemma2.51} Let $\cX$ be a classical curve of genus
$g$ such that $\omega_P=w_P$ for each $P$ (e.g. if $p=0$ or $p>2g-2$).
Then
   \begin{enumerate}
\item[\rm(1)] $2g+2\le \#\cW \le g(g^2-1);$
\item[\rm(2)] $\# \cW=2g+2$ if and only if $\cX$ is hyperelliptic$;$
\item[\rm(3)] $\# \cW=g(g^2-1)$ if and only if $\omega_P=1$ for any
$P\in\cX.$
   \end{enumerate}
   \end{lemma}
   \begin{proof} We have $g(g^2-1)=\deg(R^\cK)=\sum_Pw_P\leq \#\cW
g(g-1)/2$ by Corollary \ref{cor1.521}(1). This proves (1). (2) follows
from Corollary \ref{cor1.521}(2)(3) and Example \ref{ex2.51}. (3) is trivial.
   \end{proof}

   \begin{lemma}\label{lemma2.511} Let $(\tilde n_i:i\in\N)$ be the
Weierstrass semigroup at non-Weierstrass points. Then $n_i(P)\le \tilde
n_i$ for each $P$ and each $i.$
   \end{lemma}
   \begin{proof} Let $i$ be the minimum positive integer such that
$n_i(P)>\tilde n_i$. Then $i\ge 2$ and $n_{i-1}(P)\le \tilde n_{i-1}$ so
that $n_{i-1}(P)\leq \tilde n_{i-1}<\tilde n_i<n_i(P)$. Now we have 
$\tilde n_i=\ell_{\tilde n_i-i+1}\geq \tilde \ell_{\tilde n_i-i+1}$ by
Corollary \ref{cor2.21}(1), where $\tilde \ell_1<\tilde\ell_2<\ldots$ are
the gaps at non-Weierstrass points. Since $\ell_{\tilde n_i-i+1}\ge \tilde
n_i+1$ we have a contradiction and the result follows.
   \end{proof}

   \begin{lemma}\label{lemma2.52} The largest $\cK$-order $\epsilon_{g-1}$
is less than $\deg(\cK)=2g-2.$
   \end{lemma}
   \begin{proof} (Garcia \cite[p. 235]{garcia1}) Suppose
$\epsilon_{g-1}=2g-2$. Then for $P\not\in\cW$, $(2g-2)P$ is a canonical
divisor. In particular, $(2g-2)P\sim (2g-2)P_0$ for $P,P_0\not\in\cW\
(*)$. We consider the isogeny $i: D\mapsto (2g-2)D$ on the Jacobian
variety $\cJ$ associated to $\cX$, and the natural map $\cX\to \cJ$,
$P\mapsto [P-P_0]$. Note that $[P-P_0]=[Q-P_0]$ if and only $P=Q$ since
$g>0$. Then $(*)$ says that there are infinitely points in $\cJ$ belonging
to the kernel of $i$, a contradiction since this kernel is finite \cite[p.
62]{mumford}.
   \end{proof}

   \begin{example} (The non-classical curve of genus 3)\label{ex2.52} It
is easy to see that the only semigroups of genus two are
$\{0,3,4,5,\ldots\}$ and $\{0,2,4,5,\ldots\}$. Since a curve of genus two
must have at least one Weierstrass points, then such a curve is
hyperelliptic and hence classical.

Now let $\cX$ be a curve of genus three. We shall show a result due to
Komiya \cite{komiya}: $\cX$ is non-classical if and only if $p=3$ and
$\cX$ is $\F$-isomorphic to the non-singular plane curve of equation
$y^3+y=x^4$. If $\cX$ is non-classical, then $0<p<2g-2=4$ by Corollary
\ref{cor2.32} so that $p=2,3$. We have $\epsilon_0=0, \epsilon_1=1$ and
$\epsilon_2=3$. Then $p=3$ by the $2$-adic criterion. We have $P\in\cW
\Leftrightarrow j_0^\cK(P)=0, j_1^\cK(P)=1, j_2^\cK(P)=4 \Leftrightarrow
H(P)=\{0,3,4,6,\ldots\}$; then $\omega_P=1$ and $\cX$ has $\deg(R^\cK)=28$
Weierstrass points (note that a classical curve of genus 3 has $3\times
(3^2-1)=24$ Weierstrass points counted with their weights). Let
$P_0\in\cW, x,y\in\F(\cX)$ such that $\div_\infty(x)=3P_0$ and
$\div_\infty(y)=4P_0$. We see that $4P_0$ is a canonical divisor and so
$\cK=|4P_0|$. We also see that $x$ is a separating variable of
$\F(\cX)|\F$ so that $W^{0,1,2}_{1,x,y;x}=D^2_xy=0$ as $\epsilon_2>2$. Now
the eleven functions $1,x,y, x^2,xy,y^2,x^3,x^2y,xy^2,x^4,y^3$ belong to
$L(12P_0)$ which has dimension 10. Therefore there is a non-trivial
relation over $\F$ of type
   $$
a_{00}+a_{10}x+a_{01}y+a_{20}x^2+a_{11}xy+a_{02}y^2+a_{30}x^3+a_{21}x^2y+
a_{12}xy^2+a_{40}x^4+a_{03}y^3=0\, .
   $$
Since $v_P(x^iy^j)<12$ for $3i+4j<12$ we must have $a_{40}\neq 0$ and
$a_{03}\neq 0$. In particular we can assume $a_{40}=1$. Next we apply
$D^2_x$ to the equation above; using the fact that $D^2_xy=0$ we find:
   $$
a_{20}+a_{11}D_xy+a_{02}(D_xy)^2+a_{21}(y+2xD_xy)+a_{12}(2xyD_xy+x(D_xy)^2)
=0\, .
   $$ 
Let $v_P(D_xy)=a$. Then the valuation at $P$ of the functions
  $$
1,D_xy,(D_xy)^2,y,xD_xy,xyD_xy,x(D_xy)^2
  $$ 
are respectively 
  $$
0, a,2a,-4,-3+a,-7+a,-3+2a\, ;
  $$ 
we see that they are pairwise different and 
hence $a_{20}=a_{11}=a_{02}=a_{21}=a_{12}=0$; i.e., we have
   $$
a_{00}+a_{10}x+a_{01}y+a_{30}x^3+x^4+a_{03}y^3=0\, .
   $$
By means of $x\mapsto (x-a_{30})$ and $y\mapsto -(a_{03})^{1/3}y$ we can
assume $a_{30}=0$ and $a_{03}=-1$. Now as 
$[\F(\cX):\F(x)]=3$ the above equation is irreducible and hence
$a_{01}\neq 0$ because $x$ is a separating variable. Then by means of
$x\mapsto a_{01}^{3/8}x$ and $y\mapsto -a_{01}^{1/2}y$ we can assume
$a_{01}=1$. So we have an equation of type
   $$
y^3+y=x^4+a_{10}x+a_{00}\, .
   $$ 
Finally let $P_1$ be another Weierstrass point. Then $4P_1\sim 4P_0$ as
both divisor are canonical. So we can choose $y$ such that
$\div(y)=4P_1-4P_0$. Then $4=v_{P_1}(y)=v_{P_1}(x^4+a_{10}x+a_{00})$
implies $a_{00}=a_{10}=0$. 

Conversely if $\cX$ is defines by $y^3+y=x^4$, we have that $\cX$ is a
non-singular plane curve of genus three. Moreover there is just one point
$P_\infty$ over $x=\infty$ and $H(P_\infty)=\{0,3,4,6,\ldots\}$. This
implies that $x$ is a separating variable and we have $D_x^2y=0$; i.e.,
$\cX$ is non-classical.
   \end{example}
Further examples of non-classical linear series can be found in Neeman
\cite{neeman}. Finally we mention that Weierstrass Point Theory on schemes
was considered by Laksov and Thorup \cite{laksov-thorup}; see the
introduction there for further references.
   \section{Frobenius orders}\label{s3}

Let $\cX$ be a curve defined over $\fq$, a finite field with $q$ elements;
i.e., $\cX$ is a
curve over the algebraic closure $\bar\fq$ of $\fq$, equipped with the
action of the Frobenius morphism $\fro$ relative to $\fq$. Let $\cD\cong
\P(\cD')\subseteq |E|$ be a base-point-free $g^r_d$ on $\cX$. Assume that
$\cD$ is also defined over $\fq$; i.e., for any $D=\sum_Pn_PP\in \cD$,
$(\fro)_*(D):=\sum_Pn_P\fro(P)=D$. Let $\phi=(f_0:\ldots:f_r)$ be a
morphism over $\fq$ associated to $\cD$; i.e., its coordinates belong to
$\fq(\cX)$ and they form a $\fq$-base of $\cD'$.

The starting point of St\"ohr-Voloch's approach to the Hasse-Weil bound is
to look at points $P$ of $\cX$ such that $\phi(\fro(P))$ belongs to the
osculating hyperplane $L^{f_0,\ldots,f_r}_{r-1}(P)$ at $P$. Then Corollary
\ref{cor2.41} leads to the consideration of rational functions of type
   $$
V^{\l_0,\ldots,\l_{r-1}}_{f_0,\ldots,f_r;x}:=
{\rm det}\begin{pmatrix} f_0\circ\fro & \ldots & f_r\circ\fro \\
                D^{\l_0}_xf_0 & \ldots & D^{\l_0}_xf_r \\ 
                  \vdots & \vdots & \vdots \\
            D^{\l_{r-1}}_xf_0 & \ldots & D^{\l_{r-1}}_xf_r
\end{pmatrix}\, ,
   $$
where $x$ is a separating variable of $\bar\fq(\cX)|\bar\fq$. We set
   $$
\cB(f_0,\ldots,f_r;x):=\{(m_0,\ldots,m_{r-1})\in \N^r_0:
m_0<\ldots<m_{r-1};\ V^{m_0,\ldots,m_{r-1}}_{f_0,\ldots,f_r;x}\neq 0\}\, .
   $$

          \begin{lemma}\label{lemma3.1} Let $(m_0,\ldots,m_r)\in
\cA(f_0,\ldots,f_r;x)$ with $m_0=0$. Then there exists $0<I\le r$ such
that $(m_0,\ldots,m_{I-1},m_{I+1},\ldots,m_r)\in \cB(f_0,\ldots,f_r;x)$.
          \end{lemma}
          \begin{proof} Let $I$ be the smallest integer such that
$\phi\circ\fro:=(f_0\circ\fro,\ldots,f_r\circ\fro)$ is a $\F(\cX)$-linear
combination of $D^{m_0}_x\phi,\ldots,D^{m_I}_x\phi$. Since
$f_0,\ldots,f_r$ is a $\fq$-base of $\cD'$, then $I>0$ and the result
follows.
           \end{proof}

Since the $\cD$-order sequence $(\epsilon_0,\ldots,\epsilon_r)$ 
belongs to $\cA(f_0,\ldots,f_r;x)$ (cf. Proposition \ref{prop2.21}), 
$\cB(f_0,\ldots,f_r;x)\neq\emptyset$. Let
   $$
\cV:=(\nu_0,\ldots,\nu_{r-1})
   $$
be the minimum (in the lexicographic order) of $\cB(f_0,\ldots,f_r;x)$. 

       \begin{lemma}\label{lemma3.2} 
   \begin{enumerate}
\item[\rm(1)] $\nu_0=0;$
\item[\rm(2)] For $i=1,\ldots,r-1,$
   $$
\nu_i=\min\{s>\nu_{i-1}: \text{$\phi\circ\fro, D^{\nu_0}_x\phi,\ldots,
D^{\nu_{i-1}}_x\phi, D^s_x\phi$ are $\bar\fq(\cX)$-l.i}\}\, ;
   $$
\item[\rm(3)] Let $(m_0,\ldots,m_{r-1})\in
\cB(f_0,\ldots,f_r;x)$. Then $\nu_i\le m_i$ for each $i.$
   \end{enumerate}
       \end{lemma}
       \begin{proof} Similar to the proofs of Lemma 
\ref{lemma2.22} and Corollary \ref{cor2.21}(1).
        \end{proof}

       \begin{corollary}\label{cor3.1} There exists $0<I\le r$ such that 
   $$
\nu_i=\begin{cases}
\epsilon_i & \text{if $i<I$} ,\\
\epsilon_{i+1} & \text{if $i\ge I$} .\end{cases}
   $$
           \end{corollary}
           \begin{proof} From Proposition \ref{prop2.21}(3) and Lemma
\ref{lemma3.1}, there exists $0<I\le r$ such that 
$(\epsilon_0,\ldots,\epsilon_{I-1},\epsilon_{I+1},\ldots,\epsilon_r)\in
\cB(f_0,\ldots,f_r;x)$. Hence from Lemma \ref{lemma3.2}, $\nu_i\le
\epsilon_i$ for $i<I$ and $\nu_i\le \epsilon_{i+1}$ for $i\ge I$. Since
$D^{\nu_0}_x\phi,\ldots, D^{\nu_{I-1}}_x\phi$ are $\F(\cX)$-l.i, from 
Lemma
\ref{lemma2.22}(3) follows that $\epsilon_i\le \nu_i$ for
$i=0,\ldots,I-1$; thus 
$\nu_i=\epsilon_i$ for $i=0,\ldots,I-1$. The same argument yields  
$\epsilon_I\le \nu_I$; in fact, $\epsilon_I<\nu_I$ by the definition of
$I$ in the proof of Lemma \ref{lemma3.1}. 
Suppose that $\nu_I<\epsilon_{I+1}$. Then by Lemma \ref{lemma2.22}(3) the
vectors $D^{\nu_0}_x\phi,\ldots,D^{\nu_{I-1}}_x\phi,
D^{\epsilon_I}_x\phi,D^{\nu_I}_x\phi$ would be linearly dependent over
$\F(\cX)$ so that $D^{\nu_I}\in \langle D^{\nu_0}_x\phi,\ldots,
D^{\nu_{I-1}}_x\phi, D^{\epsilon_I}_x\phi\rangle $. This is a
contradiction because $\phi\circ\fro,
D^{\nu_0}_x\phi,\ldots,D^{\nu_{I-1}}_x\phi,D^{\nu_I}_x\phi$ are
$\bar\fq(\cX)$-linearly independent. A similar argument shows that
$\nu_i=\epsilon_{i+1}$ if $i>I$.
           \end{proof} 

We remark the following computation regarding change of basis. Let
$g_i=\sum a_{ij}f_j$ with $(a_{ij})\in M_{r+1}(\bar\fq)$. Then

     \begin{equation}\label{eq3.1}   
{\rm det}\begin{pmatrix} \tilde g_0 & \ldots & \tilde g_r \\
                D^{\l_0}_xg_0 & \ldots & D^{\l_0}_xg_r \\ 
                  \vdots & \vdots & \vdots \\
            D^{\l_{r-1}}_xg_0 & \ldots & D^{\l_{r-1}}_xg_r
   \end{pmatrix} = 
{\rm det}(a_{ij})V^{\l_0,\ldots,\l_{r-1}}_{f_0,\ldots,f_r;x}\, ,
   \end{equation} 

where $\tilde g_j=\sum_ia_{ij}f_i\circ \fro$. The following is analogous
to Proposition \ref{prop2.21}.

           \begin{proposition}\label{prop3.1} 
\begin{enumerate}
   \item[\rm(1)] If $g_i=\sum_ja_{ij}f_j$ with
$(a_{ij})\in M_{r+1}(\fq)$, then
   $$
V^{\nu_0,\ldots,\nu_{r-1}}_{g_0,\ldots,g_r;x}
=\det((a_{ij}))
V^{\nu_0,\ldots,\nu_{r-1}}_{f_0,\ldots,f_r;x}\, ;
   $$
   \item[\rm(2)] If $f\in \bar\fq(\cX)$, then
   $$
V^{\nu_0,\ldots,\nu_{r-1}}_{ff_0,\ldots,ff_r;x}=f^{q+r}
V^{\nu_0,\ldots,\nu_{r-1}}_{f_0,\ldots,f_r;x}\, ;
   $$
   \item[\rm(3)] Let $y$ be any separating variable of
$\bar\fq(\cX)|\bar\fq$. Then
   $$
V^{\nu_0,\ldots,\nu_{r-1}}_{f_0,\ldots,f_r;y}=(D^1_yx)^{\sum_i
\nu_i}
V^{\nu_0,\ldots,\nu_{r-1}}_{f_0,\ldots,f_r;x}\, .
   $$
   \end{enumerate}
   \end{proposition}
     \begin{proof} (1) follows from (\ref{eq3.1}) taking into
consideration that $a_{ij}^q=a_{ij}$. (2) and (3) follow as in 
Proposition \ref{prop2.21}. 
      \end{proof}

Now we show that $\cV$ just depend on $\cD$ and $q$. Let
$\{f_0',\ldots,f_r'\}\subseteq \fq(\cX$ be another $\fq$-base of $\cD'$
and $y$ another separating variable of $\bar\fq(\cX)|\bar\fq$. From part
(1) above, $\cV$ is the minimum for $\cB(f_0',\ldots,f_r';x)$ and from
part (3) it is also the minimum for $\cB(f_0',\ldots,f_r';y)$.

   \begin{definition*} $\cV=(\nu_0,\ldots,\nu_{r-1})$ is called the {\em
$\fq$-Frobenius orders} of $\cD$. If $\nu_i=i$ for each $i$, $\cD$ is
called {\em $\fq$-Frobenius classical}.
   \end{definition*}

Now let $P\in\cX$. We have that $v_P(E)=-\min(v_P(f_0),\ldots,v_P(f_r))$
because $\cD$ is base-point-free, cf. Lemma \ref{lemma1.22}. In
addition, the rational functions $g_i:=t^{v_P(E)}f_i$ are regular at $P$
for each $i$, where $t$ is a local parameter at $P$. Let 
$\{f_0',\ldots,f_r'\}$ and $y$ be as above. Let $f_i'=\sum_j
a_{ij}f_j$, $a_{ij}\in \fq$. Applying Proposition \ref{prop3.1} we have
    \begin{equation*}
    \begin{split}
V^{\nu_0,\ldots,\nu_{r-1}}_{f_0',\ldots,f_r';y} & = \det(a_{ij})
V^{\nu_0,\ldots,\nu_{r-1}}_{f_0,\ldots,f_r;y}\\
{} & =
\det(a_{ij})(D^1_yt)^{\sum_i\nu_i}
V^{\nu_0,\ldots,\nu_{r-1}}_{f_0,\ldots,f_r;t}\\
{} & = \det(a_{ij})(D^1_yt)^{\sum_i\nu_i}
t^{-(q+r)v_P(E)}V^{\nu_0,\ldots,\nu_{r-1}}_{g_0,\ldots,g_r;t}\, ;
    \end{split}
    \end{equation*}
i.e., 

   \begin{equation}\label{eq3.2}
V^{\nu_0,\ldots,\nu_{r-1}}_{f_0',\ldots,f_r';y}(\frac{dy}{dt})^{\sum_i\nu_i}
t^{(q+r)v_P(E)}=
\det(a_{ij})V^{\nu_0,\ldots,\nu_{r-1}}_{g_0,\ldots,g_r;t}\, .
   \end{equation}

Therefore the divisor 
   $$
S=S^{\cD,q}:=\div(V^{\nu_0,\ldots,\nu_{r-1}}_{f_0',\ldots,f_r';y})+
(\sum_{i=0}^{r-1}\nu_i)\div(dy)+(q+r)E\, ,
   $$
just depend on $\cD$ and $q$ and locally at $P$ is given by (\ref{eq3.2}).

   \begin{definition*} $S$ is called the {\em $\fq$-Frobenius divisor} of
$\cD$.
    \end{definition*}

The divisor $S$ is effective because, as we already noticed, each $g_\ell$
is regular at $P$. Note that
   $$
\deg(S)=(\sum_{i=0}^{r-1}\nu_i)(2g-2)+(q+r)d\, .
   $$
Next we study $v_P(S)$ by means of (\ref{eq3.2}); i.e. we study
   $$
v_P(S)=v_P(V^{\nu_0,\ldots,\nu_{r-1}}_{g_0,\ldots,g_r;t})\, .
   $$
We consider two cases according as $P$ is $\fq$-rational or not.

{\bf Case I:} $P\in\cX(\fq)$. Here we can assume that $f_0,\ldots,f_r$ is
a $(\cD,P)$-base; i.e, $v_P(g_\l)=j_\l$ for $\l=0,\ldots,r$. By
Proposition \ref{prop3.1}(2)
   $$
v_P(S)=v_P(g_0^{q+r}V^{\nu_0,\ldots,\nu_{r-1}}_{h_0,\ldots,h_r;t})=
v_P(V^{\nu_0,\ldots,\nu_{r-1}}_{h_0,\ldots,h_r;t})\, ,
   $$ 
where $h_\l:=g_\l/g_0$. Note that $h_0=1$ and that
$v_p(h_\l)=j_\l$. In particular,
   \begin{equation}\label{eq3.21}
V^{\nu_0,\ldots,\nu_{r-1}}_{h_0,\ldots,h_{r-1};t}=
{\rm det}\begin{pmatrix} h_1-h_1^q & \ldots & h_r-h_r^q \\
                D^{\nu_1}_th_1 & \ldots & D^{\nu_1}_th_r \\ 
                  \vdots & \vdots & \vdots \\
            D^{\nu_{r-1}}_th_1 & \ldots & D^{\nu_{r-1}}_th_r
   \end{pmatrix}\, ,
   \end{equation}
and we can made similar computations as in the proof of Lemma
\ref{lemma2.21}: Expand $h_\ell$ at $P$, $h_\ell=\sum_{s=j_\l}^\infty
c_s^\l t^s$, set $C:=\prod_{\l=1}^rc^\l_{j_\l}$; then 

   \begin{equation}\label{eq3.3}
V^{\nu_0,\ldots,\nu_{r-1}}_{h_0,\ldots,h_r;t}=C{\rm
det}(\binom{j_\l}{\nu_i})t^{\sum_{i=i}^{r-1}(j_i-\nu_{i-1})}+\ldots\, ,
   \end{equation}

where $i=0,\ldots,r-1;\l=1,\ldots,r$ in the matrix above involving the
binomial operator. Now $v_P(S)$ can be estimated via this local expansion.

{\bf Case II:} $P\not\in\cX(\fq)$. Let $h_0,\ldots,h_r$ be a
$(\cD,P)$-base. Then there exists $(a_{ij})\in M_{r+1}(\bar\fq)$ 
such that $h_i':=t^{v_P(E)}h_i=\sum_{j}a_{ij}g_j$. Then from (\ref{eq3.1})
   $$
v_P(S)=v_P(\sum_{i=0}^r(-1)^i\tilde h_i'd_i)\, ,
   $$
where the $d_i$'s are the determinants obtained by Cramer's rule. Clearly
$v_P(\tilde h_i')\ge 0$ and so
   $$
v_P(S)\geq {\rm min}\{v_P(d_0),\ldots,v_P(d_r)\}\, .
   $$
Once again we can expand each $d_i$ at $P$ as in the proof of Lemma
\ref{lemma2.21}: Let
$M:=(\binom{j_\l}{\nu_k})_{k=0,\ldots,r-1;\l=0,\ldots,r}$ and let $M_i$ be 
the matrix obtained from $M$ by deleting the $i$th column. Then
   \begin{equation}\label{eq3.4}
d_i=C_i{\rm
det}(M_i)t^{\sum_{k=0}^rj_k-j_i-\sum_{k=0}^{r-1}\nu_k}+\ldots\, ,
   \end{equation}
where $C_i\in\bar\fq^*$. Thus (\ref{eq3.3}) and (\ref{eq3.4}) imply the
following.

    \begin{proposition}\label{prop3.2} \begin{enumerate} \item[\rm(1)] For
$P\in \cX(\fq),$ $v_P(S)\geq\sum_{i=1}^{r}(j_i(P)-\nu_{i-1})$; equality
holds if and only if
$\det(\binom{j_\l(P)}{\nu_i})_{i=0,\ldots,r-1;\l=1,\ldots,r}\not\equiv
0\pmod{p};$

\item[\rm(2)] For $P\not\in\cX(\fq),$ $v_P(S)\geq
\sum_{i=1}^{r-1}(j_i(P)-\nu_i);$ if
$\det(\binom{j_\l(P)}{\nu_i})_{i,\l=0,\ldots,r-1}\equiv 0\pmod{p},$ then
the stric inequality holds$.$
   \end{enumerate}
   \end{proposition}

   \begin{proposition}\label{prop3.3} Let $\nu$ be a $\fq$-Frobenius order
such that $\nu<q.$ Let $\mu$ an integer such that
$\binom{\nu}{\mu}\not\equiv 0\pmod{p}.$ Then $\mu$ is also an
$\fq$-Frobenius order. In particular, if $\nu_i<p$ then
$(\nu_0,\ldots,\nu_i)=(0,\ldots,i).$
    \end{proposition}
   \begin{proof} Let $\nu=\nu_i$. For $j\le i$, we have
$D^{\nu_j}_t(f^q)=0$ by Remark \ref{rem2.15}. So $\nu_0,\ldots,\nu_i$ are
the first $i+1$ orders of the morphism $(h_1-h_1^q:\ldots:h_r-h^q)$, where
$h_1,\ldots,h_r$ are as in (\ref{eq3.21}). Then the resul follows from the
$p$-adic criterion (Lemma \ref{lemma2.31}).
    \end{proof}

Next we study relations between the $\fq$-Frobenius orders and
$(\cD,P)$-orders at $\fq$-rational points $P$.

   \begin{proposition}\label{prop3.4} Let $P\in\cX(\fq)$ and 
$m_0<\ldots<m_{r-1}$ be a sequence of non-negative integers such that
$\det(\binom{j_\l(P)-j_1(P)}{m_i})_{i=0,\ldots,r-1;\l=1,\ldots,r}\not\equiv
0\pmod{p}.$ Then $\nu_i\le m_i$ for each $i.$
   \end{proposition}   
   \begin{proof} Set $j_i=j_i(P)$ and let
$\phi:=(1:x^{j_2-j_1}:\ldots:x^{j_r-j_1})$, where $x$ is a separating
variable of $\bar\fq(\cX)|\bar\fq$. Let $\eta_0<\ldots<\eta_{r-1}$ be the
orders of $\phi$. Then $\eta_i\le m_i$ for each $i$ by (\ref{eq2.22}),
hypothesis and Corollary \ref{cor2.21}(1). Then, as
$\phi=(x^{j_1}:\ldots:x^{j_r})$, $\det((\binom{j_i}{\eta_\l})\not\equiv
0\pmod{p}$, and the result follows from (\ref{eq3.3}).
   \end{proof}

   \begin{remark}\label{rem3.1} From the proof above follows that the best
election of the $m_i$'s in Proposition \ref{prop3.4} are the orders of the
morphism $\phi=(x^{j_1(P)}:\ldots:x^{j_r(P)})$.
   \end{remark}

   \begin{corollary}\label{cor3.2} Let $P\in\cX(\fq)$.
\begin{enumerate}
   \item[\rm(1)] $\nu_i\leq j_{i+1}(P)-j_1(P)$ for $i=0,\ldots,r-1,$ and
so $v_P(S)\geq rj_1(P);$

   \item[\rm(2)] Suppose $a:=\prod_{1\le i<\l\le
r}(j_\l(P)-j_i(P))/(\l-i)\not\equiv 0\pmod{p}.$ Then $\cD$ is
$\fq$-Frobenius classical and $v_P(S)=r+\sum_{i=1}^r(j_i(P)-i).$
   \end{enumerate}
   \end{corollary}
   \begin{proof} Note that
$a=\det(\binom{j_\l(P)}{i})_{i=0,\ldots,r-1;\l=1,\ldots,r}$. Then 
(1) (resp. (2)) follows from Proposition \ref{prop3.4} with
$m_i=j_i(P)-j_1(P)$ (resp. from the
proof of Proposition \ref{prop3.4} with $m_i=i$, and
Proposition \ref{prop3.2}(1)). 
   \end{proof}
   
   \begin{remark}\label{rem3.2} The criterion of Corollary \ref{cor3.2}(2)
is satisfied if $j_\l(P)-j_i(P)\not\equiv 0\pmod{p}$ for $1\leq i<\l\leq
r$. In particular, the criterion is satisfied if $p\geq j_r(P)$.
   \end{remark}

   \begin{corollary}\label{cor3.3} \begin{enumerate}
   \item[\rm(1)] If $P\in\cX(\fq)$ and
$\det(\binom{j_\l(P)-j_1(P)}{\epsilon_j})
_{j=0,\ldots,r-1;\l=1,\ldots,r}\not\equiv 0\pmod{p},$ then
$\nu_i=\epsilon_i$
for $i=0,\ldots,r-1;$
   \item[\rm(2)] If $\cD$ is not $\fq$-Frobenius classical$,$ then
$j_r(P)>r$ for any $P\in\cX(\fq);$
   \item[\rm(3)] If $(\nu_0,\ldots,\nu_{r-1})\neq
(\epsilon_0,\ldots,\epsilon_{r-1}),$ then $\cX(\fq)\subseteq \supp(R).$
    \end{enumerate}
   \end{corollary}
    \begin{proof} (1) follows from Proposition \ref{prop3.4} with
$m_i=\epsilon_i$.

(2) If there exists $P\in\cX(\fq)$ such that $j_r(P)=r$, then $\nu_i=i$
for each $i$ by Corollary \ref{cor3.2}(1).

(3) Suppose that there exists $P\in\cX(\fq)\setminus \supp(R)$. Then
$j_i(P)=\epsilon_i$ for each $i$ and hence $\nu_i\leq
\epsilon_{i+1}-\epsilon_1$ by Corollary \ref{cor3.2}(1); i.e.
$\nu_i=\epsilon_i$ for each $i$, a contradiction.
   \end{proof}

   \begin{remark}\label{rem3.3} If we choose $i$ such that
$\cX(\fqi)\not\subseteq \supp(R)$, then from Corollary \ref{cor3.3}(3) we
see that the $\fqi$-order sequence of $\cD$ coincide with
$(\epsilon_0,\ldots,\epsilon_{r-1})$.
   \end{remark}

   \begin{theorem}\label{thm3.1} Let $\cX$ be a curve defined over $\fq$
that admits a base-point-free linear series $\cD=g^r_d$ defined over
$\fq$. Let $\nu_0<\ldots<\nu_{r-1}$ be the $\fq$-Frobenius orders of
$\cD$. Then
   $$
\#\cX(\fq)\le \frac{\sum_{i=0}^{r-1}\nu_i(2g-2)+(q+r)d}{r}\, .
   $$
   \end{theorem}
   \begin{proof} Let $S$ be the $\fq$-Frobenius
divisor of $\cD$. Then
$v_P(S)\ge r$ for each $P\in\cX(\fq)$ by Corollary \ref{cor3.2}(1), and so 
$\#\cX(\fq)\leq \deg(S)/r$.
   \end{proof}
   \begin{example} (The Hermitian curve over $\mathbf F_9$)\label{ex3.1}
We are looking for a curve $\cX$ of genus 3 defined over $\fq$ such that
$\#\cX(\fq)>2q+8$. Let $\epsilon_0=0<\epsilon_1=1<\epsilon_2$ (resp.
$\nu_0=0<\nu_1$) be the canonical orders (resp. canonical $\fq$-orders).
   \begin{claim*} $\cX$ is non-classical; i.e., $\epsilon_2>2.$
   \end{claim*}
Indeed, if $\epsilon_2=2$, then $\nu_1\le 2$ by Corollary \ref{cor3.1} and
Theorem \ref{thm3.1} gives $\#\cX(\fq)\le 2q+8$.

Therefore from Example \ref{ex2.52} we conclude that $q$ is a power of
three, $\epsilon_2=3$, and that $\cX$ is given by $y^3+a_{01}y=x^4$, with
$a_{01}\in\bar\fq$ (notice that the change of coordinates involving
$a_{01}$ in Example \ref{ex2.52} is not defined over $\fq$). Moreover,
the proof above also shows that $\nu_1>1$; i.e $\nu_1=3$. 
   \begin{claim*} $q=9$ and $\cX$ is $\mathbf F_9$-isomorphism to the
Hermitian curve 
$y^3+y=x^4.$ In addition, $\cX(\mathbf F_9)=\cW$ (so that $\#\cX(\mathbf
F_9)=28>2\times 9+8$)$.$
   \end{claim*}

Let $x$ and $y$ be as in Example \ref{ex2.52}. Then $V^{0,1}_{1,x,y;x}=0$
or equivalently $y-y^q=(x-x^q)D_xy\ (*)$. Then taking valuation at $P$ we
have $-4q=-3q-9$ so that $q=9$. Moreover from $(*)$ and the equation
defining $\cX$ we have $(1-a_{01}^3)y^3+(a_{10}-1)y^9=0$ so that
$a_{01}=1$. That $\cX(\mathbf F_9)\subseteq \cW$ follows from Corollary
\ref{cor3.3}(3) and equality holds since $\# \cX(\mathbf F_9)=28$ (see
Sect. \ref{s4.2}).

Finally, observe that $\#\cX(\mathbf F_9)$ attains the bound in Theorem
\ref{thm3.1}.
   \end{example}
   \begin{example}\label{ex3.2} (The Hermitian curve, I) Let $\ell$ be a
power
of a prime and $\cH$ the plane curve of equation
   \begin{equation}\label{eq3.5}
Y^\l Z+YZ^\l=X^{\l+1}\, .
   \end{equation}
It is easy to see that $\cH$ is non-singular so that it has genus
$g=\l(\l-1)/2$ by Remark \ref{rem1.22}. 
   \begin{claim*} $\#\cH(\fls)=\ell^3+1.$
   \end{claim*}
Indeed, we have $\cH\cap (Z=0)=\{(0:1:0)\}$; in $Z\neq 0$ we look for
points $(x:y:1)$ such that $y^\l+y=x^{\l+1}$. It follows that
$x\in\fls\Rightarrow y\in\fls$ and since $Y^\l+Y=x^{\l+1}$ has $\l$
different solutions for $Y$ we conclude that there are $\l^3$ such
$(x:y:1)$ points.

Now over $x:=X/Z=\infty$ there is just one point say $P_\infty$ such that
$H(P_\infty)\subseteq\langle \l,\l+1\rangle$. Since $\# (\N\setminus
\langle\l,\l+1\rangle)=\l(\l-1)/2=g$, $H(P_\infty)=\langle\l,\l+1\rangle$.
Next we consider $\cD:=|(\l+1)P_\infty|$ which is a $g^2_{\l+1}$
base-point-free on $\cH$. Since $L((\l+1)P_\infty)=\langle 1,x,y\rangle$,
where $y^\l+y=x^{\l+1}$ we see that $\cD$ is just the linear series cut
out by lines on $\cH$. Let $\epsilon_0=0,\epsilon_1=1,\epsilon_2$ (resp.
$\nu_0=0, \nu_1\in\{1,\epsilon_2\}$) denote the $\cD$-orders (resp.
$\fls$-Frobenius orders) of $\cH$.
   \begin{claim*} 
   \begin{enumerate}
\item[\rm(1)] $\epsilon_2=\nu_1=\l;$
\item[\rm(2)] $j_2(P)=\l+1$ if $P\in\cH(\fls);$ $j_2(P)=\l$ otherwise$.$
   \end{enumerate}
   \end{claim*} 
In fact, $2\#\cH(\fls)\le \nu_1(2g-2)+(\l^2+2)(\l+1)$ by Theorem
\ref{thm3.1} so that $\nu_1\geq \l$. Then $\l\leq \nu_1=\epsilon_2\le
\l+1$ and so $\l=\nu_1=\epsilon_2$ by Lemma \ref{lemma2.31} ($p$-adic
criterion). That $j_2(P)=\l+1$ whenever $P\in\cH(\fls)$ follows from
Corollary \ref{cor3.2}(1) and part (1). In particular for such points $P$,
$v_P(R)=1$. Now we have $\deg(R^\cD)=\l^3+1$ and therefore $j_2(P)=\l$ for
$P\not\cX(\fls)$.

We can write a direct proof of part (2) as follows. Let $a,b\in\bar\fl$
such that $b^\l+b=a^{\l+1}$. It is easy to see that $(x-a)$ is a local
parameter at $(a:b:1)\in\cH$ so that
$(y-b)=a^\l(x-a)+(a-a^\l)(x-a)^\l+(x-a)^{\l+1}+\ldots$. Let
   $$
f:=(y-b)-a^\l(x-a)\, .
   $$
Then
   $$
\div(f)=\l(a:b:1)+(a^{\l^2}:b^{\l^2}:1)-(\l+1)P_\infty
   $$
and part (2) follows.
   \end{example}
Further arithmetical and geometrical properties of Frobenius orders can be
read in Garcia-Homma \cite{garcia-homma}. From that paper we mention
the following.
    \begin{lemma} {\rm (\cite[Cor. 3]{garcia-homma})}\label{lemma3.3} Let
$\cV=\cE\setminus\{\epsilon_I\}$ and suppose that $I<r$. Then ${\rm
char}(\fq)$ divides $\epsilon_{I+1}$.
    \end{lemma}
   \section{Optimal curves}\label{s4}

Let $\cX$ be a curve defined over $\fq$ of genus $g$. To
study quantitative results on the number of $\fq$-rational points of
$\cX$ it is convenient to form a formal power series, the so-called 
{\em Zeta Function} of $\cX$ relative to $\fq$:
  $$
Z_{\cX,q}(t):= {\rm exp}(\sum_{i=1}^{\infty}\frac{\#\cX(\fqi)}{i} 
  t^i)\, .
  $$
By the Riemann-Roch theorem there exists a polynomial $P(t)$ of degree
$2g$ with integer coefficients, such that (see e.g.
\cite[Thm. 3.2]{moreno}, \cite[Thm. V.1.15]{sti}) 
   
   \begin{equation}\label{eq4.1}
Z_{\cX,q}(t)=\frac{P(t)}{(1-t)(1-qt)}\, .
   \end{equation}
   
   \begin{remark} {\rm (\cite[Thm. V.1.15]{sti})}\label{rem4.1}
   \begin{enumerate}
  \item[\rm(i)] Let $P(t)=\sum_{i=0}^{2g}a_it^i$. Then $a_0=1$,
$a_{2g}=q$,
and $a_{2g-i}=q^{g-i}a_i$ for $i=0,\ldots,g$. 

  \item[\rm(ii)] Set 
   $$
h(t)=h_{\cX,q}(t):=t^{2g}P(t^{-1})\, ;
   $$ 
then the $2g$ roots (counted with multiplicity)
$\alpha_1,\ldots,\alpha_{2g}$ of $h(t)$ can be arranged in such a way that
$\alpha_j\alpha_{g+j}=q$ for $j=1,\ldots,g$. Note that
$a_1=-\sum_{j=1}^{2g}\alpha_j$.
   \end{enumerate}
   \end{remark}
   
Now (\ref{eq4.1}) implies $\#\cX(\fq)=q+1+a_1$ and hence that
   $$
\#\cX(\fq)=q+1-\sum_{j=1}^{2g}\alpha_j\, ,
   $$
by Remark \ref{rem4.1}(ii). Furthermore \cite[Cor. V.1.16]{sti},
   $$
\#\cX(\fqi)=q^i+1-\sum_{j=1}^{2g}\alpha_j^i\, .
   $$
By analogy with the Riemann
hypothesis E. Artin conjectured that the absolute value of each $\alpha_i$
equals $\sqrt q$. This result was showed by Hasse for $g=1$ and for A.
Weil for arbitrary $g$ \cite{weil} (see also \cite[Cor. 2.14]{sv},
\cite{moreno}, \cite[Thm. V.2.3]{sti}). In particular, we obtain the 
Hasse-Weil bound on the number of $\fq$-rational points of $\cX$, namely
   $$
|\#\cX(\fq)-(q+1)|\le 2\sq g\, .
   $$
If $\cX$ attains the upper bound above, it is called {\em
$\fq$-maximal}; in this case $q$ must be a square.

   \begin{lemma}\label{lemma4.1} Let $q=\ell^2$. The following statements
are equivalent: 
   \begin{enumerate}
\item[\rm(1)] $\cX$ is $\fls$-maximal$;$
\item[\rm(2)] $\alpha_i=-\ell$ for $i=1,\ldots, 2g;$
\item[\rm(3)] $h_{\cX,\ell^2}(t)=(t+\ell)^{2g}.$
   \end{enumerate}
If any of these conditions hold and $\cX$ is defined over $\fl$, then
   $$
\#\cX(\mathbf F_{\ell^i})=
   \begin{cases}
\ell^i+1 & \text{if $i\equiv 1\pmod{2}$,}\\
\ell^i+1+2\sqrt{\ell^i}g & \text{if $i\equiv 2\pmod{4}$,}\\
\ell^i+1-2\sqrt{\ell^i}g & \text{if $i\equiv 0\pmod{4}$.}
   \end{cases}
   $$
  \end{lemma}
  \begin{proof} $\cX$ is $\fls$-maximal if and only if
$\sum_{i=1}^{2g}\alpha_i=\sum_{i=1}^g(\alpha_i+\bar\alpha_i)=-2\l g$. By
the Riemann-hypothesis, this is the case if and only if $\alpha_i=-\l$ for
each $i$ and the equivalences follow. Now we show the formulae on the
number
of rational points. Let 
$\#\cX(\fl)=\l+1-\sum_{j=1}^{2g}\beta_j$. Then $\beta_j^2=-\l$ for each
$j$
so that $\beta_j^i+\bar\beta_j^i=0$ for $i\equiv 1\pmod{2}$; i.e.,
$\#\cX(\mathbf F_{\ell^i})=\l^i+1$. If $i\equiv 2\pmod 4$,
$\beta_j^i=-\sqrt{\ell^i}$ 
and follows the formula for such $i$'s. Finally, if $i\equiv 0\pmod{4}$,
$\beta_j=\sqrt{\ell^i}$ and the proof is complete.
  \end{proof}
  \begin{corollary} {\rm (Ihara \cite{ihara})}\label{cor4.1} If $\cX$ is
$\fls$-maximal, then $g\leq \ell(\ell-1)/2.$
  \end{corollary}
  \begin{proof} We have $\cX(\fls)\subseteq \cX(\mathbf F_{\ell^4})$. Then
from the lemma above, $\ell^2+1+2\ell g\leq \ell^4+1-2\ell^2g,$ and the
result follows.
   \end{proof} 
   \begin{example}\label{ex4.0} (The Hermitian curve, II) The curve $\cH$
in Example \ref{ex3.2} has genus $\l(\l-1)/2$ and
$\l^3+1$ $\fls$-rational points. Hence it is $\fls$-maximal and attains
the bound in Corollary \ref{cor4.1}.

This curve is called {\em the Hermitian curve} and it is the most fancy
example of a maximal curve. By Lachaud \cite[Prop. 6]{lachaud} any curve
$\fls$-covered by a $\fls$-maximal curve is also $\fls$-maximal. Then one
obtains further examples of $\fls$-maximal curves by e.g. considering
suitable quotient curves $\cH/G$, whit $G$ a subgroup of
$\aut_{\fls}(\cH)$;
see Garcia-Stichtenoth-Xing \cite{g-sti-x}, and \cite{ckt1}, \cite{ckt2}.
As a matter of fact, all the known examples of $\fls$-maximal curves arise 
in this way.
   \begin{problem}\label{problem4.1} Is any $\fls$-maximal curve
$\fls$-covered by $\cH$?
   \end{problem}
Further properties of maximal curves can be found
in \cite{fgt}, \cite{ft2}, \cite{kt1}, \cite{kt2} and the references
therein.
   \end{example} 

If $q$ is not a square, the Hasse-Weil bound was improved by Serre
\cite[Thm. 1]{serre} as follows (see also \cite[Thm. V.3.1]{sti})
   $$
|\#\cX(\fq)-(q+1)|\le \lfloor 2\sq\rfloor g\, .
   $$
   \begin{lemma}\label{lemma4.2} The following statements are equivalent:
\begin{enumerate}
   \item[\rm(1)] $\cX$ is maximal with respect to Serre's bound$;$
   \item[\rm(2)] $\alpha_i+\bar\alpha_i=-\lfloor 2\sqrt q \rfloor$ for
$i=1,\ldots g;$
   \item[\rm(3)] $h_{\cX,q}(t)=(t^2+\lfloor 2\sqrt q\rfloor t +q)^g.$
\end{enumerate}
   \end{lemma}
   \begin{proof} $\cX$ is maximal with respect to Serre's bound if and
only if $\sum_{i=1}^g(\alpha+\bar\alpha_i)=-\lfloor 2\sqrt q\rfloor g$ if
and only if $\alpha_i+\bar\alpha_i=-\lfloor 2\sqrt q\rfloor$. Now, as we
can assume $\alpha_i\bar\alpha_i=q$ by Remark \ref{rem4.1}(ii) so that 
$h_{\cX,q}(t)=\prod_{i=1}^g(t-\alpha_i)(t-\bar\alpha_i)$, the
result follows.
   \end{proof}   
   \begin{corollary}\label{cor4.2} We have $g\leq (q^2-q)/(\lfloor 2\sqrt
q\rfloor^2+\lfloor2\sqrt q\rfloor-2q)$ whenever $\cX$ is maximal with
respect to Serre's bound. 
   \end{corollary}
   \begin{proof} As in the proof of Corollary
\ref{cor4.1} we use $\cX(\fq)\subseteq \cX(\mathbf F_{q^2})$. We have
$\alpha_i+\bar\alpha_i=-\lfloor 2\sqrt q\rfloor$ and
$\alpha_i\bar\alpha_i=q$ so that $\alpha_i^2+\bar\alpha_i^2=\lfloor 2\sqrt
q\rfloor^2-2q$; hence 
   $$
\#\cX(\fq)=q+1+\lfloor 2\sqrt q\rfloor\leq\#\cX(\mathbf F_{q^2})=q^2+1-
   (\lfloor2\sqrt q\rfloor^2-2q)g\, ,
   $$
and the result follows.
   \end{proof} 
   \begin{remark}\label{rem4.2} The proofs of the following statements are
similar to the proofs of Lemmas \ref{lemma4.1} and
\ref{lemma4.2}.
   \begin{enumerate} \item[\rm(i)] A curve $\cX$ defined over $\fls$ is
{\em $\fls$-minimal}; i.e., $\#\cX(\fls)=\ell^2+1-2\ell g$ if and only if
$h_{\cX,\ell^2}(t)=(t-\ell)^{2g}$.
 
\item[\rm(ii)] A curve $\cX$ defined over $\fq$ is minimal with respect to
Serre's bound; i.e., $\#\cX(\fq)=q+1-\lfloor2\sq\rfloor g$ if and only if
$h_{\cX,q}(t)=(t^2-\lfloor 2\sq\rfloor t+q)^g$.
    \end{enumerate}
   \end{remark}

   \begin{example} (The Klein quartic)\label{ex4.1} Let $\cX$ be the plane
curve over $\F$ defined by 
   $$
X^3Y+Y^3Z+Z^3X=0\, .
   $$
It is easy to see that $\cX$ is non-singular if and only if ${\rm
char}(\F)\neq 7$; in this case $\cX$ has genus 3. This curve was
considered by many authors since the time of Klein who showed that
$\aut(\cX)$ reaches the Hurwitz bound for the number of automorphism of
curves of genus 3 whenever ${\rm char}(\F)=0$. A connection with the Fano
plane was noticed by Pellikaan \cite{pellikaan}.
   \begin{claim*} $\cX$ defined over $\mathbf F_8$ reachs the Serre's
bound$;$ i.e$,$ $\# \cX(\mathbf F_8)=1+9+\lfloor2\sqrt{8}\rfloor 3=24.$
   \end{claim*} 
To see this we first notice that $(1:0:0), (0:1:0), (0:0:1)$ are $\mathbf
F_8$-rational points (this is true for any field where $\cX$ is defined).
Now (cf. \cite[p. 10]{pellikaan}) we look for $(x:y:1)\in \cX$ such that
$x\neq 0, y\neq 0$ and such that $x^7=1$. We have
   $$
0=x^3y++y^3+x=x^3y+x^7y^3+x=x(x^2y+(x^2y)^3+1)\, ;
   $$ 
i.e., $t^3+t+1=0\ (*)$ with $t=x^2y\ (*_1)$. Conversely, it is easy to see
that equation $(*)$ is irreducible over $\mathbf F_2$ and hence its three
roots are in $\mathbf F_8$. Then once $x\in\mathbf F_8^*$ we have $y\in
\mathbf F_8^*$ by $(*_1)$. Therefore we have 21 such points $(x:y:1)$ and
the claim follows.
 
Then $h_{\cX,8}(t)=(t^2+5t+8)^3$ by Lemma \ref{lemma4.2}. 
   \begin{claim*} $h_{\cX,2}(t)=t^6+5t^3+8;$ in particular $\#\cX(\mathbf
F_2)= 3.$
   \end{claim*}
Let $h_{\cX,2}(t)=\prod_{i=1}^3(t-\beta_i)(t-\bar\beta_i)$. Then
$\beta_i^3+\bar\beta_i^3=-5$ (cf. Lemma \ref{lemma4.2}) so that
$\beta_i^3$ and $\bar\beta_i^3$ are roots of $T^2+5T+8=0$; then
$h_{\cX,2}(t)=t^6+5t^3+8$ so that $\#\cX(\mathbf F_2)=2+1-0=3$.
    
Finally, we mention that $\cX$ is $\fls$-maximal if and only if
either $\ell=p^{6v+1}$ and $p\equiv 6\pmod{7}$, or $\ell=p^{6v+3}$ and
$p\equiv 3,5,6\pmod{7}$, or $\ell=p^{6v+5}$ and $p\equiv
6\pmod{7}$; see \cite[Cor. 3.7(2)]{akt}.
   \end{example}

   \begin{remark} (Lewittes \cite[Thm. 1(b)]{lewittes})\label{rem4.3} Let
$P\in\cX(\fq)$ and $f:\cX\to\P^1(\bar\fq)$ be the
$\fq$-rational function on $\cX$ such that $\div_\infty(f)=n_1(P)P$. Then
$\cX(\fq)\subseteq f^{-1}(\P^1(\fq))=\{P_1\}\cup f^{-1}(\fq)$ and hence
   $$
\#\cX(\fq)\leq 1+qn_1(P)\, .
   $$
   \end{remark}

Now from Corollaries \ref{cor4.1} and \ref{cor4.2} we see that neither
the Hasse-Weil bound nor Serre's bound is effective to estimate
$\#\cX(\fq)$ whenever $g$ is large with respect to $q$. So in general one
studies the number
  $$
N_q(g):={\rm max}\{\#\cY(\fq): \text{$\cY$ curve of genus $g$ defined 
over $\fq$}\}\, .  
  $$ 
For instance $N_q(0)=q+1$, and Example \ref{ex4.1} shows that $N_8(3)=24$.
The study of the actual value of $N_q(g)$ was initiated by Serre
\cite{serre} who computed $N_q(1)$ and $N_q(2)$. Further properties were
proved by Serre himself \cite{serre1}, Lauter \cite{lauter}, and
Kresh-Wetherell-Zieve \cite{kwz}. Tables for $N_q(g)$ with $q$ and $g$
small can be found in van der Geer-van der Vlugt \cite{geer-vlugt}.

   \begin{definition*} A curve $\cX$ of genus $g$ and defined over $\fq$
is called {\em optimal} (with respect to $g$ and $q$) if
$\#\cX(\fq)=N_q(g)$.
   \end{definition*} 

If $q=\ell^2$ and $\cX$ is $\fls$-maximal then $\cX$ is certainly optimal.
We already noticed (Example \ref{ex4.0}) that the Hermitian curve $\cH$ is
$\fls$-maximal whose genus attains the bound in Corollary \ref{cor4.1}.
Indeed, this property characterizes Hermitian curves:
  
   \begin{theorem} {\rm (R\"uck-Stichtenoth \cite{r-sti})}\label{thm4.1} A
$\fls$-maximal curve $\cX$ has genus $\ell(\ell-1)/2$ if and only if $\cX$
is $\fls$-isomorphic to the Hermitian curve of equation (\ref{eq3.5}).
   \end{theorem}

This result follows from Theorem \ref{thm4.21}.

Next we discuss optimal curves for $\sq \not\in\N$. Besides some curves of
small genus (see above), the only known examples of optimal curves are the
Deligne-Lusztig curves $\cS$ and $\cR$ associated to the Suzuki group
$Sz(q)$, $q=2^{2s+1}$, $s\geq 1$, and to the Ree group $R(q)$,
$q=3^{2s+1}$, $s\geq 1$, respectively \cite[Sect. 11]{deligne-lusztig}. As
a matter of terminology, $\cS$ (resp. $\cR$) will be call {\em the Suzuki
curve} (resp. {\em the Ree curve}). After the work of Hansen-Stichtenoth
\cite{hansen-sti}, Hansen \cite{hansen}, Pedersen \cite{pedersen},
Hansen-Pedersen \cite{hansen-pedersen}, the curves $\cS$ and $\cR$ can be
characterized as follows.

   \begin{theorem}\label{thm4.2} The curves $\cS$ and $\cR$ 
are the unique curves (up to $\fq$-isomorphic) $\cX$ defined over $\fq$
such that the following three conditions hold:
   
   \begin{enumerate}
\item[\rm(1)] $\#\cX(\fq)=q^2+1$ (resp. $\#\cX(\fq)=q^3+1);$
\item[\rm(2)] $\cX$ has genus $q_0(q-1)$ (resp. $3q_0(q-1)(q+q_0+1)/2$)$,$ 
where $q_0:=2^s$ (resp. $3^s$)$;$
\item[\rm(3)] $\aut_{\fq}(\cX)=Sz(q)$ (resp. $\aut_{\fq}(\cX)=R(q)$)$.$
   \end{enumerate}

Moreover, the Suzuki curve $\cS$ (resp. the Ree curve $\cR$) is the 
non-singular model of 
   $$
Y^qZ^{q_0}-YZ^{q+q_0-1}=X^{q_0}(X^q-XZ^{q-1})\, ,
   $$
(resp.
   $$
\begin{cases}
   Y^qW^{q_0}-YW^{q+q_0-1}=X^{q_0}(X^q-XW^{q-1}) & \text{}\\
   Z^qW^{2q_0}-YW^{q+2q_0-1}=X^{2q_0}(x^q-XW^{q-1}))\, .
\end{cases}
   $$
   \end{theorem}
In Sect. \ref{s4.3} we prove a stronger version of this theorem for the
Suzuki curve.

   \begin{lemma}\label{lemma4.3} Let $\cX$ be a curve defined over $\fq$
such that (1) and (2) in Theorem \ref{thm4.2} hold. Then $\cX$ is optimal;
moreover:
  \begin{enumerate}
\item[\rm(1)] If $q=2^{2s+1},$ $h_{\cX,q}(t)=(t^2+2q_0t+q)^{q_0(q-1)};$

\item[\rm(2)] If $q=3^{2s+1},$ $h_{\cX,q}(t)=(t^2+3q_0t+q)^{q_0(q^2-1)}
(t^2+q)^{q_0(q-1)(q+3q_0+1)/2}.$
  \end{enumerate}
    \end{lemma}
    \begin{proof} It is easy to see that Serre's bound is not effective to
bound $\#\cX(\fq)$; in this case one uses the so-called ``explicit
formula" (\ref{eq4.2}) of Weil \cite{serre}: (following Stichtenoth
\cite[p. 183]{sti}) Let
$h_{\cX,q}(t)=\prod_{i=1}^g(t-\alpha_i)(t-\bar\alpha_i)$,
$\alpha_i=\sqrt{q}e^{\sqrt{-1}\theta_i}$, and write
    $$
q^{-i/2}\#\cX(\fqi)=q^{i/2}+q^{-i/2}-q^{-i/2} 
\sum_{j=1}^{g}(\alpha_j^i+\bar\alpha_j^i)\, ;
    $$
this equation can we rewritten as
    $$
\#\cX(\fq)c_iq^{-i/2}=c_iq^{i/2}+c_iq^{-i/2}+c_iq^{-i/2}
\sum_{j=1}^{g}(\alpha_j^i+\bar\alpha_j^i)-(\#\cX(\fqi)-
\#\cX(\fq)c_iq^{-i/2}\, ,
    $$ 
where $c_i\in \R$. Now suppose that $c_1,\ldots,c_m$ are given real
numbers. Then from the above equation we obtain:
   \begin{equation}\label{eq4.2}
   \begin{split}
\#\cX(\fq)\lambda_m(q^{-1/2}) &=\lambda_m(q^{1/2})+\lambda_m(q^{-1/2})+g-
\sum_{j=1}^gf_m(q^{-1/2}\alpha_j)-\\
                              &{} \sum_{i=1}^m(\#\cX(\fqi)-
\#\cX(\fq))c_iq^{-i/2}\, ,
   \end{split}
   \end{equation}
where $\lambda_m(t):=\sum_{i=1}^mc_it^i$ and
$f_m(t):=1+\lambda_m(t)+\lambda_m(t^{-1})$. Note that $f_m(t)\in\R$
whenever $t\in\C$ and $|t|=1$.

Case $q=2^{2s+1}$ and $g=q_0(q-1)$. Here we choose $m=2$, $c_1=\sqrt 2/2$,
$c_2=1/4$. Then
$f_2(e^{\sqrt{-1}\theta})=1+\sqrt{2}cos\theta+cos(2\theta)/2=
(cos\theta+\sqrt{2}/2)^2\geq 0$. Then from (\ref{eq4.2}) we have
   $$
\#\cX(\fq)\lambda_2(q^{-1/2})\leq
\lambda_2(q^{1/2})+\lambda_2(q^{-1/2})+g\, ,
   $$ 
so that $\#\cX(\fq)\leq q^2+1$, and hence $\cX$ is optimal. Moreover, as
$\#\cX(\fq)=q^2+1$ we must have $f_2(q^{-1/2}\alpha_j)=0$ by (\ref{eq4.2})
so that $cos\theta_j=-\sqrt 2/2$. Then $\alpha_j+\bar\alpha_j=-2q_0$ and
the result on $h_{\cX,q}(t)$ follows.

Case $q=3^{2s+1}$ and $g=3q_0(q-1)(q+q_0+1)/2$. Here we use $m=4$,
$c_1=\sqrt{3}/2$, $c_2=7/12$, $c_3=\sqrt{3}/6$, $c_4=1/12$. Then
$f_4(e^{\sqrt{-1}\theta})=1+\sqrt 3cos\theta+7cos(2\theta)/6+\sqrt
3cos(3\theta)/3+cos(4\theta)/6=(1+\sqrt{3}cos\theta+cos2\theta)^2/3\geq
0$. Then from (\ref{eq4.2})
   $$
\#\cX(\fq)\lambda_4(q^{-1/2})\le
\lambda_4(q^{1/2})+\lambda_4(q^{-1/2})+g\, ,
   $$
so that $\cX(\fq)\le q^3+1$. Moreover,
$1+\sqrt{3}cos\theta_j+cos2\theta_j=0$ whenever $\cX(\fq)=q^3+1$. Hence
$cos\theta_j=0$ or $cos\theta_j=-\sqrt{3}/2$ so that
   $$
h_{\cX,t}(t)=(t^2+3q_0t+q)^A(t^2+q)^{g-A}\, ,
   $$
where $A$ is the number of $j$'s such that $cos\theta_j=-\sqrt{3}/2$. To
compute $A$ we use the facts that $a_1=\#\cX(\fq)-(q+1)=q^3-q$ and
$a_{2g-1}=q^{g-1}a_1$. We have $a_{2g-1}=h_{\cX,q}'(0)=3q_0q^{g-1}A$ and
hence that $A=q_0(q^2-1)$.
   \end{proof} 
\subsection{A $\fq$-divisor from the Zeta Function}\label{s4.1}

Assume now that $\cX(\fq)\neq\emptyset$, and fix a $\fq$-rational point
$P_0\in\cX$. Let $f=f^{P_0}:P\to [P-P_0]$ be the canonical map from $\cX$
to its Jacobian over $\fq$, $\cJ\cong \{D\in\Div(\cX):
\deg(D)=0\}/\{\div(x):x\in\bar\fq(\cX)^*\}$. Let $\fro'$ be the Frobenius
morphism on $\cJ$ induced by $\fro$.

We recall some facts concerning the characteristic polynomial of $\fro'$
which in fact turns out to be the polynomial $h(t)=h_{\cX,q}(t)$ which was
defined in Remark \ref{rem4.1}; see e.g. \cite[p. 205]{mumford}, or
\cite[proof of Thm. 19.1]{milne}.

For a prime $\ell$ different from ${\rm char}(\fq)$, let $\cJ_{\ell^i}$
denote the kernel of the isogeny $\cJ\to \cJ$, $P\mapsto \ell^iP$. Then
one defines the {\em Tate modulo} associated to $\cJ$ as the inverse limit
of the groups $\cJ_{\ell^i}$, $i\ge 1$, with respect to the maps
$\cJ_{\ell^{i+1}}\to \cJ_{\ell^i}$, $P\mapsto \ell P$. We have that $\#
\cJ_{\ell^i}=(\ell^i)^{2g}$ \cite[p. 62]{mumford} so that $\cJ_{\ell^i}$
is
a finite abelian group such that for all $j$, $1\leq j\leq i$ it contains
exactly $(\ell^j)^{2g}$ elements of order $\ell^j$. Therefore 
   $$
\cJ_{\ell^i}\cong (\Z/\ell^i\Z)^{2g}\qquad\text{and hence}\qquad
T_\ell(\cJ)\cong \Z_\ell^{2g}\, ,
   $$
where $\Z_\ell$ denotes the $\ell$-adic integers. Thus $T_\ell(\cJ)$ is a
free $\Z_\ell$-module of rank $2g$. Now clearly
$\fro'(\cJ_{\ell^i})\subseteq \cJ_{\ell^i}$ and hence $\fro'$ gives rise
to a $\Z_\ell$-linear map $T_\ell(\fro')$ on $T_\ell(\cJ)$. Let $\pi$ be
the
characteristic polynomial of $T_\ell(\fro')$. A priory we have that $\pi$
is a
polynomial of degree $2g$ with coefficients in $\Z_\ell$. As a matter of
fact, $\pi\in\Z[t]$ \cite[proof of Ch. IV, Thm. 4]{mumford}, and $\pi=h$
as we mentioned before. In
particular, the minimal polynomial $m$ of $T_\ell(\fro')$ has integral
coefficients. We claim that
   \begin{equation}\label{eq4.11}
m(\fro')=0\qquad \text{on $\cJ$}\, .
   \end{equation}
To see this, notice that any endomorphism $\alpha\in {\rm 
End}(\cJ):\cJ\mapsto\cJ$ acts on
$T_\ell(\cJ)$ giving rise to a $\Z_\ell$-linear map $T_\ell(\alpha)$. This
action is injective because ${\rm End}(\cJ)$ is torsion free and because
of
\cite[Ch. IV, Thm. 3]{mumford}. Now, as $m(\fro')\in {\rm End}(\cJ)$, we
have
   $$
0=m(T_\ell(\fro'))=T_\ell(m(\fro'))
   $$
and (\ref{eq4.11}) follows. Moreover, it is known that $\Q\otimes {\rm
End}(\cJ)$ is a finite dimensional semisimple algebra over $\Q$ whose
center is $\Q[\fro']$ \cite[Ch. IV, Cor. 3]{mumford}, \cite[Thm.
2(a)]{tate}. In particular, $\Q[\fro']$ is semisimple and it is not
difficult to see that $T_\ell(\fro')$ is semisimple; cf. \cite[p.
251]{mumford}. This means that 
   \begin{equation*}
m(t)=\prod_{i=1}^T h_i(t)\, ,
   \end{equation*}
where $h_1(t), \ldots, h_T(t)$ are the irreducibles $\Z$-factors of
$h(t)$. Let $U$ be the degree of $m(t)$ and let $b_1,\ldots,b_U\in \Z$
be the coefficients of $m(t)-t^U$; i.e,
   $$
m(t)=t^U+\sum_{i=1}^{U}b_it^{U-i}\, .
   $$
Thus $(\fro')^U+\sum_{i=1}^{U}b_i(\fro')^{U-i}=0$ by (\ref{eq4.11}).
Now we evaluate the left hand side of this equality at $f(P)=[P-P_0]$, and 
by using the fact that $\fro'\circ f=f\circ\fro$ we find that
   $$
f(\fro^U(P))+\sum_{i=1}^{U}a_if(\fro^{U-i}(P))=0\, ,\qquad P\in 
\cX\, ;
   $$
   \begin{equation}\label{eq4.12}
\text{i.e.,}\qquad\fro^U(P)+\sum_{i=1}^Ub_i\fro^{U-i}(P)\sim
(1+\sum_{i=1}^Ub_i)P_0=m(1)P_0\, .
   \end{equation}
This equivalence is the motivation to define on $\cX$ the linear series
   \begin{equation}\label{eq4.121}
\cD_\cX:= ||m(1)|P_0|\, ,
   \end{equation}
which is clearly independent of $P_0$ being $\fq$-rational.

   \begin{problem}\label{problem4.11} For a curve $\cX$ over $\fq$, how is
the interplay among its $\fq$-rational points, its Weierstrass points, 
its $\cD_\cX$-Weierstrass points, and the support of the $\fq$-Frobenius
divisor of $\cD_\cX$.
   \end{problem}

Next we discuss some properties of $\cD_\cX$.
    \begin{lemma}\label{lemma4.11} 
\begin{enumerate} 
   \item[\rm(1)] If $P, Q\in \cX(\fq)$, then $m(1)P\sim m(1)Q$; in
particular, $|m(1)|$ is a Weierstrass non-gap at each $P\in \cX(\fq)$.
   
   \item[\rm(2)] If $\#\cX(\fq)\ge 2g+3$, then there exists $P_1\in
\cX(\fq)$ such that $|m(1)|-1$ and $|m(1)|$ are Weierstrass non-gaps at
$P_1$. 
    \end{enumerate}  
    \end{lemma}
    \begin{proof} (1) It follows immediately from (\ref{eq4.12}). 

(2) (Following Stichtenoth-Xing \cite[Prop. 1]{sti-x}) Let $Q\in
\cX(\fq)\setminus\{P_0\}$. From (1),
there exists a 
morphism 
$x:\cX \to \P^1(\bar\fq)$ with $\div(x)=|m(1)|P_0-|m(1)|Q$. Let $n$ be the
number of $\fq$-rational points of $\cX$ which are unramified for $x$. Let 
$x^s:\cX\to \P^1(\bar\fq)$ be the separable part of $x$. We have that
$\div(x^s)=|m(1)|'P_0-|m(1)|'Q$ (here $|m(1)|'$ is the separable degree of
$x$) and from the Riemman-Hurwitz applied to $x^s$ we find that
  $$
2g-2\ge |m(1)|'(-2)+2(|m(1)|'-1)+(\#\cX(\fq)-n-2)\, ,
  $$
so that $n\ge \#\cX(\fq)-2g-2$. Thus $n\geq 1$ by hypothesis, and hence
there exists $\alpha \in \fq$, $P_1\in
\cX(\fq)\setminus \{P_0, Q\}$ such that $\div(x-\alpha)=P_1+D-mQ$ with
$P_1, Q \not\in \supp(D)$. Let $y\in
\bar\fq(\cX)$ be such that $\div(y)=|m(1)|Q-|m(1)|P_1$ (cf. (1)).
Then $\div(y(x-\alpha))=D-(|m(1)|-1)P_1$ and (2) follows.
    \end{proof}
   \begin{corollary}\label{cor4.11} 
   \begin{enumerate}
\item[\rm(1)] $\cD_\cX$ is base-point-free$;$
\item[\rm(2)] If $\#\cX(\fq)\ge 2g+3$, then $\cD_\cX$ is simple$.$ 
   \end{enumerate}
   \end{corollary}
   \begin{proof} (1) follows by Lemma \ref{lemma4.11} and Example
\ref{ex1.53}

(2) Let $P_1$ be as in Lemma \ref{lemma4.11}(2), $\phi$ a morphism
associated to $\cD_\cX$, $f_1,f_2\in
\bar\fq(\cX)$ such that $\div_\infty(f_1)=(|m(1)|-1)P_1$ and
$\div_\infty(f_2)=|m(1)|P_1$. Then $[\bar\fq(\cX):\bar\fq(f_i)]$,
$i=1,2$, divides $[\bar\fq(\cX):\bar\fq(\phi(\cX))]$ and the result
follows.
   \end{proof}

Now we study $(\cD_\cX,P)$-orders. We let
$\epsilon_0=0<\epsilon_1=1<\ldots<\epsilon_N$ (resp.  
$\nu_0=0<\ldots<\nu_{N-1}$) denote the $\cD_\cX$-orders (resp. the
$\fq$-Frobenius orders) of $\cD_\cX$, where $N:=\dim(\cD_\cX)$.  Notice
that $n_N(P)=|m(1)|$ for any $P\in\cX(\fq)$ by Lemma \ref{lemma4.11}(1).
From Example \ref{ex1.53} we obtain:

   \begin{lemma}\label{lemma4.12} For $P\in\cX(\fq)$, the   
$(\cD_\cX,P)$-orders are
   $$
j_{N-i}(P)=n_N(P)-n_i(P)\, ,\qquad i=0,1,\ldots, N\, .
   $$
   \end{lemma}
This result (for $i=1$) and Remark \ref{rem4.3} yield the following.
    \begin{corollary}\label{cor4.12} Let $P\in\cX(\fq).$ 
If $\# \cX(\fq)>q(|m(1)|-b_U)+1$, then $j_{N-1}(P)<b_U.$
    \end{corollary}

   \begin{lemma}\label{lemma4.13} Suppose 
   \begin{equation}\label{eq4.13}
b_i\geq 0\, ,\quad i=1,\ldots, U\, ,
   \end{equation}
and let $P\in\cX$ such that $\fro^i(P)\neq P$ for $i=1,\ldots,U.$ Then:
   \begin{enumerate}
\item[\rm(1)] The numbers $1, b_1,\ldots,b_U$ are $(\cD_\cX,P)$-orders$;$

\item[\rm(2)] If in addition 
   \begin{equation}\label{eq4.14}
b_1 \geq b_0:=1\quad\text{and}\quad 
b_{i+1}\geq b_i,\ \text{for}\ i=1,\ldots, U-1\, ,
   \end{equation}
then $b_U$ (resp. $b_U-1$) is a Weierstrass non-gap at $P$ whenever
$\fro^{U+1}(P)\neq P$ (resp. $\fro^{U+1}(P)=P$)$.$
    \end{enumerate}
    \end{lemma}
    \begin{proof} (1) Fix $j\in\{0,1,\ldots,U\}$, and let
$Q\in\cX$ such that $\fro^{U-j}(Q)=P\ (*)$. From (\ref{eq4.12}) we have
   $$
\sum_{i\in\{0,1,\ldots,U\}\setminus\{j\}}b_i\fro^{U-i}(Q)+b_jP\sim
m(1)P_0\, .
   $$
We claim that $\fro^{U-i}(Q)\neq P$; otherwise from $(*)$ we would have
$\fro^{i-j}(P)=P$, a contradiction. This shows (1).

(2) Applying ${\fro}_{*}$ to (\ref{eq4.12}) we have
    $$
\fro^U(P)+\sum_{i=1}^Ub_i\fro^{U-i}(P)\sim m(1)P_0\sim  
\fro^{U+1}(P)+\sum_{i=1}^Ub_i\fro^{U-i+1}(P)\, ,
    $$
so that
    $$
b_UP\sim \fro^{U+1}(P)+\sum_{i=1}^U(b_i-b_{i-1})\fro^{U-i+1}(P)\, ,
    $$
and (2) follows.
    \end{proof}
   \begin{remark}\label{rem4.11} (i) Minimal curves as well as minimal
curves with respect to Serre's bound (Remark \ref{rem4.2}) do not satisfy
(\ref{eq4.13}). However we can still use (\ref{eq4.12}) to infer that
$\sq$ is a non-gap at infinitely many points of the curve provided that
the curve is minimal. Indeed, (\ref{eq4.13}) reads $\fro(P)-\sq P\sim
(1-\sq)P_0$ so that $\sq P\sim (\sq-1)P_0+\fro(P)$. In particular, if
$g\geq \sq$, a $\fq$-minimal curve is non-classical.

(ii) The Klein curve (Example \ref{ex4.1}) defined over $\mathbf F_2$ 
satisfies (\ref{eq4.13}) but not (\ref{eq4.14}).

(iii) Other examples as in (i) and (ii) can be found in Carbonne-Henocq
\cite{carbonne-henocq}.
   \end{remark}

   \begin{corollary}\label{cor4.13} Assume (\ref{eq4.13}).
   \begin{enumerate}
\item[\rm(1)] If $P\not\in\cX(\fq)$ and $\cX(\fq)=\ldots=\cX(\mathbf
F_{q^U}),$ then $1,b_1,\ldots,b_U$ are $(\cD_\cX,P)$-orders$.$

\item[\rm(2)] The numbers $1,b_1,\ldots,b_U$ are $\cD_\cX$-orders$.$ In
particular$,$ $\dim(\cD_\cX)\ge U+1$ provided that $b_i\neq b_j$ for
$i\neq j;$

\item[\rm(3)] If in addition (\ref{eq4.14}) holds and $g\ge b_U$, then
$\cX$ is non-classical$.$
    \end{enumerate}
    \end{corollary}
    \begin{proof} Lemma \ref{lemma4.13}(1) implies (1) and (2) since there
are infinitely many points $P$ such that $\fro^i(P)\neq P$ for
$i=1,\ldots,U$.
To see (3) we take $P\in\cX$ such that $\fro^{U+1}(P)\neq P$. Then $b_U\in
H(P)$ by Lemma \ref{lemma4.13}(2). If $\cX$ were classical then
$n_1(P)=g+1$ so that $g<b_U$, a contradiction.
    \end{proof}

    \begin{corollary}\label{cor4.14} Assume (\ref{eq4.13}).
    \begin{enumerate}
\item[\rm(1)] $\epsilon_N=\nu_{N-1}=b_U;$
\item[\rm(2)] $\cX(\fq)\subseteq \supp(R^{\cD}).$
   \end{enumerate}
    \end{corollary}
    \begin{proof} (1) We have $\epsilon_{N-1}\leq j_{N-1}(P)$ for
any $P$ by Corollary \ref{cor2.21}(1); thus $\epsilon_{N-1}<b_U$ by
Corollary \ref{cor4.12}. Therefore $\epsilon_N=b_U$ by Corollary
\ref{cor4.13}(2), and so
  $$
\phi^*(L_{N-1}(P))=\fro^U(P)+\sum_{I=1}^Ub_i\fro^{U-i}(P)
  $$
by (\ref{eq4.12}), where $\phi$ is a morphism associated to $\cD_\cX$. It
follows that $\phi(\fro(P))\in L_{N-1}(P)$ so that $\nu_{N-1}=\epsilon_N$.

(2) By Lemma \ref{lemma4.12} $j_N(P)=n_N(P)=m(1)$ for each $P\in\cX(\fq)$.
Since $m(1)=1+\sum_{i=1}^Ub_i>b_U=\epsilon_N$ (cf. (1)), the result
follows.
    \end{proof}

    \begin{corollary}\label{cor4.15} Assume (\ref{eq4.14}). Then
$n_1(P)\leq b_U$ for each $P\in\cX(\fq),$ and equality holds provided
that $\#\cX(\fq)\geq qb_U+1.$
    \end{corollary} 
    \begin{proof} Let $P\in \cX(\fq)$. By Lemma \ref{lemma2.511}
$n_1(P)\leq n_1(Q)$ where $Q\not\in\cW$. Therefore $n_1(P)\leq b_U$ by
Lemma \ref{lemma4.13}(2). Now if $\#\cX(\fq)\geq qb_U+1$, then
$1+qn_1(P)\geq qb_U+1$ by Remark \ref{rem4.3} and the result follows.
    \end{proof}

\subsection{The Hermitian curve}\label{s4.2} Let $\cX$ be a $\fls$-maximal
curve of genus $g$. Recall that $g\leq \l(\l+1)/2$ by Corollary
\ref{cor4.1} and that the Hermitian curve is  
$\fls$-maximal of genus $\l(\l-1)/2$ (cf. Example \ref{ex3.2}). From
Lemma \ref{lemma4.1} and (\ref{eq4.121}), $\cX$ is equipped with the
linear series $\cD_\cX:=|(\l+1)P_0|$. By Corollary \ref{cor4.11},
$\cD_\cX$ is simple and base-point-free. We see that $\cX$ satisfies
(\ref{eq4.14}) (and hence (\ref{eq4.13})); in particular $1,\l$ are
$\cD_\cX$ orders so that $N:=\dim(\cD_\cX)\geq 2$.

   \begin{theorem} {\rm (\cite[Thm. 2.4]{ft2})}\label{thm4.21} Let $\cX$
be a $\fls$-maximal curve of genus $g$. The following statements are
equivalent:
   \begin{enumerate}
\item[\rm(1)] $\cX$ is $\fls$-isomorphic to the Hermitian curve $\cH$ of
equation (\ref{eq3.5})$;$
\item[\rm(2)] $g>(\l-1)^2/4;$
\item[\rm(3)] $N=2.$
   \end{enumerate}
   \end{theorem} 
   \begin{proof} (1) implies (2) because the genus of $\cH$ is
$\l(\l-1)/2$. Assume (2) and suppose that $N\geq 3$. Then Castelnuovo's
genus bound (Remark \ref{rem1.21}) applied to $\cD_\cX$ would yield
$g\leq
(\l-1)^2/4$, a contradiction. Finally let $N=2$. By (\ref{eq4.12})
$(\l+1)P\sim (\l+1)P_0$ for any $P\in\cX(\fls)$ and hence we can assume
that $\l,\l+1\in H(P_0)$ by Lemma \ref{lemma4.11}(2); in this case, as
$N=2$, $n_1(P_0)=\l$ and $n_2(P_0)=\l+1$. Let
$\epsilon_0=0<\epsilon_1=1<\epsilon_2$ (resp. $\nu_0=0<\nu_1$) denote the
$\cD_\cX$-orders (resp. $\fls$-orders) of $\cX$. Then
$\epsilon_2=\nu_1=\l$ by Corollary \ref{cor4.14}. Let $x,y\in \fls(\cX)$
such that $\div_\infty(x)=\l P_0$ and $\div_\infty(y)=(\l+1)P_0$. We have
that$x$ is a separating variable (Lemma \ref{lemma1.520}) and therefore
   \begin{equation*}
V^{0,1}_{1,x,y,;x}={\rm det}\begin{pmatrix} 1& x^{\l^2} & y^{\l^2} \\
                                             1& x        & y \\ 
                                             0& 1        & D^1_xy
                              \end{pmatrix}=
(x-x^{\l^2})D^1_xy-(y-y^{\l^2})=0\, .\tag{$*$}
   \end{equation*}
   \begin{claim*} There exists $f\in \bar\fls(\cX)$ such that
$D^1_xy=f^\l.$
   \end{claim*} 
To proof this we have to show that $D^i_x(D^1_xy)=0\, (*_1)$ for $1\leq
i<\l$ by Remark \ref{rem2.15}(ii). We apply $D^1_x$ to $(*)$: 
$(x-x^{\l^2})D^1_x(D^1_xy)=0$ and so $(*_1)$ holds for $i=1$. Suppose that
$(*_1)$ is true for $i=1,\ldots,j$, $1\leq j\le \l-2$. We apply
$D^{j+1}_x$ to $(*)$ and using the inductive hypothesis and Remark
\ref{rem2.15}(i) we find that $(x-x^{\l^2})D^{j+1}_x(D^1_xy)=D^{j+1}_xy$.
It
turns out that
   $$
W^{0,1,j+1}_{1,x,y;x}=\begin{pmatrix} 1& x & y\\
                                      0& 1 & D^1_xy\\
                                      0& 0 & D^{j+1}_xy
                                    \end{pmatrix}=D^{j+1}_xy=0\, ,
   $$
since $\epsilon_2=\l$, and the claim follows.
    \begin{claim*} $\#x^{-1}(x(P))=\l$ for $P\neq P_0.$
    \end{claim*}
From $(*)$ $v_{P_0}(D^1_xy)=-\l^2$. Let $t$ be a local parameter at $P_0$.
Then $v_{P_0}(D^1_tx)=\l^2-l-2$ since $D^1_ty=D^1_txD^1_xy$ by the chain
rule (\ref{eq2.13}). We have
that $\deg(dx)=2g-2$ (see Example \ref{ex1.11}) and that $v_P(x)\ge 0$ for
$P\neq P_0$. Therefore $2g-2\geq \l^2-l-2$; i.e., $g\geq l(l-2)/2$; i.e.
$g=\l(\l-1)/2$ by Corollary \ref{cor4.1}. It follows that $v_P(dx)=0$ for
$P\neq P_0$ and so the claim.

We conclude that $D^1_xy=f^{\l}$ with $\div_infty f=\l P_0$; moreover
$f\in\fq(\cX)$ since $D^1_xy\in\fq(\cX)$. Then $f=a+bx$ with $a,b\in\fls$
and $(*)$ gives a relation of type 
  $$
(y_1^\l+y_1-x_1^{\l+1})^\l=y_1^\l+y_1^\l-x_1^{\l+1}\, .
  $$
Finally we have that $y_1^\l+y_1-x_1^{\l+1}=c\in\mathbf F_\l$ and with
$y_2:=y_1+\lambda$, $\lambda^\l+\lambda=a$, we have that (\ref{eq3.5})
holds; i.e., $\cX$ is $\fls$-isomorphic to $\cH$.
   \end{proof}
   \begin{corollary} {\rm (\cite{ft1})}\label{cor4.21} The genus $g$ of a
$\fls$-maximal curve satisfies
   $$
\text{either}\qquad g\leq(\l-1)^2/4\qquad\text{or}\qquad g=\l(\l-1)/2\, .
   $$
   \end{corollary}
   \begin{remark}\label{rem4.21} This result was improved in \cite{kt2}
where it is shown that $g\le (\l^2-\l+1)/6$ whenever $g<(\l-1)^2/4$.
   \end{remark}

\subsection{The Suzuki curve}\label{s4.3} Set $q_0:=2^s$, $s\in\N$,
$q:=2q_0^2$. Let $\cX$ be a curve defined over $\fq$ of genus $g$ such
that
   \begin{equation}\label{eq4.31}
g=q_0(q-1)\qquad\text{and}\qquad\#\cX(\fq)=q^2+1\, .
   \end{equation}
The main result of this sub-section is the following theorem 
which improves Theorem \ref{thm4.2} for the Suzuki curve $\cS$.

   \begin{theorem}\label{thm4.31} A curve $\cX$ defined over $\fq$ 
is $\fq$-isomorphic to the Suzuki curve
$\cS$ if and only if (\ref{eq4.31}) hold true.
   \end{theorem}

   \begin{problem}\label{problem4.31} Can we expect a similar result for
the Ree curve?
   \end{problem}

If (\ref{eq4.31}) hold, then $h_{\cX,q}(t)=(t^2+2q_0t+q)^g$ by Lemma
\ref{lemma4.3}(1), and from (\ref{eq4.121}) we see that $\cX$
is equipped with the linear series
   $$
\cD_\cX=|(q+2q_0+1)P_0|\, ,\qquad P_0\in\cX(\fq)\, .
   $$
The results of Sect. \ref{s4.1} applied to this
case are summarized in the following proposition. Let $N:=\dim(\cD_\cX)$,
$\epsilon_0=0<\epsilon_1=1<\ldots<\epsilon_N$ (resp.
$\nu_0=0<\ldots<\nu_{N-1}$) be the $\cD_\cX$-orders (resp. $\fq$-Frobenius
orders) of $\cX$.
   \begin{proposition}\label{prop4.31}
\begin{enumerate}
   \item[\rm(1)] $j_N(P)=n_N(P)=q+2q_0+1$ for any $P\in\cX(\fq);$ in
addition, there exists $P_1\in\cX(\fq)$ such that $n_{N-1}(P_1)=q+2q_0;$
   \item[\rm(2)] $\cD_\cX$ is simple and base-point-free$;$
   \item[\rm(3)] $2q_0$ and $q$ are $\cD_\cX$-orders so that $N\geq 3;$
   \item[\rm(4)] $\epsilon_N=\nu_{N-1}=q;$
   \item[\rm(5)] $n_1(P)=q$ for any $P\in\cX(\fq).$
\end{enumerate}
   \end{proposition}
From (5) and (1) above and Lemma \ref{lemma4.12},  
$j_{N-1}(P)=j_N(P)-n_1(P)=2q_0+1$ for any $P\in\cX(\fq)$
so that
  $$
2q_0\le \epsilon_{N-1}\leq 2q_0+1\, .
  $$
   \begin{lemma}\label{lemma4.31} $\epsilon_{N-1}=2q_0.$
   \end{lemma}
   \begin{proof} Suppose that $\epsilon_{N-1}> 
2q_0$. Then $\epsilon_{N-2}=2q_0$ and 
$\epsilon_{N-1}=2q_0+1$. By Corollary \ref{cor3.2}(1) $\nu_{N-2}\leq
j_{N-1}(P)-j_1(P)\leq 2q_0=\epsilon_{N-2}$, and thus the $\fq$-Frobenius
orders of $\cD_\cX$ would be
$\epsilon_0,\epsilon_1,\ldots,\epsilon_{N-2}$,  
and $\epsilon_N$. Now from Proposition \ref{prop3.2}(1) 
   \begin{equation}\label{eq4.311}
v_P(S)\geq \sum_{i=1}^{N}(j_i(P)-\nu_{i-1})\geq
(N-1)j_1(P)+1+2q_0\ge N+2q_0\ ,
   \end{equation}
for $P\in\cX(\fq)$ so that
$\deg(S)=(\sum_i\nu_i)(2g-2)+(q+N)(q+2q_0+1)\geq (N+2q_0)\#\cX(\fq)$. From 
the identities $2g-2=(2q_0-2)(q+2q_0+1)$ and 
$\#\cX(\fq)=(q-2q_0+1)(q+2q_0+1)$ we would have
   $$
\sum_{i=1}^{N-2}\nu_i=\sum_{i=1}^{N-2}\epsilon_i\geq (N-1)q_0\, .  
   $$ 
Now, as $\epsilon_i+\epsilon_j\leq \epsilon_{i+j}$ for $i+j\leq N$ by
Corollary \ref{cor2.23},  
   $$ 
(N-1)2q_0=(N-1)\epsilon_{N-2}\geq 2\sum_{i=0}^{N-2}\epsilon_i\geq
2(N-1)q_0\, ,  
   $$ 
and hence $\epsilon_i+\epsilon_{N-2-i}=\epsilon_{N-2}$ for
$i=0,\ldots,N-2$. In particular, 
$\epsilon_{N-3}=2q_0-1$ and by the $p$-adic criterion (Lemma
\ref{lemma2.31}) we would have $\epsilon_i=i$ for $i=0,1,\ldots,N-3$. 
Then $N=2q_0+2$. Now from Castelnuovo's genus bound (Remark \ref{rem1.21}) 
   $$
2g=2q_0(q-1)\leq (q+2q_0-(N-1)/2)^2)/(N-1)\, ;
   $$
i.e., $2q_0(q-1)< (q+q_0)^2/2q_0=q_0q+q/2+q_0/2$, a contradiction.
   \end{proof}
   \begin{corollary}\label{cor4.31} There exists $P_1\in \cX(\fq)$ such
that  
   $$
\left\{ \begin{array}{ll}
j_1(P_1)=1         & {} \\
j_i(P_1)=\nu_{i-1}+1 & \mbox{if}\ i=2,\ldots, N-1.
\end{array}\right.
   $$
   \end{corollary}
   \begin{proof} Since we already observed that $v_P(S)\geq
(N-1)j_1(P)+2q_0+1\geq N+2q_0$ for
$P\in\cX(\fq)$, it is enough to show that there exists
$P_1\in \cX(\fq)$ such that $v_{P_1}(S)=N+2q_0$. Suppose that $v_P(S)\geq
N+2q_0+1$ for any $P\in \cX(\fq)$. Then by Theorem \ref{thm3.1} 
   $$
\sum_{i=0}^{N-1}\nu_i \geq q+Nq_0+1\ ,
   $$
so that
   $$
\sum_{i=0}^{N-1}\epsilon_i \geq Nq_0+2\, ,
   $$
because $\epsilon_1=1$, $\nu_{N-1}=q$ and $\nu_i\leq \epsilon_{i+1}$. Then
from Corollary \ref{cor2.23} we would have 
$N\epsilon_{N-1}\geq 2Nq_0+4$; i.e., $\epsilon_{N-1}>2Nq_0$, a
contradiction by Lemma \ref{lemma4.31}.
   \end{proof}
   \begin{lemma}\label{lemma4.32} 
   \begin{enumerate} 
\item[\rm(1)] $\nu_1>\epsilon_1=1;$
\item[\rm(2)] $\epsilon_2$ is a power of two$.$
   \end{enumerate} 
   \end{lemma}
   \begin{proof} If $\nu_1>\epsilon_1=1$, then $\nu_1=\epsilon_2$ and
it must be a power of two by the $p$-adic criterion (Lemma
\ref{lemma2.31}): i.e., (1) implies (2). 
Suppose now that $\nu_1=1$. Then from Corollary \ref{cor4.31} 
there exists a point $P_1\in \cX(\fq)$ such 
that $j_1(P_1)=1, j_2(P_1)=2$; thus
   $$
H(P_1)\subseteq H:=\langle q, q+2q_0-1, q+2q_0, q+2q_0+1\rangle\, ,
   $$
by Proposition \ref{prop4.31}(1)(5) and Lemma \ref{lemma4.12}. In
particular $g=q_0(q-1)\leq \tilde g:=\#(\N_0\setminus H)$. This is a
contradiction as follows immediately from the claim below.
    \begin{claim*} $\tilde g=g-q_0^2/4.$
    \end{claim*}
In fact, $L:=\cup_{i=1}^{2q_0-1}L_i$ is a complete system of residues
module $q$, where
   $$
\begin{array}{lll}
L_i & = &
\{iq+i(2q_0-1)+j: j=0,\ldots,2i\}\quad  \mbox{if}\ \ 1\le i\le q_0-1,\\
L_{q_0} & = & \{q_0q+q-q_0+j:j=0,\ldots,q_0-1\},\\
L_{q_0+1} & = & \{(q_0+1)q+1+j:j=0,\ldots,q_0-1\},\\ 
L_{q_0+i} & = & 
\{(q_0+i)q+(2i-3)q_0+i-1+j: 
j=0,\ldots,q_0-2i+1\}\cup\\
          &   & \{(q_0+i)q+(2i-2)q_0+i+j: j=0,\ldots q_0-1\}\quad  
\mbox {if}\ \ 2\le i\le q_0/2,\\
L_{3q_0/2+i} & = & 
\{(3q_0/2+i)q+(q_0/2+i-1)(2q_0-1)+q_0+2i-1+j:\\
             &  &  
j=0,\ldots,q_0-2i-1\}\quad \mbox {if}\ \ 1\le i\le q_0/2-1.
\end{array}
    $$
Moreover, for each $\l \in L$, $\l \in H$ and $\l-q\not\in H$. Hence
$\tilde g$ can be computed by summing up the coefficients of $q$ from the 
above list (see e.g. \cite[Thm. p.3]{selmer}); i.e.,
   \begin{equation*}
\begin{array}{lll}
\tilde g & = & \sum_{i=1}^{q_0-1}i(2i+1)+q_0^2+(q_0+1)q_0+
\sum_{i=2}^{q_0/2}(q_0+i)(2q_0-2i+2)+\\
         &   & 
\sum_{i=1}^{q_0/2-1}(3q_0/2+i)(q_0-2i)=q_0(q-1)-q_0^2/4\, .
\end{array}
   \end{equation*}
   \end{proof}
In the remaining part of this sub-section we let $P_0=P_1$ be a
$\fq$-rational point satisfying Corollary \ref{cor4.31}; we set $n_i:=
n_i(P_1)$ and $v:=v_{P_1}$.

Lemma \ref{lemma4.32}(1) implies $\nu_i=\epsilon_{i+1}$ for 
$i=1,\ldots, N-1$. Therefore from Corollary \ref{cor4.31} and Lemma
\ref{lemma4.12} we have 
    \begin{equation}\label{eq4.32}
\left\{\begin{array}{ll}
n_i=2q_0+q-\epsilon_{N-i} & \mbox{if}\ i=1,\ldots N-2\\
n_{N-1}=2q_0+q,\ \ n_N=1+2q_0+q. & {}
\end{array}\right.
     \end{equation}
Let $x, y_2,\ldots, y_N\in \fq(\cX)$ be such that $\div_\infty(x)=n_1P_1$,
and $\div_\infty (y_i)=n_i P_1$ for $i=2,\ldots, N$. The fact that
$\nu_1>1$ means that the following matrix has rank two (see Sect.
\ref{s3}) 
   $$
\left( \begin{array}{ccccc}
1 & x^q & y_2^q &\ldots &y_r^q\\
1 & x   & y_2   &\ldots &y_r\\
0 & 1   & D^1_xy_2   &\ldots& D^1_xy_r
\end{array} \right)\, .
    $$ 
In particular, 
   \begin{equation}\label{eq4.33}
y_i^q-y_i= D^1_xy_i(x^q-x) \qquad \text{for}\ \ i=2,\ldots, N\, .
   \end{equation}
   \begin{lemma}\label{lemma4.33} 
   \begin{enumerate}
\item[\rm(1)] $(2g-2)P$ is canonical for any $P\in\cX(\fq);$ i.e.$,$ the
Weierstrass semigroup at such a $P$ is symmetric$;$ 
\item[\rm(2)]  Let $m\in H(P_1)$ such that $m<q+2q_0$. Then $m\le
q+q_0;$
\item[\rm(3)] There exists $g_i\in \fq(\cX)$ such that $
D^1_xy_i=g_i^{\epsilon_2}$ for $i2,\ldots,N.$ Furthermore$,$ 
$\div_\infty(g_i)=\frac{qm_i-q^2}{\epsilon_2}P_1.$
    \end{enumerate}
    \end{lemma}
    \begin{proof} (1) By the identity $2g-2=(2q_0-2)(q+2q_0+1)$ and
(\ref{eq4.12}) we can assume $P=P_1$.  Now the case $i=N$ of 
Eqs. (\ref{eq4.33}) implies $v(dx)=2g-2$ and the result follows since
$v_P(dx)\geq 0$ for $P\neq P_1$.

(2) From (\ref{eq4.32}), $q, q+2q_0$ and $q+2q_0+1\in H(P_1)$. Then the
numbers
   $$
(2q_0-2)q+q-4q_0+j\qquad j=0,\ldots,q_0-2
   $$
are also non-gaps at $P_1$. Therefore, by the symmetry of $H(P_1)$,
   $$
q+q_0+1+j\qquad j=0,\ldots,q_0-2
    $$
are gaps at $P_1$ and the proof follows.

(3) Set $f_i:= D^1_xy_i$. We have 
$D^j_xy_i=(x^q-x)D^j_xf_i+D^{(j-1)}_xf_i$ for $1\le j<q$ by the product
rule applied to (\ref{eq4.33}). Then, $D^j_xf_i=0$ for $1\le 
j<\epsilon_2$, because the matrices
   $$
\left( \begin{array}{ccccc}
1 & x   & y_2   &\ldots &y_N\\
0 & 1   & D^1_xy_2   &\ldots& D^1_xy_N\\
0 & 0   & D^j_xy_2   &\ldots& D^j_xy_N
\end{array} \right),
    \quad 2\le j<\epsilon_2
   $$ 
have rank two (see Sect. \ref{s2.2}). Consequently, as $\epsilon_2$
is a power of two by Lemma \ref{lemma4.32}(2)), from Remark
\ref{rem2.15}(2), $f_i=g_i^{\epsilon_2}$ for some $g_i\in \fq(\cX)$.
Finally, from the proof of (1) we have that $x-x(P)$ is a local parameter
at $P$ if $P\neq P_1$. Then, by the election of the $y_i$'s, $g_i$ has no
pole but in $P_1$, and from (\ref{eq4.33}),
$v(g_i)=-(qm_i-q^2)/\epsilon_2$.
    \end{proof}
    \begin{lemma}\label{lemma4.34} $N=4$ and $\epsilon_2=q_0.$
    \end{lemma}
    \begin{proof} We know  that $N\ge 3$. We claim that $N\ge 4$ otherwise  
we would have $\epsilon_2=2q_0$, $n_1=q$, $n_2=q+2q_0$,
$n_3=q+2q_0+1$, and hence $v(g_2)=-q$ (with $g_2$ being as in Lemma  
\ref{lemma4.33}(3)). Therefore, after some $\fq$-linear transformations,
the case $i=2$ of (\ref{eq4.33}) reads 
   $$ 
y_2^q-y_2=x^{2q_0}(x^q-x)\, .
   $$
Now the function $z:= y_2^{q_0}-x^{q_0+1}$ satisfies
$z^q-z=x^{q_0}(x^q-x)$ and we find that $q_0+q$ is
a non-gap at $P_1$ (cf. \cite[Lemma 1.8]{hansen-sti}). This contradiction
eliminates the possibility $N=3$.

Let $N\geq 4$ and $2\le i\leq N$. By Lemma \ref{lemma4.33}(3)  
$(qn_i-q^2)/\epsilon_2\in H(P_1)$, and since $(qn_i-q^2)/\epsilon_2\ge
n_{i-1}\ge q$, by (\ref{eq4.32}) we have 
  $$
2q_0\ge \epsilon_2 +\epsilon_{N-i}\qquad \mbox{for}\ i=2,\ldots,N-2\, .
  $$
In particular, $\epsilon_2\le q_0$. On the other hand, by Lemma  
\ref{lemma4.33}(2) we must have $n_{N-2}\le q+q_0$ and so, by
(\ref{eq4.32}) we find that $\epsilon_2\ge q_0$; i.e., $\epsilon_2=q_0$. 

Finally we show that $N=4$. $\epsilon_2=q_0$ 
implies $\epsilon_{N-2}\leq q_0$. Since $n_2\le q+q_0$ (cf. Lemma 
\ref{lemma4.33}(2)), by (\ref{eq4.32}), we have $\epsilon_{N-2}\geq q_0$.
Therefore $\epsilon_{N-2}=q_0=\epsilon_2$ so that $N=4$.
    \end{proof}

{\em Proof of Theorem \ref{thm4.31}.} Let $P_1\in \cX(\fq)$ be as above.
By (\ref{eq4.33}), Lemma \ref{lemma4.33}(3) and Lemma \ref{lemma4.34} we
have the following equation   
   $$
y_2^q-y_2=g_2^{q_0}(x^q-x)\ ,
   $$
where $g_2$ has no pole except at $P_1$. Moreover, by (\ref{eq4.32}), 
$n_2=q_0+q$ and so $v(g_2)=-q$ (cf. Lemma \ref{lemma4.33}(3)). Thus 
$g_2=ax+b$ with $a,b\in \fq$, $a\neq 0$, and after some $\fq$-linear
transformations (as those in the proof of Theorem \ref{thm4.21}) the
result follows.
 
   \begin{remark}\label{rem4.31} (i) From the above computations we
conclude that the Suzuki curve $\cS$ is equipped with
a complete, simple and base-point-free $g^4_{q+2q_0+1}$, namely
$\cD_\cS=|(q+2q_0+1)P_0|$, $P_0\in \cS(\fq)$. Such a linear series is an
$\fq$-invariant. The orders of $\cD_\cS$ (resp. the $\fq$-Frobenius
orders) are $0, 1, q_0, 2q_0$ and $q$ (resp. $0, q_0, 2q_0$ and $q$).

(ii) There exists $P_1\in \cS(\fq)$ such that the $(\cD_\cS,P_1)$-orders
are $0,1,q_0+1, 2q_0+1$ and $q+2q_0+1$ (Corollary \ref{cor4.31}). Now 
we show that the above sequence is, in fact, the 
$(\cD_\cS,P)$-orders for each $P\in \cS(\fq)$. To see this, notice that
  $$
\deg(S)=(3q_0+q)(2g-2)+(q+4)(q+2q_0+1)=(4+2q_0)\#\cS(\fq).
  $$ 
Let $P\in \cS(\fq)$. By (\ref{eq4.311}) we conclude that 
$v_P(S^\cD)=\sum_{i=1}^{4}(j_i(P)-\nu_{i-1})=4+2q_0$ and so, by
Proposition \ref{prop3.2}(1) that $j_1(P)=1$, $j_2(P)=q_0+1$,
$j_3(P)=2q_0+1$, and $j_4(P)=q+2q_0+1$.

(iii) Then, by Lemma \ref{lemma4.12} $H(P)$ contains the 
semigroup $H:= \langle q,q+q_0,q+2q_0,q+2q_0+1\rangle$ whenever
$P\in\cS(\fq)$. Indeed $H(P)=H$ since $\#(\mathbb N_0\setminus
H)=g=q_0(q-1)$ (this can be proved as in the claim in the proof of Lemma
\ref{lemma4.32}(1); see also \cite[Appendix]{hansen-sti}). 

(iv) We have  
   $$
\deg(R)=\sum_{i=0}^{4}\epsilon_i(2g-2)+5(q+2q_0+1)=(2q_0+3)\#\cS(\fq)\, ,
   $$
and $v_P(R)=2q_0+3$ for $P\in\cS(\fq)$ as follows from (i), (ii) 
and Sect. \ref{s2.2}. Therefore the set of $\cD_\cS$-Weierstrass points of
$\cS$ is equal to $\cS(\fq)$. In particular, the $(\cD,P)$-orders for
$P\not\in \cS(\fq)$ are $0, 1, q_0, 2q_0$ and $q$. 

(v) We can use the above computations to obtain information on orders
for the canonical morphism on $\cS$. By using the fact that
$(2q_0-2)\cD_\cS$ is canonical (cf. Lemma 
\ref{lemma4.33}(1)) and (iv), we see that the set $
\{a+q_0b+2q_0c+qd: a+b+c+d \le 2q_0-2\}$
is contained in the set of orders of $\cK_\cS$ at non-rational points. (By
considering first order differentials on $\cS$, similar computations were 
obtained in \cite[Sect. 4]{garcia-sti}.) 

(vi) Finally, we remark that $\cS$ is non-classical for the canonical
morphism: We have two different proofs for
this fact: loc. cit. and Corollary \ref{cor4.13}(3).
    \end{remark}
    \begin{remark} (A. Cossidente)\label{rem4.32} Recall that an {\em
ovoid} in $\P^N(\fq)$ is a set of points $P$ no three of which are
collinear and such that for each $P$ the union of the tangent lines
at $P$ is a hyperplane; see
\cite{h}. We are going to related the Suzuki-Tits ovoid $\cO$ in
$\P^4(\fq)$ with the $\fq$-rational points of the Suzuki curve $\cS$.

It is known that any ovoid in $\P^4(\fq)$ that contains the point
$(0:0:0:0:1)$ can be defined by
  $$ 
\{(1:a:b:f(a,b):af(a,b)+b^2): a, b \in \fq\}\cup\{(0:0:0:0:1)\},
  $$ 
where $f(a,b):=a^{2q_0+1}+b^{2q_0}$; cf. \cite{tits},
\cite[p.3]{penttila-williams}. 

Let $\phi=(1:x:y:z:w)$ be the morphism associated to $\cD_\cS$ such that 
$\div_\infty(x)=qP_0$, $\div_\infty(y)=(q+q_0)P_0$,
$\div_\infty(z)=(q+2q_0)P_0$ and
$\div_\infty(w)=q+2q_0+1$; see Remark \ref{rem4.31}(iii).
  
   \begin{claim*} $\cO=\phi(\cS(\fq)).$
   \end{claim*}
Indeed we have $\phi(P_0)=(0:0:0:0:1)$; in addition the coordinates of
$\phi$ can be choosen such that $y^q-y=x^{q_0}(x^q-x)$, $z:=
x^{2q_0+1}+y^{2q_0}$, and $w:=
xy^{2q_0}+z^{2q_0}=xy^{2q_0}+x^{2q+2q_0}+y^{2q}$ (see 
\cite[Sect. 1.7]{hansen-sti}). For $P\in \cS(\fq)\setminus\{P_0\}$ set
$a:=x(P)$, $b:=y(P)$, and $f(a,b):= z(a,b)$. Then $w(a,b)=af(a,b)+b^2$ and
the claim follows.
   \end{remark}

   \begin{remark}\label{rem4.33} The morphism $\phi$ in the previous
remark is an embedding. To see this, as $j_1(P)=1$ for any $P\in\cS$ (
Remarks \ref{rem4.31}(ii)(iv)), it is enough to show that 
$\phi$ is injective. We have
   \begin{equation}\label{eq4.34}
(q+2q_0+1)P_0\sim q\fro^2(P)+2q_0\fro(P)+P
   \end{equation}
so that the points $P\in\cS$ where $\phi$ could not
be injective satisfy either $P\not\in \cS(\fq)$, or $ \fro^3(P)=P$ or
$\fro^2(P)=P$. Now from the Zeta function of $\cS$ one sees that
$\#\cS(\mathbb F_{q^3})=\#\cS(\mathbb F_{q^2})=\#\cS(\fq)$, and the remark
follows. 
   \end{remark}

   \begin{remark}\label{rem4.34} From the claim in Remark \ref{rem4.32}, 
(\ref{eq4.34}) and \cite{henn} we have
   $$
{\rm Aut}_{\bar\fq}(\cS)={\rm Aut}_{\fq}(\cS)\cong
\{A\in {\rm PGL}(5,q): A\cO=\cO\}\, .
   $$
   \end{remark}

\section{Plane arcs}\label{s5}

In this section we show how to apply Sections \ref{s2} and \ref{s3} to
study the size of plane arcs. The approach is from 
Hirschfeld-Korchm\'aros \cite{hk1}, \cite{hk2} and Voloch \cite{voloch1},
\cite{voloch2}. Our exposition follows \cite{perugia1}.

A {\em $k$-arc} in $\P^2(\fq)$ is a set $\cK$ of $k$ points no three of
which are collinear. It is {\em complete} if it is not properly
contained in 
another arc. For a given $q$, a basic problem in Finite Geometry is to
find the values of $k$ for which a complete $k$-arc exists. Bose
\cite{bose} showed that
   $$
k\le m(2,q):=\begin{cases}
q+1 & \text{if $q$ is odd}\, ,\\
q+2 & \text{otherwise}\, .
\end{cases}
   $$  
For $q$ odd the bound $m(2,q)$ is attained if and only if $\cK$ is an
irreducible conic \cite{segre}, \cite[Thm. 8.2.4]{h}. For $q$ even the
bound is attained by the union of an irreducible conic and its nucleus,
and not every $(q+2)$-arc arises in this way; see \cite[Sect. 8.4]{h}. Let
$m'(2,q)$ denote the second largest size that a complete arc in
$\P^2(\fq)$ can have. Segre \cite{segre}, \cite[Sect. 10.4]{h} showed that 
    \begin{equation}\label{eq5.1}
m'(2,q)\le\begin{cases}
q-\frac{1}{4}\sq+\frac{7}{4} & \text{if $q$ is odd},\\
q-\sq+1                      & \text{otherwise}.
\end{cases}
    \end{equation}
Besides small $q$, namely $q\le 29$ \cite{chao-kaneta}, \cite{h},
\cite{h-storme}, the only case where $m'(2,q)$ has been determined is
for $q$ an even square.  Indeed, for $q$ square, examples of complete
$(q-\sq+1)$-arcs \cite{boros}, \cite{coss}, \cite{ebert},
\cite{fisher-h-thas}, \cite{kestenband} show that
  \begin{equation}\label{eq5.2}
m'(2,q)\ge q-\sq+1\, ,
   \end{equation}
and so the bound (\ref{eq5.1}) for an even $q$ square is sharp. This
result has been recently extended by Hirschfeld and Korchm\'aros
\cite{hk3} who showed that the third largest size that a complete arc
can have is upper bounded by $q-2\sq+6$. 

If $q$ is not a square, Segre's bounds were notably improved by
Voloch \cite{voloch1}, \cite{voloch2}. 

If $q$ is odd, Segre's bound was slightly improved to $m'(2,q)\le
q-\sq/4+25/16$ by Thas \cite{thas}. If $q$ is an odd square and large
enough, Hirschfeld and Korchm\'aros \cite{hk2} significantly improved the
bound to
   \begin{equation}\label{eq5.3}
m'(2,q)\le q-\frac{1}{2}\sq+\frac{5}{2}\, . 
   \end{equation}
Inequalities (\ref{eq5.2}) and (\ref{eq5.3}) suggest the following
problem, which seems to be difficult and has remained open since the 60's.
   \begin{problem}\label{prob5.1}
For $q$ an odd square, is it true that $m'(2,q)=q-\sq+1$?
   \end{problem}
The answer is negative for $q=9$ and affirmative for $q=25$ 
\cite{chao-kaneta}, \cite{h}, \cite{h-storme}. So Problem \ref{prob5.1} is
indeed open for $q\ge 49$.

\subsection{B. Segre's fundamental theorem: Odd case}\label{s5.1} We
recall a fundamental theorem of Segre which is the link between arcs and
curves.  

Let $\cK$ be an arc in $\P^2(\fq)$. Segre associates to $\cK$ a plane
curve $\cC$ in the dual plane of $\P^2(\bar\fq)$. This curve is defined
over $\fq$ and it is called {\em the envelope of} $\cK$. For $P\in
\P^2(\bar\fq)$, let $\ell_P$ denote the corresponding line in the dual
plane. A line $\ell$ in $\P^2(\fq)$ is called an {\em $i$-secant} of $\cK$
if $\# \cK\cap\ell=i$. The following result summarizes the main properties
of $\cC$ for the odd case.

    \begin{theorem} {\rm (B. Segre \cite{segre}, \cite[Sect. 10]{h})}
\label{thm5.11} If $q$ is odd, then the following statements hold:
    \begin{enumerate}

\item[\rm(1)] The degree of $\cC$ is $2t$, with $t=q-k+2$ being the number
of $1$-secants through a point of $\cK$.

\item[\rm(2)] All $kt$ of the $1$-secants of $\cK$  belong to $\cC$.

\item[\rm(3)] Each $1$-secant $\ell$ of $\cK$ through a point $P\in\cK$ is
counted twice in the intersection of $\cC$ with $\ell_P$; i.e., 
$I(\cC,\ell_P;\ell)=2$.

\item[\rm(4)] The curve $\cC$ contains no $2$-secant of $\cK$.

\item[\rm(5)] The irreducible components of $\cC$ have multiplicity at
most two, and $\cC$ has at least one component of multiplicity one.

\item[\rm(6)] For $k>(2q+4)/3$, the arc $\cK$ is incomplete if and only if
$\cC$ admits a linear component over $\fq$. For $k>(3q+5)/4$, the arc
$\cK$ is a conic if and only if it is complete and $\cC$ admits a
quadratic
component over $\fq$.
   \end{enumerate}
   \end{theorem} Next we show some properties of $\cC$. Recall that a
non-singular point $P$ of a plane curve $\cA$ is called an {\em inflexion
point of $\cA$} if $I(\cA,\l;P)>2$, with $\ell$ being the tangent line of
$\cA$ at $P$.
  \begin{definition*} A point $P_0$ of $\cC$ is called {\em
special} if the following conditions hold:
   \begin{enumerate} 
\item[(i)] it is non-singular; 
\item[(ii)] it is $\fq$-rational;
\item[(iii)] it is not an inflexion point of $\cC$.
   \end{enumerate}
Then, by (i), a special point $P_0$ belongs to an 
unique irreducible component of the envelope which will be called {\em the
irreducible envelope} associated to $P_0$ or {\em an irreducible envelope
of $\cK$}. 
   \end{definition*}
   \begin{lemma}\label{lemma5.11} Let $\cC_1$ be an irreducible envelope
of $\cK$. Then
   \begin{enumerate}
\item[\rm(1)] $\cC_1$ is defined over $\fq;$
\item[\rm(2)]  if $q$ is odd and the $k$-arc $\cK$, with $k>(3q+5)/4$, is 
complete and different from a conic, then the degree of $\cC_1$ is at
least three.
   \end{enumerate}
   \end{lemma}
   \begin{proof} (1) Let $\cC_1$ be associated to $P_0$, let $\mathbf \Phi$ be
the Frobenius morphism (relative to $\fq$) on the dual plane of
$\P^2(\bar\fq)$, and suppose that $\cC_1$ is not defined over $\fq$. Then,
since the envelope is defined over $\fq$ and $P_0$ is $\fq$-rational,
$P_0$ would belong to two different components of the envelope, namely
$\cC_1$ and ${\mathbf \Phi}(\cC_1)$. This is a contradiction because the
point is non-singular. 

(2) This follows from Theorem \ref{thm5.11}(6). 
    \end{proof}
The next result will show that special points do exist provided that $q$ is
odd and the arc is large enough.
     \begin{proposition}\label{prop5.11} Let $\cK$ be an arc in
$\P^2(\fq)$ of size $k$ such that $k>(2q+4)/3$. If $q$ is odd, then the
envelope $\cC$ of $\cK$ has special points.
     \end{proposition}
  \begin{remark}\label{rem5.11} The hypothesis $k>(2q+4)/3$ in the
proposition is equivalent
to $k>2t$, with $t=q-k+2$. Also, under this hypothesis, the envelope
$\cC$ is uniquely determined by $\cK$, see \cite[Thm. 10.4.1(i)]{h}.
  \end{remark}
To prove Proposition \ref{prop5.11} we need the following lemma.
      \begin{lemma}\label{lemma5.12} 
Let $\cA$ be a plane curve defined over $\bar\fq$ and suppose that it 
has no multiple components. Let $\alpha$ be the degree of $\cA$ and $s$
the number of its singular points. Then,
   $$
s\le \binom{\alpha}{2}\, ,
   $$
and equality holds if $\cA$ consists of $\alpha$ lines no three 
concurrent.  
      \end{lemma}
      \begin{proof}
That a set of $\alpha$ lines no three concurrent satisfies the bound is
trivial. Let $G=0$ be the equation of $\cA$, let $G=G_1\ldots G_r$ be the 
factorization of $G$ in $\bar\fq[X,Y]$, and let $\cA_i$ be the curve given by
$G_i=0$. For 
simplicity we assume $\alpha$ even, say $\alpha=2M$. Setting
$\alpha_i:=\deg(G_i)$, 
$i=1, \ldots, r$ and $I:=\sum_{i=1}^{r-1}\alpha_i$ we have
$\alpha_r=2M-I$. The
singular points of $\cA$ arise from the singular points of each component
and from the points in $\cA_i \cap \cA_j $, $i \neq j$. Recall that an
irreducible plane curve of degree $d$ has at most $\binom{d-1}{2}$ singular
points, and that $\#\cA_i \cap \cA_j \le a_ia_j$, $i \neq j$ (B\'ezout's
Theorem). So
   \begin{equation*}
\begin{split}
s  & \le  \sum_{i=1}^{r-1}\binom{\alpha_i
-1}{2}+\binom{2M-I-1}{2}+\sum_{1\le
i_1<i_2\le r-1}\alpha_{i_1}\alpha_{i_2}+\sum_{i=1}^{r-1}(2M-I)\alpha_i\\
 & =
\sum_{i=1}^{r-1}\frac{\alpha_i^2-3\alpha_i+2}{2}+\frac{4M^2-4MI+I^2-6M+3I+2}{2}
+\sum_{1\le i_1 < i_2 \le r-1}\alpha_{i_{1}}\alpha_{i_{2}}+(2M-I)I\\
 & =
\frac{1}{2}[
\sum_{i=1}^{r-1}\alpha_{i}^{2}-3I+2(r-1)+4M^2-4MI+I^2-6M+3I+2+\\
{} &\quad \ 2\sum_{1\le i_1 
< i_2 \le r-1}\alpha_{i_{1}}\alpha_{i_{2}}+4MI-2I^2]\\
 & \le 2 M^2 - 3M + \alpha = 2M^2 -M\, .
\end{split}
    \end{equation*}
          \end{proof}
          \begin{proof} ({\em Proposition \ref{prop5.11}}) Let $F=0$ be
the equation of $\cC$ over $\fq$. By Theorem \ref{thm5.11}(5), $F$ admits
a factorization in $\bar\fq[X,Y,Z]$ of type
  $$
G_1\ldots G_r H_1^2\ldots H_s^2\, ,
  $$
with $r\ge 1$ and $s\ge 0$. Let $\cA$ be the plane curve given by 
  $$
G:=G_1\ldots G_r=0\, .
  $$
Then $\cA$ satisfies the hypothesis of Lemma \ref{lemma5.12} and it has  
even degree by Theorem \ref{thm5.11}(1). From Theorem \ref{thm5.11}(3)  
and B\'ezout's theorem,  for each line $\ell_P$ (in the dual plane)  
corresponding to a point $P\in \cK$, we have  
$$
\#(\cA \cap \ell_P) \ge M\, ,
$$
where $2M = \deg(G)$, and so at least $kM$ points corresponding 
to unisecants of $\cK$ belong to $\cA$. Since $k>2t$ (see
Remark \ref{rem5.11}) and $2t\ge 2M$, then $kM>2M^2$ and from Lemma 
\ref{lemma5.11} we have that at least one
of the unisecant points in $\cA$, says $P_0$, is non-singular. Suppose
that $P_0$ passes through $P\in \cK$. The point $P_0$ is clearly
$\fq$-rational and 
$P_0$ is not a point of the curve of equation $H=0$: 
otherwise $I(P_0, \cC\cap\ell_P) > 2$ (see Theorem \ref{thm5.11}(3)).
Then,  
$I(P_0, \cC\cap\ell_P)=I(P_0,\cA\cap\ell_P)=2$ and so $\ell_P$ is
the tangent of $\cC$ at $P_0$. Therefore $P_0$ is not an inflexion 
point of $\cC$, and the proof of Proposition \ref{prop5.11} is complete. 
             \end{proof}
Let $\cC_1$ be an irreducible envelope associated to a special point
$P_0$, and 
   $$
\pi: \cX\to \cC_1\, ,
   $$
the non-singular model of $\cC_1$. Then by Lemma
\ref{lemma5.11}(1) we can
assume that $\cX$ and $\pi$ are both defined over $\fq$. In particular,
the linear series $\Sigma_1$ cut out by lines of 
$\P^2(\bar\fq)^*$ on $\cX$ is $\fq$-rational. Also, there is just one
point 
$\tilde P_0\in \cX$ such that $\pi(\tilde P_0)=P_0$. 
   \begin{lemma}\label{lemma5.13} Let $q$ be odd. Then,
\begin{enumerate}
  \item[\rm(1)] the $(\Sigma_1,\tilde P_0)$-orders are $0, 1, 2;$
  \item[\rm(2)] the curve $\cX$ is classical with respect to $\Sigma_1$.
\end{enumerate}
    \end{lemma}
\begin{proof} (1) follows from the proof of Proposition
\ref{prop5.11} while (2) from (1) and Corollary \ref{cor2.21}(1).
   \end{proof}
   \begin{remark}\label{rem5.12} The hypothesis $q$ odd in Lemma
\ref{lemma5.13} (as well as in Proposition \ref{prop5.11}) 
is necessary. In fact, from \cite{fisher-h-thas} and \cite{thas} follow
that the envelope associated to the cyclic $(q-\sq+1)$-arc, with $q$ an even
square, is irreducible and $\fq$-isomorphic to the curve of equation 
$XY^{\sq}+X^{\sq}Z+YZ^{\sq}=0$. It is not difficult to see that this curve
is $\bar\fq$-isomorphic to the Hermitian curve $\cH$ in Example
\ref{ex3.2} (see e.g. \cite[p. 4711]{ckt2}) so that it is $\Sigma_1$ 
non-classical.
   \end{remark}
Next consider the following sets:
\begin{align*}
\cX_1(\fq):= & \{P\in \cX: \pi(P)\in \cC_1(\fq)\}\, ,\\
\cX_{11}(\fq):= & \{P\in \cX_1(\fq): j^1_2(P)=2j^1_1(P)\}\, ,\\
\cX_{12}(\fq):= & \{P\in \cX_1(\fq): j^1_2(P)\neq 2j^1_1(P)\}\, ,
\end{align*}
and the following numbers:
   \begin{equation}\label{eq5.11} 
M_q=M_q(\cC_1):=\sum_{P\in \cX_{11}(\fq)} j^1_1(P)\, , \qquad
M_q'=M_q'(\cC_1):=\sum_{P\in \cX_{12}(\fq)} j^1_1(P)\, ,
   \end{equation}
where $0<j^1_1(P)<j^1_2(P)$ denotes the $(\Sigma_1,P)$-order sequence. We
have that 
$$
M_q+M_q'\ge \#\cX_1(\fq) \ge \#\cX(\fq)\qquad\text{and}\qquad
\#\cX_1(\fq)\ge \#\cC_1(\fq)\, .
$$
    \begin{proposition}\label{prop5.12} Let $\cK$ be an arc of size $k$
and $d$ the degree of an irreducible envelope of 
$\cK$. For $M_q$ and $M_q'$ as above we have
$$
2M_q+M_q'\ge kd\, .
$$
    \end{proposition}
To prove this proposition we first prove the following lemma.
	\begin{lemma}\label{lemma5.14}
Let $\cK$ be an arc and $\cC_1$ an irreducible envelope of $\cK$.  Let $Q
\in \cK$ and $\cA_Q$ be the set of points of $\cC_1$ corresponding to
unisecants of $\cK$ passing through $Q$. Let $u:=\#\cA_Q$ and $v$ be the
number of points in $\cA_Q$ which are non-singular and inflexion points
of $\cC_1$. Then
  $$
2(u-v)+v \ge d\, ,
  $$
where $d$ is the degree of $\cC_1$.
	\end{lemma}
       \begin{proof} 
Let $P'\in \cA_Q$.  Suppose that it is  non-singular and an inflexion point
of $\cC_1$. Then, from Theorem \ref{thm5.11}(3) and the definition of
$\cA_Q$, we have that $\ell_Q$ is not the tangent line of $\cC_1$ at $P'$,
i.e. we have that $I(P', \cC_1\cap \ell_Q)=1$. Now suppose that $P'$ is
either singular or a non-inflexion point of $\cC_1$. Then from Theorem
\ref{thm5.11}(3) we have $I(P',\cC_1\cap \ell_Q)\le 2$ and the result
follows from B\'ezout's theorem applied to $\cC_1$ and $\ell_Q$. 
       \end{proof}
{\em Proof of Proposition \ref{prop5.12}.} Let $Q\in \cK$ and $\cA_Q$ be
as in Lemma \ref{lemma5.14}. Set 
\begin{align*}
\cY_Q  & :=\{P \in \cX_1(\fq) :  \pi(P) \in \cA_Q \}\, ,\\
\intertext{and}
m(Q)   & :=2\sum_{P \in \cX_{11}(\fq) \cap \cY_Q}j^1_1(P)+
\sum_{P \in \cX_{12}(\fq) \cap \cY_Q}j^1_1(P)\, .
\end{align*}
We claim that $m(Q)\ge d$. Indeed, this claim implies the proposition
since, from Theorem \ref{thm5.11}(4),
$$
\cY_Q \cap \cY_{Q_1} = \emptyset\qquad 
\text{whenever}\qquad Q \neq Q_1\, .
$$
To prove the claim we distinguish four types of points in $\cY_Q$, namely
\begin{align*}
\cY_Q^1:= & \{P\in \cY_Q: \text{$\pi(P)$ is non-singular and non-
inflexion point of $\cC_1$}\}\, , \\
\cY_Q^2:= & \{P\in \cY_Q:  \text{$\pi(P)$ is a non-singular 
inflexion point of $\cC_1$}\}\, ,\\
\cY_Q^3:= & \{P\in \cY_Q: \text{$\pi(P)$ is a singular
point of $\cC_1$ such that $\# \pi^{-1}(\pi(P))=1$}\}\, ,\\
\cY_Q^4:= & \{P\in \cY_Q: \text{$\pi(P)$ is a singular 
point of $\cC_1$ such that $\#\pi^{-1}(\pi(P)) >1$}\}\, .
\end{align*}
Observe that $\cY_Q^1 \subseteq \cX_{11}(\fq)$ and so
$$
m(Q)\ge 2\sum_{P\in \cY_Q^1}j^1_1(P)+
\sum_{P\in \cY_Q^2}j^1_1(P)+\sum_{P\in \cY_Q^3}j^1_1(P)+\sum_{P\in
\cY_Q^4}j^1_1(P)\, .
$$
Since $j^1_1(P) > 1$ for all $P\in \cY_Q^4$, the above inequality becomes
$$
m(Q)\ge 2\#\cY_Q^1+2\#\cY_Q^4+\#\cY_Q^3+\#\cY_Q^2\, .
$$ 
Therefore, as to each singular non-cuspidal point of
$\cC_1$ in $\cA_Q$ corresponds at least two points in $\cY_Q^3$, it
follows that
\begin{align*}
m(Q) & \ge 2\#\{P' \in  \cA_Q : \text{$P'$ is either singular or not 
an inflexion point of $\cC_1$} \}+\\ 
{} &\quad  \#\{P' \in  \cA_Q : \text{$P'$ is a nonsingular inflexion
point of $\cC_1$}\}\, .
\end{align*}
Then the claim follows from Lemma \ref{lemma5.14} and the proof of
Proposition \ref{prop5.12} is complete. 

\subsection{The work of Hirschfeld, Korchm\'aros and Voloch}\label{s5.2}

Throughout the whole sub-section we fix the following notation:
    \begin{itemize}
\item $q$ is a power of an odd prime $p$;

\item $\cK$ is a complete arc of size $k$ such that $(3q+5)/4<k\le
m'(2,q)$; therefore the degree of any irreducible envelope of $\cK$ is 
at least three by Theorem \ref{thm5.11}(6);

\item $P_0$ is a special point of the envelope $\cC$ of $\cK$
and the plane curve $\cC_1$ of degree $d$ is an irreducible envelope
associated to $P_0$; 

\item $\pi: \cX\to \cC_1$ is the normalization of $\cC_1$ which is
defined over $\fq$; as a matter of terminology, $\cX$ will be also called
an irreducible envelope of $\cK$.

\item $\tilde P_0$ is the only point in $\cX$ such that
$\pi(\tilde P_0)=P_0$; $g$ is the genus of $\cX$ (so that $g\le
(d-1)(d-2)/2$);

\item The symbols $M_q$ and $M_q'$ are as in Sect. \ref{s5.1};

\item $\Sigma_1$ is the linear series $g^2_d$ cut out by lines of
$\P^2(\bar\fq)^*$ on $\cX$; $\Sigma_2$ is the linear
series $g^5_{2d}$ cut out by conics of $\P^2(\bar\fq)^*$ on $\cX$; then 
$\Sigma_2=2\Sigma_1$. Notice that $\dim(\Sigma_2)=5$ because $d\geq 3$ and
that $\Sigma_1$ and $\Sigma_2$ are base-point-free;

\item $S$ is the $\fq$-Frobenius divisor associated to $\Sigma_2$;

\item $j_5(\tilde P_0)$  is the 5th positive 
$(\Sigma_2,\tilde P_0)$-order;  $\epsilon_5$ is the 5th positive 
 $\Sigma_2$-order; 
$\nu_4$ is the 4th positive $\fq$-Frobenius order of $\Sigma_2$. 
    \end{itemize}

We apply the results in Sects. \ref{s2} and \ref{s3} to $\Sigma_1$ and
$\Sigma_2$. We have already noticed that the $(\Sigma_1,\tilde
P_0)$-orders, as well as the $\Sigma_1$-orders, are 0,1 and 2; see Lemma
\ref{lemma5.13}. Then, the $(\Sigma_2, \tilde P_0)$-orders are 0,1,2,3,4
and $j_5(\tilde P_0)$, with $5\le j_5(\tilde P_0)\le 2d$, and the
$\Sigma_2$-orders are 0,1,2,3,4 and $\epsilon_5$ with $5\le \epsilon_5\le
j_5(\tilde P_0)$.

Then, we compute the $\fq$-Frobenius orders of $\Sigma_2$. We apply
Proposition \ref{prop3.2}(1) to $\tilde P_0$ and infer that this sequence
is 0,1,2,3 and $\nu_4$, with 
   $$ \nu_4\in\{4,\epsilon_5\}\, . 
   $$ 
Therefore 
    \begin{align*}
\deg(S) & =(6+\nu_4)(2g-2)+(q+5)2d\, ,\\
\intertext{and}
v_P(S)& \ge 5j^2_1(P), \qquad \text{for each $P\in \cX_1(\fq)$}\, ,
    \end{align*}
where $j^2_1(P)$ stands for the first positive $(\Sigma_2,P)$-order. 
   \begin{claim*} 
$j^2_1(P)$ equals $j_1^1(P)$ (the first positive $(\Sigma_1,P)$-order).
   \end{claim*}
   \begin{proof} Let $\Sigma_1=\{E+\div(f):f\in
\Sigma_1'\setminus\{0\}\}$. From Sect. \ref{s2.2} we can assume that
$\Sigma_1'=\langle 1,x,y\rangle$ where 
   \begin{equation*}
j_1^1(P)=v_P(E)+v_P(x)\qquad\text{and}\qquad 
j_2^1(P)=v_P(E)+v_P(y)\, .\tag{$*$}
   \end{equation*}
Now $\Sigma_2=\{2E+\div(f):f\in \Sigma_2'\setminus\{0\}\}$, where
$\Sigma_2'=\langle
1,x,y,xy,x^2,y^2\rangle$, and there exists $f\in \Sigma_2'$ such that 
  $$
j_1^2(P)=v_P(2E)+v_P(f)\, .
  $$
Let $f=a_0+a_1x+a_2y+a_3x^2+a_4xy+a_5y^2$. From Lemma \ref{lemma1.22},
  $$
v_P(2E)=-{\rm min}\{v_P(1),v_P(x),v_P(y),v_P(x^2),v_P(xy),v_P(y^2)\}\, .
  $$
Suppose that $0\leq v_P(x)$ and $0\leq v_P(y)$. Then $v_P(2E)=0$ so that
$v_P(f)=j_1^2(P)>0$ and hence $a_0=0$. Then the result follows from $(*)$.
Now suppose that $0>v_P(x)$ or $0>v_P(y)$. Then $v_P(2E)<0$ and hence
$a_i\neq 0$ for some $i\in\{1,\ldots,5\}$. Then the result follows from
$(*)$ and the fact that $v_P(f)\geq {\rm min}\{v_P(x),
v_P(y),v_P(x^2),v_P(xy),v_P(y^2)\}$.
   \end{proof}
We then have 
   $$
\deg(S)\ge 5(M_q+M_q')\, ,
   $$
where $M_q$ and $M_q'$ were defined in (\ref{eq5.11}). 
   \begin{proposition}\label{prop5.21}
Let $\cK$ be a complete arc of size $k$ such that $(3q+5)/4< k\le
m'(2,q)$. Then
$$
k\le \min\{q-\frac{1}{4}\nu_4+\frac{7}{4}\, ,\   
\frac{28+4\nu_4}{29+4\nu_4}q+\frac{32+2\nu_4}{29+4\nu_4}\}\, ,
$$
where $\nu_4$ is the 4th positive $\fq$-Frobenius order of the linear 
series $\Sigma_2$ defined on an irreducible envelope of $\cK$.
    \end{proposition}
    \begin{proof} From the computations above and Proposition 
\ref{prop5.12}, 
   $$
\deg(S)=(6+\nu_4)(2g-2)+(q+5)2d\ge 5(M_q+M_q')\ge \frac{5}{2}kd\, .
   $$
Now $d(d-3)\ge 2g-2$ and $d\le 2t=2(q+2-k)$ (Theorem \ref{thm5.11}(1)).
Then $k(29+\nu_4)\le (28+4\nu_4)q+(32+2\nu_4)$. On the other
hand, $\nu_4\le j_5(\tilde P_0)-1\le 2d-1$ (Proposition \ref{prop3.2}(1))
and hence $k\le q-\nu_4/4+7/4$. 
    \end{proof}
Next we consider separately the cases $\nu_4=4$ and $\nu_4=\epsilon_5$.

{\bf Case $\nu_4=4$.} In this case, the corresponding irreducible 
envelope will be called {\em Frobenius 
classical}. Proposition \ref{prop5.21} becomes the following.
    \begin{corollary}\label{cor5.21}  
Let $\cK$ be a complete arc of size $k$ such that $(3q+5)/4<k\le
m'(2,q)$. Suppose that $\cK$ admits a Frobenius classical 
irreducible envelope. Then 
$$
k\le \frac{44}{45}q+\frac{40}{45}\, .
$$     
    \end{corollary}
The bound in the corollary holds in the following cases:

(A) (Voloch \cite{voloch2}) Whenever $q=p$ is an odd prime;

(B) (Giulietti \cite{giu}) The arc is cyclic of Singer type whose size $k$
satisfies $2k\not\equiv -2,1,2,4 \pmod{p}$, where $p>5$.
 
For the sake of completeness let us prove (A): Let 
$\cC_1$ be an irreducible envelope of
$\cK$ and $d$ the degree of $\cC_1$. If $p< 2d$, then $p<4t=4(p+2-k)$
so that $k<(3p+8)/4$ and the result follows. So let $p\ge 2d$. Then from
Remark \ref{rem3.2} we have that $\cC_1$ is Frobenius classical and
(A) follows from Proposition \ref{prop5.21}. 

Next we show that, for $q$ square and $k=m'(2,q)$, Corollary  
\ref{cor5.21} can only hold for $q$ small.
    \begin{corollary}\label{cor5.22} Let $\cK$ be an arc of size
$m'(2,q)$ and suppose that $q$ is a square. Then,
   \begin{enumerate}
\item[\rm(1)] if $q>9$, $\cK$ has irreducible envelopes;
\item[\rm(2)] if $q>43^{2}$, any irreducible envelope of $\cK$ is
Frobenius
non-classical. 
   \end{enumerate}
    \end{corollary}
    \begin{proof} (1) As we mentioned in (\ref{eq5.2}), $m'(2,q)\ge
q-\sq+1$. Since  $q-\sq+1>(2q+4)/3$ for $q>9$, (1) follows from
Proposition \ref{prop5.11}.

(2) If existed a Frobenius classical irreducible envelope of $\cK$,
then from Lemma \ref{lemma5.21} and (\ref{eq5.2}) we would have
$$
q-\sq+1\le m'(2,q)\le 44q/45+40/45\, .
$$
so that $q\le 43^2$. 
   \end{proof}

{\bf Case $\nu_4=\epsilon_5$.} Here, from Lemma \ref{lemma3.3} we
have that $p$ divides $\epsilon_5$. More precisely we have the following
result. 
    \begin{lemma}\label{lemma5.21} 
Either $\epsilon_5$ is a power of $p$ or $p=3$ and $\epsilon_5=6$.
    \end{lemma}
    \begin{proof} We can assume $\epsilon_5>5$. If $\epsilon_5$ is not a
power of $p$, by the $p$-adic criterion (Lemma \ref{lemma2.31}) we have 
$p\le 3$ and $\epsilon=6$.
    \end{proof} 
From Proposition \ref{prop5.21}, the case $\nu_4=\epsilon_5=6$ provides
the following bound:
  \begin{lemma}\label{lemma5.22}
Let $\cK$ be a complete arc of size $k$ such that $(3q+5)/4<k\le
m'(2,q)$. Suppose that $\cK$ admits an irreducible envelope such that 
$\nu_4=\epsilon_5=6$. Then $p=3$ and 
$$
k\le \frac{52}{53}q+\frac{44}{53}\, .
$$
  \end{lemma}
As in the case $\nu_4=4$, for $q$ an even power of 3 and
$k=m'(2,q)$ the case $\nu_4=\epsilon_5=6$ occur only for $q$ small. More
precisely, we have the following result.  
  \begin{corollary}\label{cor5.23} Let $\cK$ be an arc of size $m'(2,q)$. 
Suppose that $q$ is an even power of $p$ and that $\cK$ admits an
irreducible envelope with $\nu_4=\epsilon_5=6$. Then $p=3$ and $q\le 3^6$.
  \end{corollary}
   \begin{proof} From the $p$-adic criterion (Lemma \ref{lemma2.31}),
$p=3$. Then from Proposition \ref{prop5.21} and (\ref{eq5.2}) we have
   $$
q-\sq+1\le m'(2,q)\le 52q/53+44/53\, ,
   $$
and the result follows.
  \end{proof} 
From now on we assume 
$$
\nu_4=\epsilon_5=\text{a power of $p$}\, . 
$$
Then, the bound 
   \begin{equation}\label{eq5.21}
k\le q-\frac{1}{4}\nu_4+\frac{7}{4}
   \end{equation}
in Proposition \ref{prop5.21} and
Segre's bound (\ref{eq5.1}) provide motivation to consider three cases
according as $\nu_4>\sq$, $\nu_4< \sq$, or $\nu_4=\sq$. 

{\bf Case $\nu_4>\sq$.} Since $\nu_4$ is a power of $p$, here we
have that $\nu^2 \ge pq$ and so from (\ref{eq5.21}) the following
holds: 
   \begin{lemma}\label{lemma5.23} Let $\cK$ be a complete arc of size $k$
such that $(3q+5)/4<k\le m'(2,q)$. Suppose that $\cK$ admits an
irreducible envelope such that $\nu_4$ is a power of $p$ and
that $\nu_4>\sq$. Then
$$
k\le \begin{cases}
q-\frac{1}{4}\sqrt{pq}+\frac{7}{4} & \text{if $q$ is not a square}\, ,\\
q-\frac{1}{4}p\sq+\frac{7}{4} & \text{otherwise}\, .
     \end{cases}
$$
   \end{lemma}
If $q$ is a square and $k=m'(2,q)$, then $\nu_4>\sq$ can only occur in
characteristic 3:
   \begin{corollary}\label{cor5.24} Let $\cK$ be an arc of size $m'(2,q)$.
Suppose that $q$ is an even power of $p$ and that $\cK$ admits an
irreducible envelope with $\nu_4$ a power of $p$ and $\nu_4>\sq$.
Then $p=3$, $\nu_4=3\sq$, and 
  $$
k\le q-\frac{3}{4}\sq+\frac{7}{4}\, .
  $$
   \end{corollary}
    \begin{proof} From Lemma \ref{lemma5.23} and (\ref{eq5.2}) follow that
$\sq(p-4)\le 3$ and so that $p=3$. From $\nu_4\le
2d-1$ and $2d\le 4t=4(q+2-m'(2,q))\le 4\sq+4$ we have that $\nu_4\le
4\sq+3$ and it follows the assertion on $\nu_4$. The bound on $k$ follows
from Lemma \ref{lemma5.23}.
   \end{proof}
{\bf Case $\nu_4<\sq$.} Let 
$$
F(x):= (2x+32-q)/(4x+29)\, .
$$
Then the bound 
$$
k\le \frac{28+4\nu_4}{29+4\nu_4}q+\frac{32+2\nu_4}{29+4\nu_4}
$$ 
in Proposition \ref{prop5.21} can be written as 
   \begin{equation}\label{eq5.22}
k\le q+F(\nu_4)\, .
   \end{equation} 
For $x>0$, $F(x)$ is an increasing function so that 
$$
F(\nu_4)\le\begin{cases} 
F(\sqrt{q/p})= -\frac{1}{4}\sqrt{pq}+\frac{29}{16}p+\frac{1}{2}+R & 
\text{if $q$ is not a square}\, ,\\
F(\sq/p)=-\frac{1}{4}p\sq+\frac{29}{16}p^2+\frac{1}{2}+R &
\text{otherwise}\, ,
\end{cases}
$$
where
$$
R=\begin{cases} 
-\frac{841p-280}{16(4\sqrt{q/p}+29)} & \text{if $q$ is not a square}\, ,\\
-\frac{841p^2-280}{16(4\sq/p+29)}    & \text{otherwise}\, .
\end{cases}
$$
Then from (\ref{eq5.22}) and since $R<0$ we have the following result.
   \begin{lemma}\label{lemma5.24} Let $\cK$ be a complete arc of size $k$
such that $(3q+5)/4<k\le m'(2,q)$. Suppose that $\cK$ admits an
irreducible envelope such that $\nu_4$ is a power of $p$ and
that $\nu_4<\sq$. Then
$$
k<\begin{cases}
q -\frac{1}{4}\sqrt{pq}+\frac{29}{16}p+\frac{1}{2} &  
\text{if $q$ is not a square}\, ,\\
q-\frac{1}{4}p\sq+\frac{29}{16}p^2+\frac{1}{2} &
\text{otherwise}\, .
\end{cases}
$$
    \end{lemma}
     \begin{corollary}\label{cor5.25} Let $\cK$ be a complete arc of size
$m'(2,q)$. Suppose that $q$ is an even power of $p$ and that $\cK$ admits
an irreducible envelope with $\nu_4$ a power of $p$ and $\nu_4<\sq$.
Then one of the following statements holds:
   \begin{enumerate}
\item[\rm(1)] $p=3$, $\nu_4=\sq/3$, and $m'(2,q)$ satisfies Lemma
\ref{lemma5.24}.
\item[\rm(2)] $p=5$, $q=5^4$, $\nu_4=5$, and $m'(2,5^4)\le 613$;
\item[\rm(3)] $p=5$, $q=5^6$, $\nu_4=5^2$, and $m'(2,5^6)\le 15504$;
\item[\rm(4)] $p=7$, $q=7^4$, $\nu_4=7$, and $m'(2,7^4)\le 2359$.
   \end{enumerate}
     \end{corollary}
     \begin{proof} Let $q=p^{2e}$; so $e\ge 2$ as $p\le \nu_4<p^e$. 
From (\ref{eq5.2}) and Lemma \ref{lemma5.24} we have
that 
$$
(p-4)p^e/4<29p^2/16-0.5\, ,
$$
so that $p\in \{3,5,7,11\}$. 

Let $p=3$. If $\nu_4\le \sq/9$ (so $e\ge 4$), then from (\ref{eq5.2})
and $m'(2,q)\le q+F(\sq/9)$ we would have that 
$$
q-\sq+1\le q-9\sq/4+2357/16-67841/16(43^{e-2}+29)\, ,
$$
which is a contradiction for $e\ge 4$. 

Let $p=11$. Then $p^e\le 125$ and $e=2$ and $\nu_4=11$. Thus from
Proposition \ref{prop5.21} we have $m'(2,11^4)\le 11^4+F(11)$, i.e. 
$m'(2,11^4)\le 14441$. This is a contradiction since by (\ref{eq5.2}) we
must have $m'(2,11^4)\ge 14521$. This eliminates the possibility $p=11$. 

The other cases can be handled in an analogous way.
     \end{proof}
{\bf Case  $\nu_4=\sq$.} In this case, according to (\ref{eq5.21}),
we just obtain Segre's bound (\ref{eq5.1}).

Next we study geometrical properties of irreducible envelopes associated
to large complete arcs in $\P^2(\fq)$, $q$ odd. In doing so we use the
bounds obtained above and divide our study in two cases according as $q$
is a square or not.

{\bf Case $q$ square.} Let $\cX$ be an irreducible envelope
associated to an arc of size $m'(2,q)$. Then from Lemma \ref{lemma5.13}, 
and Corollaries \ref{cor5.22}, \ref{cor5.23}, \ref{cor5.24},
\ref{cor5.25}, we have the following result.
   \begin{proposition}\label{prop5.22} If $q$ is an odd square and
$q>43^2$, then 
$\cX$ is $\Sigma_1$-classical. The $\Sigma_2$-orders 
are $0,1,2,3,4, \epsilon_5$ and the 
$\fq$-Frobenius $\Sigma_2$-orders are $0,1,2,3,\nu_4$, with
$\epsilon_5=\nu_4$, where also one of the following holds:
   \begin{enumerate}
\item[\rm(1)] $\nu_4\in \{\sq/3,3\sq\}$ for $p=3$;
\item[\rm(2)] $(\nu_4,q)\in \{(5,5^4), (5^2,5^6), (7,7^4)\}$;
\item[\rm(3)] $\nu_4=\sq$ for $p\ge 5$.
   \end{enumerate}
   \end{proposition}
{\bf Case $q$ non-square.} In this case there is no analogue to 
bound (\ref{eq5.2}). From Corollary \ref{cor5.21} and Lemmas
\ref{lemma5.22}, \ref{lemma5.23},
\ref{lemma5.24}, and taking into consideration
({\ref{eq5.22}) we have the following result.
   \begin{proposition}\label{prop5.23} Let $q>43^2$ and $q=p^{2e+1}$,
$e\ge 1$. Then, apart from the
values on $\nu_4$, 
the curve $\cX$, $\nu_4$ and $\epsilon_5$ are as in Proposition
\ref{prop5.22}. In this case 
$$
m'(2,q)>q-3\sqrt{pq}/4+7/4
$$
implies
  \begin{enumerate}
\item[\rm(1)] $\nu_4=\sqrt{q/p};$
\item[\rm(2)] $m'(2,q)<q-\sqrt{pq}/4+29p/16+1/2$.
  \end{enumerate}
  \end{proposition}
In particular, our approach just
gives a proof of Segre's bound (\ref{eq5.1}) and Voloch's bound
\cite{voloch2}. However, both
propositions above show the type of curves associated to
large complete arcs. The study of such curves, for $q$ square and large enough,
allowed Hirschfeld and Korchm\'aros
\cite{hk1}, \cite{hk2} to improve Segre's bound (\ref{eq5.1}) to the bound
in (\ref{eq5.3}). 

Next we stress here the main
ideas from \cite{hk2} necessary to deal with Problem \ref{prob5.1}. Due to 
Proposition \ref{prop5.12}, the main strategy is to bound from above the 
number $2M_q+M_q'$ (which is defined via 
(\ref{eq5.11})). For instance, if one could prove that
  \begin{equation}\label{eq5.23}
2M_q+M_q'\le d(q-\sq+1)\, ,
   \end{equation}
where $d$ is the degree of the irreducible envelope whose normalization is
$\cX$, then from Proposition \ref{prop5.12} would follow immediately an 
affirmative answer to Problem \ref{prob5.1}. However, since we know 
the answer to be negative for $q=9$ and $d\le 2t=2(q+2-m'(2,q))$, then
one
can assume that $d$ is bounded by a linear function on $\sq$ and should  
expect to prove (\ref{eq5.23}) only under certain conditions on $q$. 
   \begin{lemma}\label{lemma5.25}
Let $q$ be an odd square. If (\ref{eq5.23}) holds true for $d\le
2\sq-\alpha$ with $\alpha\ge 0$, then $m'(2,q)<q-\sq+2+\alpha/2 $. In
particular, if (\ref{eq5.23}) holds 
true for $d\le 2\sq$, then the answer to Problem \ref{prob5.1} is
positive; i.e,, $m'(2,q)=q -\sq+1$.
    \end{lemma}
    \begin{proof} If $m'(2,q)\ge q-\sq+2+\alpha/2$, then from $d\leq 
2(q+2-m'(2,q))$ we would have that $d\le 2\sq-\alpha$ and so, from
Proposition \ref{prop5.12} and (\ref{eq5.23}), that $m'(2,q)\le q-\sq+1$,
a contradiction.
     \end{proof} 
Now, in \cite{hk1}, (\ref{eq5.23}) is proved for 
$d\le \sq-3$ and $q$ large enough, and so (\ref{eq5.3}) follows. More
precisely we have the following.
    \begin{theorem} {\rm (Hirschfeld-Korchm\'aros \cite[Thm.
1.3]{hk2})}\label{thm5.21} 
Let $q$ be a square, $q>23^2$, $q\neq 3^6$. Let $3\leq d\leq \sq-3$.
Suppose that $\Sigma_1$ is classical, that $0,1,2,3,4,\sq$ are the  
$\Sigma_2$-orders, and that $0,1,2,3,\sq$ are the $\fq$-Frobenius orders
of $\Sigma_2$. Then (\ref{eq5.23}) holds.
    \end{theorem}
    \begin{proof} (Sketch) Suppose that $2M_q+M_q'\geq d(q-\sq+1)$. We are
going to show that $2M_q+M_q'=d(q-\sq+1)$. Notice that $d\geq (\sq+1)/2$
by Corollary \ref{cor3.2}(1). Let
$\phi=(f_0:\ldots:f_5)$ be a morphism
associated to $\Sigma_2$. From Lemma \ref{lemma2.22} there exist
$z_0,\ldots,z_5\in \bar\fq(\cX)$, not all zero, such that
$\sum_{i=0}^{5}z_i^{\sq}f_i=0$. 
Set
    $$
\cZ:=(z_0:\ldots:z_5)(\cX)\, .
    $$
(This curve is related to the dual curve of $\phi(\cX)$ since it 
is easy to see that $\sum_{i=0}^5z_i^{\sq}(P)X_i=0$ is the hyperplane
tangent at $P$ for infinitely many $P$'s.) 

We have \cite[Props. 8.3, 8.4, 8.5]{hk2}
   \begin{enumerate}
\item[\rm(I)] $\sq\deg(\cZ)\leq d(2d+q+3)-(2M_q+M_q')$;
\item[\rm(II)] $\deg(\cZ)\geq \sq j_1(P)$ for any $P\in\cX$;
\item[\rm(III)] $\deg(\cZ)\ge 2\sq$ whenever $\cC_1$ is singular.
   \end{enumerate}
It follows from (I) and (II) that $j_1(P)\leq 2$ since $d\leq \sq-3$. Now
from Corollary 
\ref{cor2.24} and the hypothesis on $d$ there are three possibilities for
$(\Sigma_1,P)$-orders:
   \begin{enumerate}
\item[\rm(A)] $j_2(P)=2j_1(P)$;
\item[\rm(B)] $j_2(P)=(\sq+j_1(P))/2$;
\item[\rm(C)] $j_2(P)=\sq-j_1(P)$.
   \end{enumerate}
We see that points of type (C) cannot occur since $j_1(P)\leq 2$ and
$d\leq \sq-3$. Now from the proof of \cite[Prop. 9.4]{hk2} we have that
   $$
\sq\deg(\cZ)=2(dq+d-2M_q-M_q')\leq 2d\sq, ,
   $$
so that $\deg(\cZ)<2\sq$ as $d\leq \sq-3$. It follows from (III) that
$\cC_1$ is non-singular; i.e., $\cX=\cC_1$. In particular the
$\Sigma_1$-Weierstrass points are of type (B) and we have
   $$
\deg(R_1)=3d(d-2)=(\sq-3)/3\tau\, ,
   $$ 
where $R_1$ is the ramification divisor of $\Sigma_1$ and $\tau$ is the
number of points of type (B). Now we use the
following relation between $\deg(\cZ)$ and $\tau$ \cite[Prop. 9.3]{hk2}:
   \begin{enumerate}
\item[\rm(IV)] $3\deg(\cZ)=2\tau$.
   \end{enumerate} 
Since we already notice that $\deg(\cZ)\leq 2d$ it follows that $d\leq
(\sq+1)/2$; i.e., $d=(\sq+1)/2$. Next we show that $\tau=M_q'$. For $P$ of
type (B), the $(\Sigma_2,P)$-orders are $0,1,2,(\sq+1)/2,(\sq+3)/2,\sq+1$.
Suppose that $P\not\in\cX(\fq)$. Then $2\ell_P$ is the tangent hyperplane
$L_4(P)$ at $P$ with respect to $\Sigma_2$, where $\l_P$ is the tangent
line at $P$ with respect to $\Sigma_1$. It is easy to see that $\fro(P)\in
L_4(P)$ so that $\fro(P)\in \l_P$. This implies $d>(\sq+1)/2$, a
contradiction. Thus $M_q'=3(\sq+1)/2$. Finally by means of
   $$
\deg(S_1)=d(q+d-1)=2M_q+\frac{\sq+1}{2}M_q'\, ,
   $$
where $S_1$ is the $\fq$-Frobenius divisor associated to $\Sigma_1$, we
find that \\$M_q=(\sq+1)(q-\sq-2)/4$, and one easily checks that
$2M_q+M_q'=d(q-\sq+1)$.
    \end{proof}
  \begin{remark}\label{rem5.21} The plane curve $\cX$ of degree
$d=(\sq+1)/2$ in the above proof satisfies 
  $$
\#\cX(\fq)=M_q+M_q'=q+1+\sq (d-1)(d-2)\, ;
  $$
i.e, it is $\fq$-maximal. If $q\geq 121$, such a curve is $\fq$-isomorphic
to the Fermat curve $X^{(\sq+1)/2}+Y^{(\sq+1)/2}+Z^{(\sq+1)/2}=0$; see
\cite{chkt}.
    \end{remark}
Recently, Aguglia and Korchm\'aros \cite{a-k} proved a weaker version of
(\ref{eq5.23}) for $d=\sq-2$ and $q$ large enough, namely
   $$
2M_q+M_q'\leq d(q-\sq/2-9/2)-3\, .
   $$ 
From this inequality and Proposition \ref{prop5.12} one slightly improves
(\ref{eq5.3}) to $m'(2,q)\leq q-\sq/2-11/2$ whenever $d=\sq-2$ and $q$ is
large enough. Therefore the paper \cite{a-k}, as well as \cite{hk1} or
\cite{hk2}, is a good guide toward the proof of (\ref{eq5.23}) for
$\sq-2\le d\le 2\sq$.

\end{document}